\documentclass{article}

\usepackage{hyperref}
\usepackage[left=2.5cm,right=2.5cm,top=2.5cm,bottom=2.5cm]{geometry}
\usepackage{graphicx}
\usepackage{amssymb,amsmath,amscd}
\usepackage{mathtools, color}
\usepackage{cleveref}
\usepackage{multirow}
\usepackage{soul}

\newcommand{\C}{\mathbb{C}}
\newcommand{\mS}{\mathbb{S}}
\newcommand{\mB}{\mathbb{B}}

\newcommand{\N}{\mathbb{N}}

\newcommand{\R}{\mathbb{R}}

\newcommand{\mRed}[1]{{\color{red} \ensuremath{\mathbf{#1}}}}
\newcommand{\aNorm}[1]{\mid\!\mid\!#1\!\mid\!\mid}

\newcommand{\lp}{\left}
\newcommand{\rp}{\right}

\crefformat{equation}{(#2#1#3)}

\title{Saddle Invariant Objects and their Global Manifolds
  in a Neighborhood of a Homoclinic Flip Bifurcation of  Case B}
\author{Andrus Giraldo\footnotemark[2] , Bernd
  Krauskopf\footnotemark[2] , and Hinke M. Osinga\footnotemark[2]}

\begin{document}
\maketitle

\renewcommand{\thefootnote}{\fnsymbol{footnote}}
\footnotetext[2]{Department of Mathematics, The University of
  Auckland, Private Bag 92019, Auckland 1142, New Zealand
  (\href{mailto:a.giraldo@auckland.ac.nz)}{a.giraldo@auckland.ac.nz},
  \href{mailto:b.krauskopf@auckland.ac.nz)}{b.krauskopf@auckland.ac.nz},
  \href{mailto:h.m.osinga@auckland.ac.nz)}{h.m.osinga@auckland.ac.nz})}

\renewcommand{\thefootnote}{\arabic{footnote}}

\begin{abstract}
When a real saddle equilibrium in a three-dimensional vector field
undergoes a homoclinic bifurcation, the associated two-dimensional
invariant manifold of the equilibrium closes on itself in an orientable or
non-orientable way,  provided the corresponding genericity conditions. We are
interested in the interaction between global invariant  
manifolds of saddle equilibria and saddle periodic orbits for a vector
field close to a codimension-two homoclinic flip
bifurcation, that is, the point of transition between having an
orientable or non-orientable two-dimensional surface.  Here, we focus on homoclinic flip bifurcations 
of case $\textbf{B}$, which is characterized by the fact that the
codimension-two point gives rise to an additional homoclinic bifurcation,
namely, a two-homoclinic orbit.  To explain how the global
manifolds organize phase space, we consider  Sandstede's
three-dimensional vector field model, which features inclination and
orbit flip bifurcations. We compute global invariant manifolds and their intersection sets with a
suitable sphere, by means of continuation of suitable two-point
boundary problems, to understand their role as separatrices of basins
of attracting periodic orbits.  We show  
representative images in phase space and on the sphere, such that we
can identify topological properties of the manifolds in
the different regions of parameter space and at the homoclinic
bifurcations involved.  We find heteroclinic orbits between saddle
periodic orbits and equilibria, which give rise to regions of
infinitely many heteroclinic orbits. Additional equilibria exist in
Sandstede's model and we compactify phase space to capture how
equilibria may emerge from or escape to infinity. We present 
images of these bifurcation diagrams, where we outline different
configurations of equilibria close to homoclinic flip bifurcations of
case $\textbf{B}$; furthermore, we characterize  the dynamics of
Sandstede's model at infinity.
\end{abstract}

\begin{section}{Introduction}
Dynamical systems has been an active area of
research since the work of Henri Poincar\'e on
celestial mechanics \cite{Poinc1}.  More recently, the bifurcation theory of
dynamical systems has become a tool for understanding different phenomena,
as far ranging as the excitation of neurons
\cite{Bar1,Hux1, Sneyd1}, turbulence in fluid flows \cite{Tak1,
  Swinney}, and the dynamics of laser systems \cite{Raj1,Wie1};  more applications can
be found, for example, in
\cite{Guck1, Str1}. The models arising in such applications are typically
vector fields of the form 
\begin{equation}
\dot x = f(x,\mu), \label{eq:genP}
\end{equation} 
where $x \in \R^n$ is the state, $\mu \in \R^m$ is a (multi)parameter
and $f: \R^{n} \times \R^{m} \rightarrow \R^{n}$ is a sufficiently 
smooth function. For any fixed value of $\mu$, equation \cref{eq:genP} defines a flow
$\phi^t$ on the phase space $\R^n$ for all $t \in
\R$. In bifurcation theory, one wants to
understand how the phase portrait of this flow $\phi^t$ changes topologically when
$\mu$ is varied.  One way for
such topological changes to occur is through changes of stabilities
of equilibria and periodic orbits in phase space. These
are known in the literature as \emph{local bifurcations}; and they have been
studied in detail by normal forms and desingularization
techniques \cite{Guck1,Kuz1,perk,Str1}. \emph{Global bifurcations}, on the other
hand, are topological changes arising from
interactions between global invariant manifolds of saddle equilibria and  saddle
periodic orbits, which can re-arrange to  change the phase space globally.
In particular, the existence of  homoclinic or
heteroclinic orbits, which are connecting orbits between saddle
equilibria and/or saddle periodic orbits, can have dramatic effects,
for example, regarding the existence and size of basins of attractions. 

We are interested in gaining a better understanding of a special type
of global bifurcation that is know as \emph{homoclinic flip
  bifurcation}; it can occur in vector fields of dimension three or higher.  This
bifurcation concerns a real saddle equilibrium and a homoclinic orbit, that is, a connecting orbit of the
equilibrium back to itself, such that the associated invariant
manifolds are neither orientable or non-orientable; see 
\cref{sec:hFB} for details.  To study this bifurcation, we work
with the three-dimensional vector field
\begin{equation}
X^s(x,y,z):
\begin{cases} 
\dot x = P^1(x,y,z) := ax + by -ax^2+(\tilde{\mu}-\alpha z)x(2-3x), \\
\dot y = P^2(x,y,z) := bx +ay -\frac32 bx^2-\frac32 axy-2y(\tilde{\mu}-\alpha
z), \\
\dot z = P^3(x,y,z) := cz +\mu x +\gamma xz +\alpha \beta  (x^2(1-x)-y^2).
\end{cases}
\label{eq:san}
\end{equation}
It was introduced by Sandstede in \cite{san1}, who studied
this model with additional $z$-dependent terms in the equations for $x$
and $y$, which are controlled by a parameter $\delta$ in \cite{san1} that we set to
$0$ in \cref{eq:san}.  We choose the parameters such that  the origin $\mathbf{0} \in
\R^3$ is a saddle equilibrium of system \cref{eq:san} whose
linearization has  two different negative and one
positive eigenvalues  $\lambda^{ss}<\lambda^s<0<\lambda^u$; see
\cref{sec:san}. The other
case of a real saddle, that is, an equilibrium with two positive
and one negative eigenvalues can be reduced to this case by reversing time. Since
the origin is hyperbolic, the Stable Manifold Theorem \cite{Taken1} implies the
existence of an immersed two-dimensional stable manifold $W^s(\textbf{0})$ and
an immersed one-dimensional unstable manifold $W^u(\textbf{0})$; the stable manifold
$W^s(\textbf{0})$ is a surface foliated by orbits that converge to
$\textbf{0}$ as $t \rightarrow \infty$, and $W^u(\textbf{0})$
consist of two orbits that converge to $\textbf{0}$
as $t \rightarrow -\infty$.  
\begin{figure}
\centering
\includegraphics{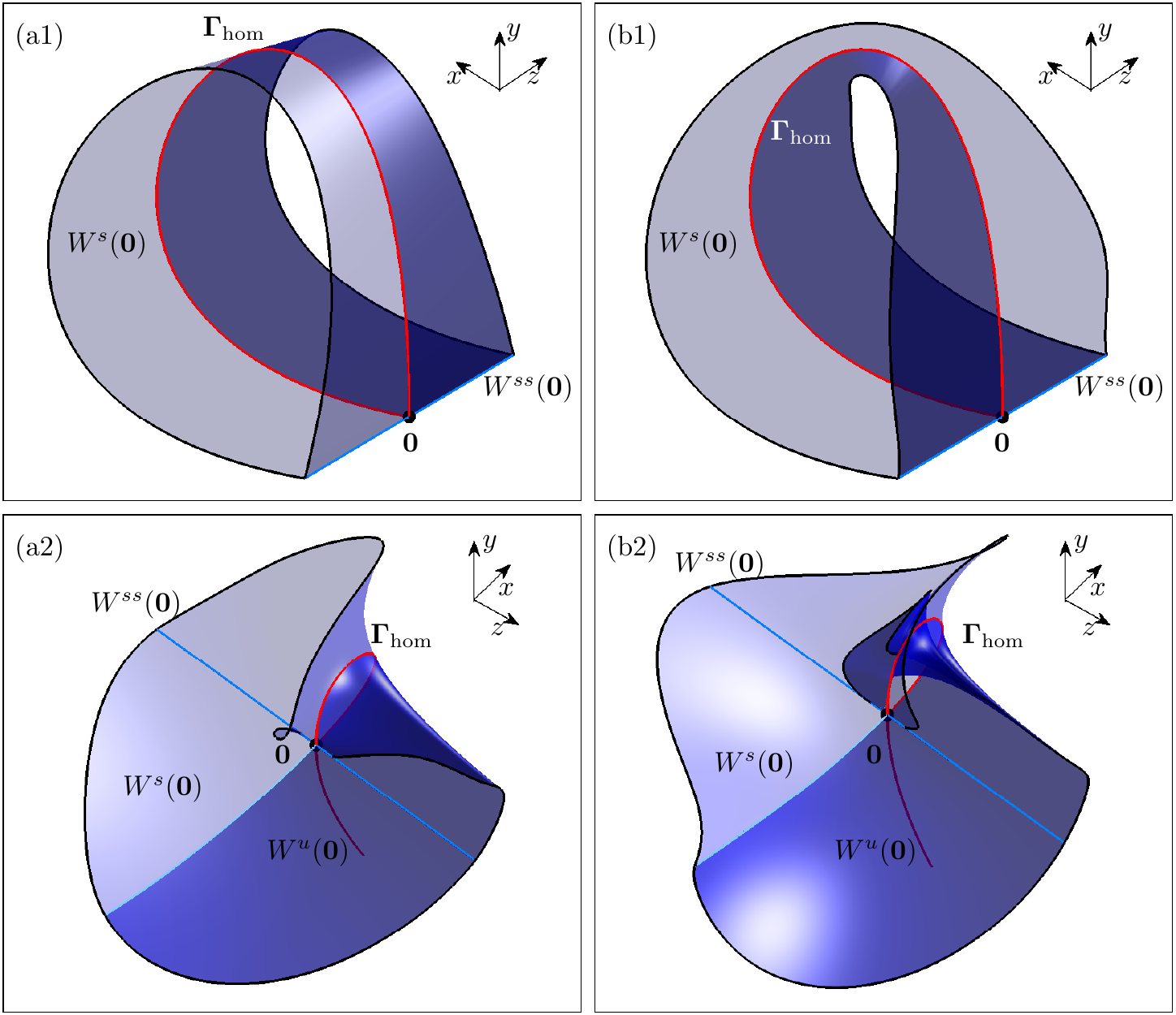}
\caption{The stable manifold $W^{s}(\mathbf{0})$ of system
  \cref{eq:san} at a codimension-one 
  homoclinic bifurcation in $\R^3$. Column (a) shows an orientable and 
  column (b) a non-orientable homoclinic orbit.  The top row
  illustrates the tangent space of $W^{s}(\mathbf{0})$ around  $\mathbf{\Gamma_{\rm hom}}$. The bottom row shows a portion of 
  $W^{s}(\mathbf{0})$ in phase space. Shown are $W^{s}(\mathbf{0})$
  as a rendered surface  with one half colored 
  dark-blue and the other half light-blue,  the one-dimensional strong stable manifold
  $W^{ss}(\mathbf{0})$ as a light blue curve and the unstable manifold
  $W^u(\mathbf{0})$ as  a red curve. Column (a) 
  is for $(a,b,c,\alpha,\beta,\gamma,\mu,\tilde{\mu})=  
  (0.22,1,-2,0.3,1,2,0,0)$ and column (b)  for
  $(a,b,c,\alpha,\beta,\gamma,\mu,\tilde{\mu})= 
  (0.22,1,-2,0.65,1,2,0,0)$. } \label{fig:ori}
\end{figure}  

We can choose $\alpha, \mu$ and $\tilde{\mu}$ such that
a homoclinic bifurcation occurs: this means that there exists an
orbit $\mathbf{\Gamma_{\rm hom}}$ that converges to $\textbf{0}$  both
as $t \rightarrow \infty$ and as $t 
\rightarrow -\infty$;  that is, one of the branches of
$W^u(\textbf{0})$ lies entirely in $W^s(\textbf{0})$. Consequently, under
certain genericity conditions that are outlined in \cref{sec:hFB}, the respective local part of $W^s(\textbf{0})$ 
closes back on itself and is either topologically equivalent to a  
cylinder or a M\"obius band, which classifies the homoclinic
bifurcation as orientable or non-orientable, respectively
\cite{Deng1,Hom1,kis1,Shil1}. \cref{fig:ori} shows two homoclinic 
orbits of system \cref{eq:san}, together with the associated stable and
unstable manifolds of $\textbf{0}$;  they share the same
parameter values except that $\alpha=0.3$ in column (a) and
$\alpha=0.65$ in column (b). The top row shows a linear approximation
of $W^s(\mathbf{0})$ around the homoclinic orbit $\mathbf{\Gamma_{\rm hom}}$, which
its computed as the span of the tangent vectors of $W^s(\mathbf{0})$
around $\mathbf{\Gamma_{\rm hom}}$. The bottom row shows
$W^s(\mathbf{0})$ in a larger region of phase space. The stable manifold
$W^s(\mathbf{0})$ is rendered in two shades of blue to illustrate the
orientability properties of the manifold. In
\cref{fig:ori} panel (a1) we see that $W^s(\textbf{0})$, locally near
$\mathbf{\Gamma_{\rm hom}}$, is a topological cylinder; while in panel (b1), this local part of
$W^s(\mathbf{0})$ is a topological M\"obius band. Indeed \cref{fig:ori}
illustrates how $W^s(\mathbf{0})$ closes on
itself along the strong stable manifold $W^{ss}(\mathbf{0})$ at the
moment of an orientable homoclinic bifurcation in column (a); and the
non-orientable case in column (b).

Generically, homoclinic orbits as shown in \cref{fig:ori} exist at codimension-one bifurcations
\cite{Kuz1, Shil2}.  We are interested in the case when one of the genericity conditions is not
valid, such that the homoclinic bifurcation has
codimension-two. More precisely, we study the case where
$W^s(\mathbf{0})$ transitions from being orientable to being
non-orientable. This codimension-two point is called a \emph{homoclinic
  flip bifurcation} and it may be an \emph{inclination} or an
\emph{orbit flip bifurcation} \cite{san3}.  There are three different
codimension-two unfoldings of homoclinic flip bifurcation, called cases \textbf{A, B} and
\textbf{C}, for both an inclination and an orbit flip
bifurcation; these have been studied theoretically with methods including return maps \cite{Deng1, Hom1},  
Shilnikov variables \cite{kis1} and Lin's method \cite{san4}.  The
theoretical results describe the unfoldings of the dynamics locally in
a small tubular neighbourhood of the homoclinic orbit. A ``more
global'' approach, which relies on numerical computations, 
has been used in \cite{Agu1}  to understand how the global manifolds
re-arrange phase space for the simplest case \textbf{A}. Already for
this case an extra bifurcating branch of heteroclinic folds was found that had 
previously not been identified.
\begin{figure}
\centering
\includegraphics{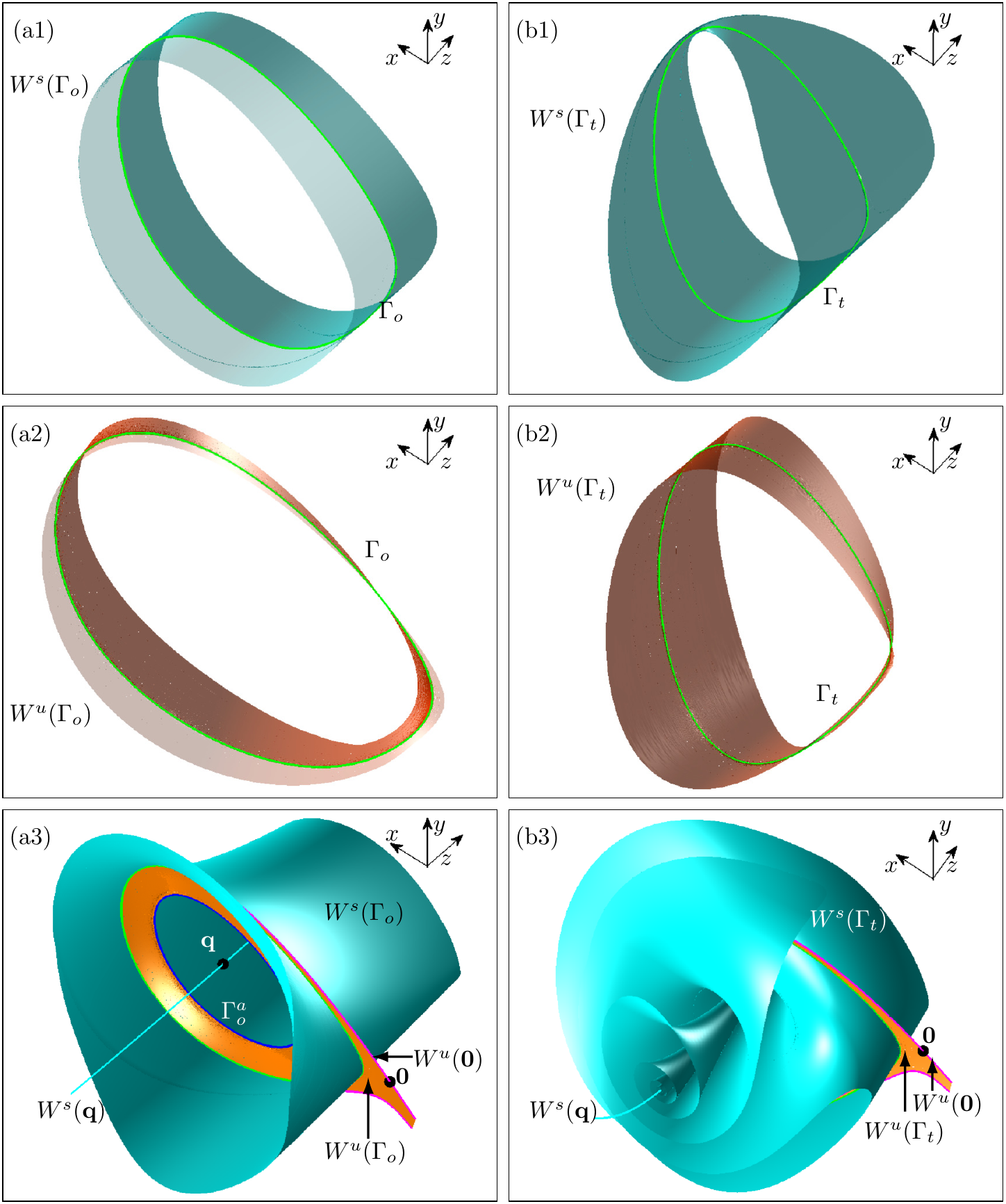}
\caption{The stable manifolds $W^s(\Gamma_{o})$ and $W^s(\Gamma_{t})$,
  and unstable manifolds $W^u(\Gamma_{o})$ and $W^u(\Gamma_{t})$ of periodic orbits
  $\Gamma_{o}$ and $\Gamma_{t}$ in  system~\cref{eq:san}.  Column (a)  shows the orientable and
  column (b) the non-orientable case. The top row shows 
  $W^s(\Gamma_{o})$ and  $W^s(\Gamma_{t})$ in a small tubular
  neighborhood of $\Gamma_{o}$ and $\Gamma_{t}$, respectively, while the middle row 
  shows the  $W^u(\Gamma_{o})$ and $W^u(\Gamma_{t})$ in this same small tubular
  neighborhood. The bottom  row, shows a large portion of these manifolds and how they are
  interacting with other one-dimensional manifolds in phase space. Shown are
  $W^s(\Gamma_{o})$ and $W^s(\Gamma_{t})$  as cyan surfaces,
  $W^u(\Gamma_{o})$ and $W^u(\Gamma_{t})$   as orange surfaces,
  $\Gamma_{o}$ and $\Gamma_{t}$ as green curves,
  $\Gamma^a_o$ as a blue curve, $W^u(\mathbf{0})$ as a pink
  curve, and $W^s(\mathbf{q})$ as a cyan curve. Column (a)  is
  for 
  $(a,b,c,\alpha,\beta,\gamma,\mu,\tilde{\mu})=(0.22,1,-2,0.3,1,2,0.004,0)$, 
  and  column (b) is for
  $(a,b,c,\alpha,\beta,\gamma,\mu,\tilde{\mu})=
  (0.22,1,-2,0.65,1,2,-0.004,0)$.}\label{fig:per}
\end{figure}

Compared with case \textbf{A}, cases \textbf{B} and \textbf{C} are
richer with respect to the invariant objects that are created and destroyed close to
the homoclinic flip bifurcation. In case \textbf{A} a single
attracting (or repelling) periodic orbit is created. The
unfolding of case \textbf{B}, on the other hand, involves saddle periodic orbits, a period doubling and an additional
homoclinic bifurcation curve; see \cref{sec:Inc}
and \cref{fig:BDInc}. Finally, in the unfolding of case \textbf{C} it
has been proved that there exists a period-doubling
cascade, region of horseshoe dynamics, $n$-homoclinic orbits (for any $n \in \N$) and
strange attractors \cite{Deng1,Hom2,san3,Nau1,Nau2}.
Additionally, case \textbf{C} has been identified as an organizing
center for the creation of spikes of periodic orbits in the
Hindmarsh-Rose model that describes the essential 
spiking behaviour of a neuron \cite{Lina1}.  

We focus our attention on case \textbf{B}. It is the next  step in
understanding a more complicated case, namely, case \textbf{C},
and also the main ingredient in the homoclinic-doubling cascade that
appears close to bifurcations of higher codimension \cite{Hom4}. As
mentioned before, the unfolding of case \textbf{B} involves saddle
periodic orbits; these have two-dimensional stable and unstable manifolds that may  or may not be
orientable; this is illustrated in \cref{fig:per} for parameter values close to
case \textbf{B} for system~\cref{eq:san}.  The first column shows an 
orientable saddle periodic orbit $\Gamma_o$, and the second column a
non-orientable (twisted) periodic orbit $\Gamma_t$.  \cref{fig:per}
shows portions of the two-dimensional stable manifolds $W^s(\Gamma_{o})$ and
$W^s(\Gamma_{t})$, and unstable manifolds $W^u(\Gamma_{o})$ and
$W^u(\Gamma_{t})$ of the saddle periodic orbits  $\Gamma_o$ and
$\Gamma_t$, respectively. The first row shows $W^s(\Gamma_{o})$ and
$W^s(\Gamma_{t})$ as contained in a small tubular neighborhood with 
radius $d=0.01$ around $\Gamma_{o}$ and $\Gamma_{t}$. The second row shows
$W^u(\Gamma_{o})$ and $W^u(\Gamma_{t})$  in the same corresponding
tubular neighborhood. The third 
row illustrates a larger portion of these manifolds in phase
space. Panels (a1) and (a2) illustrate that both $W^s(\Gamma_o)$ and
$W^u(\Gamma_o)$ are homeomorphic to a cylinder. In panel (a3), we see that
$W^s(\Gamma_o)$ is unbounded but remains a topological cylinder; 
one side of $W^u(\Gamma_o)$ accumulates on an attracting periodic
orbit denoted $\Gamma^a_o$ and the other side is bounded by $W^u(\mathbf{0})$.
Since $W^s(\Gamma_o)$ is unbounded and orientable, it 
acts as a separatrix that bounds the basin of attraction of
$\Gamma^a_o$. In contrast, \cref{fig:per} (b1)  and (b2) illustrate
that $W^s(\Gamma_t)$ and $W^u(\Gamma_t)$ are a topological
M\"obius band. The non-orientable nature of $W^s(\Gamma_t)$
and $W^u(\Gamma_t)$ is hard to appreciate in panel (b3), but we can
see that $W^s(\Gamma_t)$  spirals around the one-dimensional stable manifold
$W^s(\mathbf{q})$ of an additional equilibrium denoted $\mathbf{q}$, and the unstable manifold
$W^u(\Gamma_t)$ is bounded by $W^u(\mathbf{0})$.  

One of the biggest advantage of using numerical techniques to understand the
behaviour of these manifolds is the possibility to study their
interactions with other invariant objects and to
determine how they can organize  phase space; this is shown in
\cref{fig:manT}  for the orientable case in panel (a) and the non-orientable  case in panel
(b). Here, we also show the two-dimensional stable manifold $W^s(\mathbf{0})$ and the 
two-dimensional unstable manifold  $W^u(\mathbf{q})$ of
$\mathbf{q}$. \cref{fig:manT}(a) illustrates how  $W^s(\mathbf{0})$
spirals towards the topological cylinder formed by 
$W^s(\Gamma_o)$.  Moreover, $W^s(\Gamma_o)$ does not interact with
$W^u(\mathbf{q})$, which accumulates onto $\Gamma^a_o$; indeed,
$W^u(\mathbf{q})$ lies in the basin of attraction 
of $\Gamma^a_o$ and $W^s(\mathbf{0})$ does not. The non-orientable
case in \cref{fig:manT}(b) is quite different, $\Gamma_o^a$ does not exist and $W^s(\mathbf{0})$ 
together with $W^s(\Gamma_o)$ rolls around
$W^s(\mathbf{q})$.  As a consequence, $W^s(\Gamma_t)$ intersects 
$W^u(\mathbf{q})$ transversally, implying the existence of a
heteroclinic orbit from $\mathbf{q}$ to $\Gamma_t$. 

The main purpose of this paper is to understand how the different 
manifolds of periodic orbits and equilibria organize the phase space
and basins of attraction close to a homoclinic flip bifurcation of
case \textbf{B}. For this reason, we choose parameter values in each
open region of parameter plane, close to the homoclinic flip
bifurcation point, to provide representative figures of
phase space; here we render each invariant object as in
\cref{fig:manT} and analize their transition as a set of parameters is
varied.  For the purpose of understanding  the nature of the basins of attracting periodic orbits, we
also provided figures of the intersection sets of the 
stable manifolds with a suitable sphere.  This allows us to describe
such basins when parameters are varied.  As in case \textbf{A}, 
the existence of an additional saddle focus equilibrium \textbf{q} in Sandstede's model
creates additional dynamics in phase space; these include the
existence of a fold curve of (structurally stable) heteroclinic orbits
from $\mathbf{q}$ to $\mathbf{0}$ in parameter plane. New for case
\textbf{B} is that the existence of \textbf{q} creates regions where there
are infinitely many heteroclinic orbits in phase space; these are
consequence of structurally stable heteroclinic orbits from $\mathbf{q}$
to $\Gamma_t$. We investigate the role and bifurcations of the additional equilibrium \textbf{q} in
system~\cref{eq:san}. Moreover, we analyze all equilibria in Sandstede's model
and consider their bifurcations for parameters  close to the homoclinic
bifurcation.  We find that some equilibria
disappear at infinity.  Therefore, we utilize \emph{Poin\-ca\-r\'e
  com\-pac\-ti\-fi\-ca\-tion} \cite{Dum1,Vel, Jes1} to characterize Sandstede's
model at infinity  and to complete the bifurcation diagram of  these
equilibria; see \cref{secApp:Com} for details. 
\begin{figure}
\centering
\includegraphics[width=14cm,height=18.4cm]{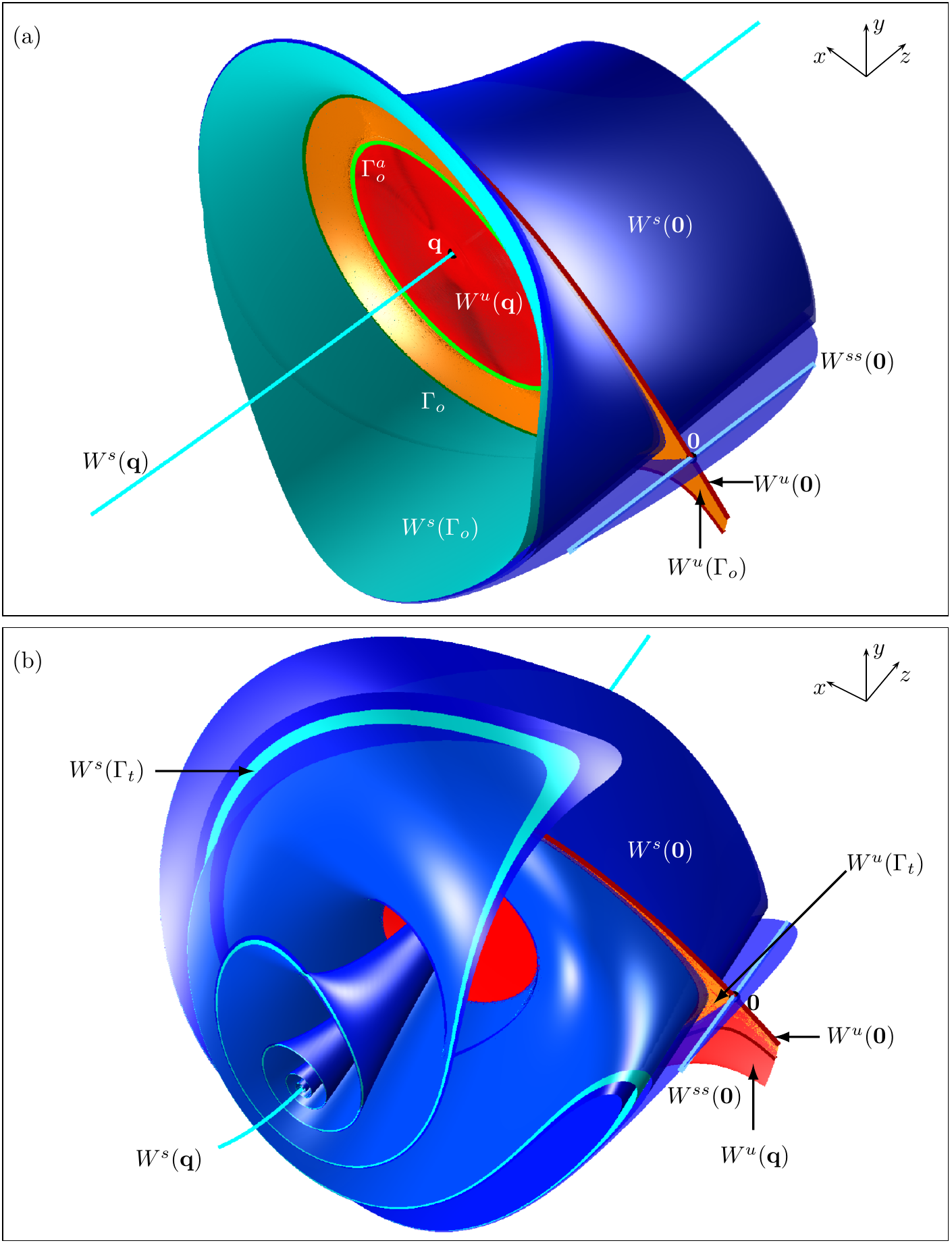}
\caption{Interaction of the different manifolds of system~\cref{eq:san} in
  $\R^3$. Shown are $W^s(\mathbf{0})$ as a dark-blue surface,
  $W^{ss}(\mathbf{0})$ as a blue curve, $W^{u}(\mathbf{0})$ as 
  a pink curve, $W^{u}(\mathbf{q})$ as a red surface,
  $W^{s}(\mathbf{q})$ as a cyan curve,  $W^s(\Gamma_o)$ and
  $W^s(\Gamma_t)$ as cyan surfaces, and $W^u(\Gamma_o)$ and
  $W^u(\Gamma_t)$  as orange surfaces. Panel (a) is for 
  $(a,b,c,\alpha,\beta,\gamma,\mu,\tilde{\mu})=(0.22,1,-2,3,1,2,0.004,0)$,
  and panel (b) for
  $(a,b,c,\alpha,\beta,\gamma,\mu,\tilde{\mu})=
  (0.22,1,-2,0.65,1,2,-0.004,0)$; compare with \cref{fig:per}.} \label{fig:manT}
\end{figure}

The computations in this paper are performed with the software package
\textsc{Auto} \cite{Doe1,Doe2} and its 
extension \textsc{HomCont} \cite{san2}.  In particular, the global
manifolds are computed with a  two-point boundary value problem (2PBVP)
set-up; see \cite{Call1, Kra2} for details. 

The organization of this paper is as follows. In \cref{sec:Not}
we introduce notation and background material; here, we also present
the parameter values that we use to unfold the two types of homoclinic
bifurcations in system~\cref{eq:san}.
In \cref{sec:Dyc_Inf}, we study the bifurcation diagram of the
equilibria via a compactified version of system~\cref{eq:san}.  The
codimension-two inclination and orbit 
flip bifurcations for case \textbf{B} are the subjects of \cref{sec:Inc}  and \cref{sec:Or}, respectively.  \cref{sec:Dis}
contains the discussion of the results and an outlook for future
research.  \cref{secApp:Com} give a brief summary of  Poin\-ca\-r\'e
  com\-pac\-ti\-fi\-ca\-tion and an analytic study of
  system~\cref{eq:san} at infinity. Finally, \cref{secApp:BVP}
  introduces  the 2PBVP-formulations for the computation of a two-dimensional
manifold inside a tubular section, and for curves
along which the Floquet multipliers of a periodic
orbit become complex conjugates.
\end{section}


\begin{section}{Notation and set-up} \label{sec:Not}
Recall that we consider system~\cref{eq:san} as a representative
example of a  three-dimensional  vector field of the form
\cref{eq:genP} with a hyperbolic real saddle 
equilibrium at $\mathbf{0} \in \R^3$.  We choose parameters such
that the Jacobian $Df(\mathbf{0})$ of $\textbf{0}$ has two stable
and one unstable eigenvalues, $\lambda^{ss} < \lambda^{s} < 0 < \lambda^{u}$;
we denote by $e^{ss}$, $e^{s}$ and $e^{u}$ the respective associated
eigenvectors. The global stable and unstable manifolds of $\textbf{0}$ are defined as
\begin{align*}
      W^s(\bf{0}) &:= \{x \in \R^3 : \phi^{t}(x) \rightarrow \textbf{0} \text{ as } t  \rightarrow
\infty \}, \text{ and}\\
      W^u(\bf{0}) &:= \{x \in \R^3 : \phi^{t}(x) \rightarrow \textbf{0}  \text{ as } t  \rightarrow
-\infty \}. 
\end{align*} 
The Stable Manifold Theorem \cite{Taken1} guarantees that both
$W^s(\mathbf{0})$ and $W^u(\mathbf{0})$ are immersed manifolds that
are as smooth as $f$ and tangent at $\textbf{0}$ to the linear
eigenspaces $E^s(\textbf{0})=\text{span} \{e^s,e^{ss}\}$, and 
$E^u(\textbf{0})=\text{span}\{e^u\}$, respectively.  Furthermore, $W^{s}(\mathbf{0})$ has
a one-dimensional strong stable manifold $W^{ss}(\mathbf{0})$, defined as the
subset of points on $W^{s}(\mathbf{0})$  that converge to $\mathbf{0}$ tangentially
to $e^{ss}$.

The stability and invariant manifolds of a periodic orbit $\Gamma$ of
system~\cref{eq:san} are defined in a very similar way.  We denote
its two nontrivial Floquet multipliers by 
$\Lambda_1,\Lambda_2 \in \C$; they are the eigenvalues of the
variational equation along $\Gamma$ over the period of $\Gamma$. Note
that there is also the trivial Floquet multiplier $1$ associated with the
tangent direction of $\Gamma$. In a 
three-dimensional vector field $\Lambda_1$ and $\Lambda_2$ are  always such that their real parts
have the same sign; moreover each
has an associated eigenfunction that is referred as the Floquet bundle
\cite{Leonid1}. 

If $\Lambda_1,\Lambda_2\in\R$ and $0<\Lambda_1<1<\Lambda_2$ then
one speaks of an orientable saddle periodic orbit,  which we denote by
$\Gamma_o$. Its stable and unstable
manifolds $W^s(\Gamma_o)$ and $W^u(\Gamma_o)$, respectively, are
locally a cylinder \cite{Hin1, Leonid1};  if  $\Lambda_2<-1<\Lambda_1<0$, then the
saddle periodic orbit is non-orientable, denoted 
$\Gamma_t$, and $W^s(\Gamma_t)$ and $W^u(\Gamma_t)$ are locally
a M\"obius band \cite{Hin1, Leonid1}. The associated stable and
unstable manifolds of a saddle periodic orbit are two-dimensional
immersed manifolds that are tangent to the Floquet bundle of the periodic orbit associated
with $\Lambda_1$ and $\Lambda_2$, respectively. 

On the other hand, if $\Lambda_1,\Lambda_2 \in \C$ such that $|\Lambda_i|<1$ for
$i=1,2$ then we speak of an attracting periodic orbit, which we denoted
by $\Gamma^a$. This implies the existence of an open set $U \subset \R^3$ that satisfies 
\begin{equation}
\forall t \geq 0, \; \phi^t(U) \subset U \text{ and }
\bigcap_{t>0}\phi^t(U) = \Gamma^{a}. \label{eq:basC}
\end{equation}
Furthermore, the basin of attraction
$\mathcal{B}(\Gamma^{a})$ of $\Gamma^{a}$ is defined as 
the set of all points in phase space that converge to $\Gamma^{a}$, 
that is, $\mathcal{B}(\Gamma^{a})=\bigcup_{t \leq 0}
\phi^t(U)$.  If $|\Lambda_1|<|\Lambda_2|<1$ are both real then we define the strong stable manifold $W^{ss}(\Gamma^{a})$ of $\Gamma^a$
as the set of points that converge to $\Gamma^{a}$ tangent to the
Floquet bundle associated with $\Lambda_1$; this strong stable manifold $W^{ss}(\Gamma^{a})$
is a two-dimensional immersed manifold. Using the same
terminology and notation, we denote the periodic orbit $\Gamma^a_o$
and its strong stable manifold $W^{ss}(\Gamma^a_o)$ if
$0<\Lambda_1<\Lambda_2<1$ and they are orientable; and
as $\Gamma^a_t$ and $W^{ss}(\Gamma^a_t)$ if
$-1<\Lambda_2<\Lambda_1<0$ and they are non-orientable.

\subsection{Homoclinic Flip Bifurcations} \label{sec:hFB}
Let $\mathbf{\Gamma_{\rm hom}}$ be a homoclinic orbit of $\textbf{0}$, that 
is, $\mathbf{\Gamma_{\rm hom}} \subset W^s(\mathbf{0}) \cap
W^u(\mathbf{0})\neq\emptyset$ converges in forward and backward time to
$\mathbf{0}$. The homoclinic orbit 
$\mathbf{\Gamma_{\rm hom}}$ is of codimension one,
provided the following conditions hold \cite{Hom1,kis1}.
\begin{itemize}
\item[(\textbf{G1})] (Non-resonance) $|\lambda^s| \not = \lambda^u$;
\item[(\textbf{G2})] (Principal homoclinic orbit) In positive time the homoclinic
  trajectory approaches the origin tangent to the weakest stable direction $e^s$;
\item[(\textbf{G3})] (Strong inclination) The tangent space
  $\rm{T}W^s(\mathbf{0})$ of the 
  stable manifold, followed along $\mathbf{\Gamma_{\rm hom}}$
  backward in time, converges to 
  span$\{e^{ss},e^u\}$.
\end{itemize}
For any codimension-one homoclinic orbit, a portion of
$W^s(\mathbf{0})$ folds over and closes up along
$W^{ss}(\mathbf{0})$; hence, the immersion of $W^s(\mathbf{0})$ in the three-dimensional
phase space becomes orientable or non-orientable close to
$\mathbf{\Gamma_{\rm hom}}$; see panels (a1) and (b1) in \cref{fig:ori}.

If precisely one of the genericity conditions is not
fulfilled then the homoclinic orbit is of
codimension two, leading to different kinds of unfoldings. If
\textbf{(G1)} fails, then one speaks of a \emph{resonant homoclinic
  bifurcation} \cite{Hom2}.  We focus on the inclination flip
\textbf{(IF)} and orbit flip  bifurcations \textbf{(OF)}  
that occur when conditions \textbf{(G2)} or \textbf{(G3)} fail,
respectively. In both cases, one speak of a homoclinic flip
bifurcation, which is of codimension two, if additional genericity
conditions are satisfied. Then there exists a curve of homoclinic
orbits in any suitable two-parameter plane along which $W^s(\mathbf{0})$ changes from orientable
to non-orientable at the codimension-two flip bifurcation point \cite{Hom1,kis1}.

The unfolding of a flip bifurcation depends on  the eigenvalues of
$\textbf{0}$.  Three cases have been identified for the 
inclination flip and the orbit flip bifurcations; they are denoted \textbf{A},
\textbf{B} and \textbf{C}.  The unfolding of these cases are
topologically the same for both \textbf{IF} and \textbf{OF}  
but they satisfy different conditions \cite{san3}.  More specifically, the eigenvalue conditions and the
unfoldings for the respective cases  are:
\begin{itemize}
\item[\textbf{A.}] If $|\lambda^{s}|>\lambda^u$ then a single
  attracting periodic orbit $\Gamma^{a}$ is created, for both the orbit flip and
  inclination flip bifurcations.
\item[\textbf{B.}] Suppose the following respective conditions for the inclination
  and orbit flip bifurcations are satisfied:
\begin{tabbing}
\textbf{(IF)} \= \hspace{0.1cm} \= $\lambda^u/2<|\lambda^{s}|<\lambda^u$ and
  $|\lambda^{ss}|>\lambda^u$, or \\
\textbf{(OF)} \> \> $|\lambda^{s}|<\lambda^u$ and
  $|\lambda^{ss}|>\lambda^u$.
\end{tabbing}
Then the unfolding contains a homoclinic doubling
bifurcation, a period-doubling bifurcation and saddle-node bifurcation of periodic orbits.
\item[\textbf{C.}] Suppose the following respective conditions for the inclination
  and orbit flip bifurcation are satisfied:
\begin{tabbing}
\textbf{(IF)} \= \hspace{0.1cm} \=  $|\lambda^{s}|<\lambda^u$ and
  $|\lambda^{ss}|<\lambda^u$, or $|\lambda^{s}|<\lambda^u/2$ and
  $|\lambda^{ss}|>\lambda^u$, or \\
\textbf{(OF)} \> \>  $|\lambda^{s}|<\lambda^u$ and $|\lambda^{ss}|<\lambda^u$.
\end{tabbing}
Then the unfolding contains $k$-homoclinic bifurcations \cite{Hom1},
for any $k \in \N$, and a region with horseshoe 
dynamics exists. Two different bifurcation diagrams can arise
depending on extra genericity conditions regarding the geometry of the
stable manifold $W^s(\mathbf{0})$; details can be found \cite{san3}.  
\end{itemize}
For both \textbf{IF} and \textbf{OF}, the unfolding and eigenvalue conditions
for cases \textbf{A} and \textbf{B}  were proven for any smooth vector field of dimension $ n \geq 3$; see
\cite{kis1, san4}. On the other hand, for case \textbf{C}  it has been
proved that regions of horseshoe dynamics, cascades of period-doubling and homoclinic
bifurcations, and strange attractors exist \cite{Deng1,Hom1,Hom4,kis1,Nau1,Nau2}; however, these results have only
been proved for three-dimensional vector fields. Moreover, our understanding
of the exact nature of the unfoldings of case \textbf{C} it is not as
complete as cases \textbf{A} and \textbf{B}. 

Since the eigenvalues of an equilibrium depend continuously on the parameters for smooth
vector fields, the transitions between cases \textbf{A}, \textbf{B} and
\textbf{C}, for both inclination and orbit flip bifurcations,
are codimension-three phenomena caused by resonance, that is, a
violation of condition \textbf{(G1)}; such a resonant
  homoclinic flip bifurcation  was studied in \cite{Hom2} and explored numerically in
\cite{OldKra1}. Furthermore, in \cite{Mor1} it was shown that
$C^1$-near a vector field exhibiting an orbit flip of case
\textbf{C}, there is a vector field with an inclination flip; this
approximation result for $C^1$-flows provides further 
insight into the similarities between these  two types of flip 
bifurcations. 

We remark that the conditions for the homoclinic flip bifurcation
of a hyperbolic equilibrium have been studied for the
non-hyperbolic case, namely, for the case of a transcritical bifurcation \cite{Liu1};
the authors show that the non-hyperbolic
equilibrium gives rise to new heteroclinic orbits and that its unfolding is
different from the hyperbolic case.  Reference \cite{Hom3} explores the creation of a
\emph{Lorenz-like attractor} in homoclinic loop configurations that
exhibit homoclinic flip bifurcations;  this happens when
two homoclinic orbits connect to the same equilibrium, which
in \cite{Hom3} is studied by looking at systems with reflectional symmetry.

\subsection{Sandstede's Model} \label{sec:san}
Sandstede \cite{san1} introduced a model vector field that
exhibits codimension-two flip bifurcations and is particularly
suitable for studying their unfoldings. The system has been constructed in such a way that the
different cases of both inclination and orbit flip bifurcations occur,
and the homoclinic orbits always involve the  
equilibrium $\textbf{0}$. 

As mentioned in the introduction, we set $\delta=0$ in the original
model introduced in \cite{san1}, such that the
$z$-axis is invariant. That is, we work with the vector field
$X^s(x,y,z)$ as defined by system~\cref{eq:san}. Note that
$\textbf{0}$ is an equilibrium of  $X^s$ for all 
parameter values; furthermore, since $\delta=0$ (or if $\mu=0$), the eigenvalues of
$\textbf{0}$ are given by
\begin{equation*}
\lambda_{1,2}= a \pm \sqrt{b^2+4\tilde{ \mu}^2} \text{ and }
\lambda_3=c,
\end{equation*}
and the eigenvector associated with $\lambda_3$ points in the
$z$-direction.  Taking into account the parameter ranges found in \cite{OldKra1} for 
\textbf{IF} and \textbf{OF} of case \textbf{B}, we choose the
following values for the other parameters:
\begin{itemize}
\item[\bf{(IF)}] For $(a,b,c,\beta,\gamma)=
  (0.22,1,-2,1,2)$, there is an inclination flip at
  $(\alpha,\mu,\tilde{\mu})= (\alpha_B,0,0)$ where $\alpha_B \approx
  0.4664012$. At this point $\lambda_1=1.22$, $\lambda_2=-0.78$ 
  and $\lambda_3=-2$. Note that $\lambda^{ss}=\lambda_3$, so for this
  and nearby parameter values $W^{ss}(\mathbf{0})$ is the
  $z$-axis. The codimension-two bifurcation is unfolded by 
  $\alpha$ and $\mu$. 
\end{itemize}

\begin{itemize}
\item[\bf{(OF)}] For $(a,b,c,\beta,\gamma)=
  (-0.5,2.5,-1,0,0)$ , there is an orbit flip at
  $(\alpha,\mu,\tilde{\mu})=(1,0,0)$. At this point $\lambda_1=2$, $\lambda_2=-3$
  and $\lambda_3=-1$. The codimension-two bifurcation is unfolded by
  $\mu$ and $\tilde{\mu}$. 
\end{itemize}
Note that the choice $\gamma=2$ for the case \textbf{IF} differs from
the values taken in \cite{Agu1,OldKra1}.  In \cite{Agu1}, the value $\gamma=0$ was used, but for
case \textbf{B} this value does not give an inclination flip.  In \cite{OldKra1}, the value $\gamma=3$ was
used, but it turns out that the choice $\gamma=3$ is rather
unfortunate with respect to the existence of additional equilibria; we
justify our choice of $\gamma=2$ in the next subsection.
\end{section}

\subsubsection{Configuration of Equilibria} \label{sec:Dyc_Inf} 
We start by determining the equilibria of system~\cref{eq:san} and
their stability. Because system~\cref{eq:san} is a polynomial 
vector field, we use Poincar\'e compactification to project the
phase space into the three-dimensional open ball $\mB:=\mB^2(2)$ of
radius $2$. For the compactified model of \cref{eq:san} the sphere
$\mS:=\mS^2(2)$ bounding $\mB$ represents the dynamics at infinity 
\cite{Vel, Jes1, Jes2}; \cref{secApp:Com} gives more details on the
compactification, and shows the respective compactified vector
field~\cref{eq:ultVectField} of system~\cref{eq:san} and discusses its behaviour at
$\mS$. This compactification allows us to consider all equilibria and
continue them with \textsc{Auto}, even when they
interact with infinity. 
\begin{figure}
\centering
\includegraphics{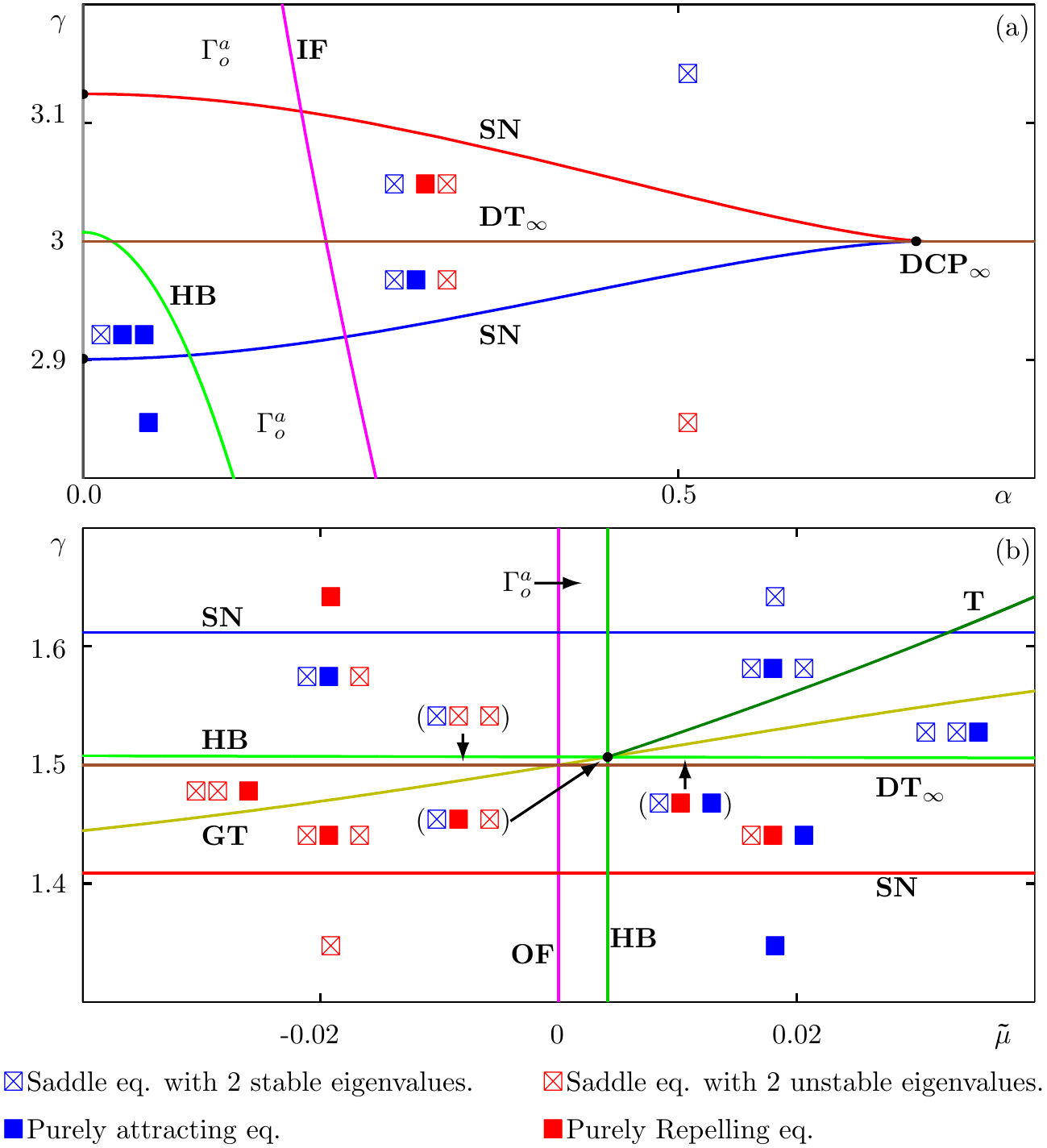}
\caption{Bifurcation diagram of equilibria for system
  \cref{eq:san}. Panel (a) shows the $(\alpha, 
  \gamma)$-plane with $(a,b,c,\beta,\mu,\tilde{\mu})=(0.22,1,-2,1,0,0)$
  and panel (b) the $(\tilde{\mu},\gamma)$-plane with
  $(a,b,c,\alpha,\beta,\mu)=(-0.5,2.5,-1,1,0,0)$.  Each equilibrium is
  represented by a square, where color and filling indicates its eigenvalue configuration as
  explained in the legend. Shown are saddle-node
  bifurcations $\mathbf{SN}$ as red and blue curves, inclination flip 
  $\mathbf{IF}$ and orbit flip bifurcations $\mathbf{OF}$ as pink curves,
  Hopf   bifurcation $\mathbf{HB}$ as a light-green curve, torus
  bifurcation $\mathbf{T}$ as a dark-green curve, transcritical bifurcation
   $\mathbf{GT}$ as a yellow curve, and nongeneric transcritical bifurcation at
  infinity $\mathbf{DT_\infty}$ as a brown curve; the point
  $\mathbf{DCP_\infty}$  is a nongeneric cusp  
   point at infinity. An attracting periodic orbit $\Gamma^a_o$ exist
  between $\mathbf{HB}$ and $\mathbf{IF}$ in panel (a), and between $\mathbf{HB}$ and
  $\mathbf{OF}$ in panel (b).  }\label{fig:BifInf} 
\end{figure}

\cref{fig:BifInf} shows the corresponding bifurcation diagrams of the
equilibria of system~\cref{eq:san}, as determined by using its
compactification~\cref{eq:ultVectField},  for both the inclination and orbit flip
cases. Specifically, panel (a) shows the inclination flip case
\textbf{IF} in the $(\alpha,\gamma)$-plane of the compactified 
model~\cref{eq:ultVectField} for 
$(a,b,c,\beta,\mu,\tilde{\mu})=(0.22,1,-2,1,0,0)$, and panel (b) the
orbit flip case \textbf{OF} in the $(\tilde{\mu},\gamma)$-plane for 
$(a,b,c,\alpha,\beta,\mu)=(-0.5,2.5,-1,1,0,0)$.  Starting with \textbf{IF}
in panel (a), we focus on the region near
$\gamma=3$.  The origin is always an equilibrium and the squares
indicate the number and stability of the additional
equilibria. A blue square corresponds to an equilibrium with at least
two stable eigenvalues; otherwise the square is red. Squares with a
cross refer to saddles and solid squares to sinks (blue) or
sources (red). System~\cref{eq:san} and its compactification are symmetric under
the transformation $(z,\alpha) \rightarrow (-z,-\alpha)$ when $\mu=0$; hence, there
is no need to show the negative values
of $\alpha$ in \cref{fig:BifInf}.  The blue curve represents a saddle-node 
bifurcation, labeled $\textbf{SN}$, that gives rise to a sink and a
saddle; similarly, the red curve $\textbf{SN}$ gives rise to a source and a saddle. The line
$\gamma=3$ (brown) represents a  degenerate transcritical bifurcation
at infinity, labeled \textbf{DT$_\infty$}. Along the
bifurcation curve \textbf{DT$_\infty$}, an equilibrium reaches
infinity, namely, at the non-hyperbolic equilibrium
$(0,0,2) \in \mS$. After the equilibrium crosses the curve
\textbf{DT$_\infty$}, its $z$-coordinate and  all its
eigenvalues change sign, that is, it reappears at $(0,0,-2)$ with the
opposite stability.  The
curves labeled $\textbf{SN}$ meet in a degenerate cusp   
point \textbf{DCP$_\infty$} on \textbf{DT$_\infty$}.  The curves \textbf{IF} (pink)
and \textbf{HB} (green) are bifurcations of homoclinic or periodic
orbits. Since $\mu=0$ there exists a homoclinic orbit
irrespective of the choices for $\alpha$ and $\gamma$ \cite{san1}.
The homoclinic orbit is orientable for small $\alpha$ and changes type
at the inclination flip curve \textbf{IF}, which is the pink curve; the homoclinic orbit
is non-orientable for values of $\alpha$ to the right of
\textbf{IF}.  The Hopf bifurcation \textbf{HB}  
gives rise to an orientable attracting periodic orbit $\Gamma^a_o$
that merges with the orientable homoclinic orbit and disappears in the
homoclinic flip bifurcation \textbf{IF}.

We choose to focus on the same situation that was studied in
\cite{Agu1}, namely, where there is a single extra
saddle-focus equilibrium $\textbf{q}\in\R^3$, with a two-dimensional unstable
manifold $W^u(\textbf{q})$ and a one-dimensional stable manifold
$W^s(\textbf{q})$. For this reason, we fix $\gamma=2$ and vary $\alpha>0$, which is
equivalent to the situation shown along the horizontal line $\gamma=2.8$ in
\cref{fig:BifInf}(a). 

\cref{fig:BifInf}(b) shows that there is a similar configuration of equilibria
for parameters $\tilde{\mu}$ and $\gamma$ for the \textbf{OF} case.
We again find a degenerate transcritical  
bifurcation \textbf{DT$_\infty$} at infinity at
$\gamma=1.5$,  two curves of Hopf
bifurcation \textbf{HB} and two saddle-node bifurcation curves
\textbf{SN}.  We also find a  curve of
torus bifurcation \textbf{T} (dark-green curve) and a generic transcritical
bifurcation \textbf{GT} (light-green).   There exists a curve \textbf{OF} of orbit flip bifurcations 
at $\tilde{\mu}=0$ in the $(\tilde{\mu},\gamma)$-plane. However,  the
homoclinic orbit that goes through these orbit flip bifurcations
cannot be found in this parameter plane.

For the case $\mathbf{OF}$, we again consider the situation where system~\cref{eq:san} has an additional equilibrium $\mathbf{q}$ with the
same properties as described before. For this reason, we
can study the orbit flip bifurcation by setting $\gamma=0$, which is equivalent
to the horizontal line $\gamma=1.3$ in \cref{fig:BifInf}(b).

\section{Inclination flip of case \textbf{B}} \label{sec:Inc}
\begin{figure}
\centering
\includegraphics[width=350pt]{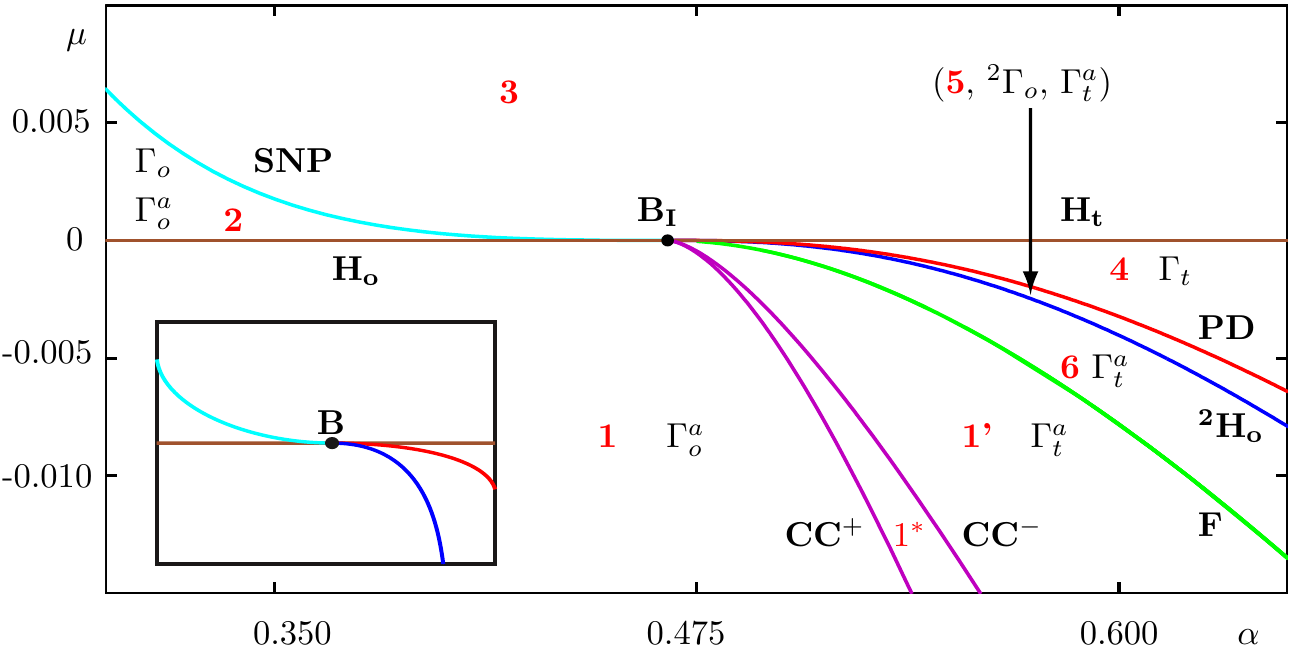}
\caption{Bifurcation diagram in the $(\alpha, \mu)$-plane near an inclination flip
    bifurcation  $\mathbf{B_I}$ of system~\cref{eq:san} for other
    parameters as given in \cref{sec:san}.  The inset shows
    only the curves of the theoretical unfolding of case $\mathbf{B}$
    \cite{kis1}. Shown are the homoclinic bifurcations $\mathbf{H_o}$
    and $\mathbf{H_t}$ as brown curves, the homoclinic bifurcation
    $^2\mathbf{H_o}$ as a blue curve, the saddle-node bifurcation $\mathbf{SNP}$
    of periodic orbit as a cyan curve, the period doubling bifurcation $\mathbf{PD}$
    as a red curve, the fold bifurcation $\mathbf{F}$ of heteroclinic
    orbits as green curve, and the loci $\mathbf{CC^{\pm}}$ as purple curves.}\label{fig:BDInc}
\end{figure}

We denote the inclination flip \textbf{IF} of type \textbf{B} by $\mathbf{B_I}$. On the level of the
codimension-one homoclinic bifurcation, $\mathbf{B_I}$ marks
the transition from an orientable homoclinic bifurcation to a non-orientable one by breaking
condition \textbf{(G3)}.  \cref{fig:BDInc} shows the unfolding of
$\mathbf{B_I}$ in the $(\alpha,\mu)$-plane for system~\cref{eq:san}
with the other parameters as stated in
\cref{sec:san}. The bifurcation curves that emanate from the
codimension-two point are a codimension-one orientable homoclinic
bifurcation \textbf{H$_o$} (brown curve), a codimension-one
non-orientable homoclinic bifurcation  \textbf{H$_t$} (brown 
curve), a saddle-node bifurcation of periodic
orbits \textbf{SNP} (cyan curve), a period-doubling bifurcation
\textbf{PD} (red curve) and a codimension-one homoclinic bifurcation
$\mathbf{^2H_o}$ (blue curve), as proven in \cite{kis1}. We also find an additional curve of fold
bifurcation of  heteroclinic orbits  \textbf{F} and  curves
$\mathbf{CC^{\pm}}$ (purple curves) that represent the moment that the Floquet
multipliers of  an attracting periodic orbit becomes complex
conjugates.  These curves divide the $(\alpha,\mu)$-plane in to open regions,
which are labeled by red numbers. Even though the curves $\mathbf{CC^{\pm}}$  are not
bifurcation curves, they bound region \mRed{1^*} in \cref{fig:BDInc}
where the Floquet multipliers of $\Gamma^a$ are complex
conjugates. Crossing through $\mathbf{CC^{\pm}}$ and this region
results in the transition of $\Gamma^a$  having a non-orientable to
having an orientable strong stable manifold, so that this attracting periodic orbit can
bifurcate at the curves \textbf{SNP} and \textbf{PD}, respectively.

Starting from region \mRed{1}, where an 
orientable attracting periodic orbit $\Gamma^a_o$ exists,  we move to
region \mRed{2} through \textbf{H$_o$}.  This homoclinic orbit 
creates an orientable saddle periodic orbit  $\Gamma_o$ in region
\mRed{2}, which disappears with $\Gamma^a_o$ in the \textbf{SNP}
bifurcation as we cross in to region \mRed{3}.  The
transition between regions \mRed{3} and \mRed{4}  is the 
\textbf{H$_t$} bifurcation. As in region \mRed{2},  the homoclinic orbit  becomes
a saddle periodic orbit $\Gamma_t$ in 
region \mRed{4}, but this saddle periodic orbit is
non-orientable.  As we move to region \mRed{5}, the periodic orbit $\Gamma_t$ undergoes
the \textbf{PD}  bifurcation and
becomes the non-orientable attracting periodic orbit $\Gamma^a_t$; furthermore, an orientable
saddle periodic orbit $^2\Gamma_o$ with twice the period of
$\Gamma_t$ is created.  Next, the transition between regions \mRed{5}
and \mRed{6} is characterized by the disappearance of  $^2\Gamma_o$
in $\mathbf{^2H_o}$ (blue curve) as it becomes an homoclinic orbit.
As shown in \cref{fig:BDInc}. The curve  \textbf{F} delimits
region~\mRed{6} and marks the creation of a pair of
heteroclinic orbits from $\mathbf{q}$ to $\mathbf{0}$ that exists in
regions \mRed{4}, \mRed{5} and \mRed{6}.  These heteroclinic orbits represent the
transverse intersection between $W^s(\mathbf{0})$ and
$W^u(\mathbf{q})$, which becomes tangent at
\textbf{F}  so that the two heteroclinic orbits merge and then
disappear in region \mRed{1'}.  Region \mRed{1'}  is topologically
equivalent to region \mRed{1} but the attracting periodic
orbit $\Gamma^a_t$ is non-orientable instead of orientable for
$\Gamma^a_o$ in region \mRed{1}. The transition from region~\mRed{1'}
to region~\mRed{1} occurs via a crossing of the curves $\mathbf{CC^-}$ and 
$\mathbf{CC^+}$ where  the Floquet multipliers of $\Gamma^a_t$ ($\Gamma^a_o$) in region
\mRed{1'} (\mRed{1}) change from being 
real positive (negative) to complex conjugate. In region~$\mRed{1^*}$,
bounded by $\mathbf{CC^+}$ and $\mathbf{CC^-}$,  the 
 periodic orbit does not have a strong stable
manifold.   \cref{secApp:BVP}  gives details on the
computation of the curves \textbf{F}, 
$\mathbf{CC^+}$ and $\mathbf{CC^-}$.  

Our goal is now to characterize the 
topological properties of the global manifolds in a neighborhood of the
inclination flip bifurcation.  We use the 
bifurcation diagram in \cref{fig:BDInc} as a 
reference to describe the changes in the organization of the manifolds of
system~\cref{eq:san} in phase space, as $\alpha$ and $\mu$ vary
between the different regions in the $(\alpha,\mu)$-plane.
\cref{tab:Inc1} provides an overview of the representative values of 
$\alpha$ and $\mu$ we selected from each region. We also illustrate the
manifolds for representative parameter points approximately on $\mathbf{H_o}$, $\mathbf{H_t}$,
$\mathbf{^2H_o}$, \textbf{F} and at $\mathbf{B_I}$; these values are
given in \cref{tab:Inc2}. We first present their phase portrait
in $\R^3$, where the orbit segments that forms the  two-dimensional
stable (unstable) manifolds are computed by restricting one end point
to lie in the sphere $\mS^*:=\{x \in \R^{3} : \| x-c\|=R\}$ with
$c:=(c_x,c_y,c_z)=(0.5,0,0)$ and $R=0.6$.

\begin{table}
\begin{center}
\begin{tabular}{|c|r@{.}l|r@{.}l|r@{.}l|r@{.}l|r@{.}l|r@{.}l|r@{.}l|} 
\hline
\textbf{Region} & \multicolumn{2}{c|}{\textbf{1}} &
\multicolumn{2}{c|}{\textbf{2}} & \multicolumn{2}{c|}{\textbf{3}} &
\multicolumn{2}{c|}{\textbf{4}} & \multicolumn{2}{c|}{\textbf{5}} &
\multicolumn{2}{c|}{\textbf{6}} & \multicolumn{2}{c|}{\textbf{1'}} \\ \hline
$\alpha$ & 0 & 300 &  0 & 300 & 0 & 650 & 0 & 650 & 0 & 650 & 0 & 650 & 0 & 650\\ \hline
$\mu$ & $-0$ & 004 &  0 & 004 & 0 & 004 & $-0$ & 004 & $-0$ & 007 &
$-0$ & 010 & $-0$ & 014\\  
\hline
\end{tabular}
\vspace{2mm}
\caption{Chosen representative parameter values for the different open
  regions in \cref{fig:BDInc}.} \label{tab:Inc1}
\end{center}

\begin{center}
\begin{tabular}{|c|r@{.}l|r@{.}l|r@{.}l|r@{.}l|r@{.}l|}
\hline
\textbf{Curve} & \multicolumn{2}{c|}{$\mathbf{H_o}$} &
\multicolumn{2}{c|}{$\mathbf{B_I}$} & \multicolumn{2}{c|}{$\mathbf{H_t}$} &
\multicolumn{2}{c|}{$\mathbf{^2H_o}$} & \multicolumn{2}{c|}{\textbf{F}} \\ \hline
$\alpha$ & 0 & 3000000 &  0 & 4664012& 0 & 6500000 & 0 & 6500000 & 0 &
6500000\\\hline 
$\mu$ & \multicolumn{2}{c|}{0} &  \multicolumn{2}{c|}{0} & \multicolumn{2}{c|}{0} & -0 & 0079047 & $-0$& 0134990\\\hline
\end{tabular}
\vspace{2mm}
\caption{Chosen representative parameter values at selected bifurcations in
  \cref{fig:BDInc}.} \label{tab:Inc2}
\end{center}
\end{table}

\begin{figure}
\centering
\includegraphics{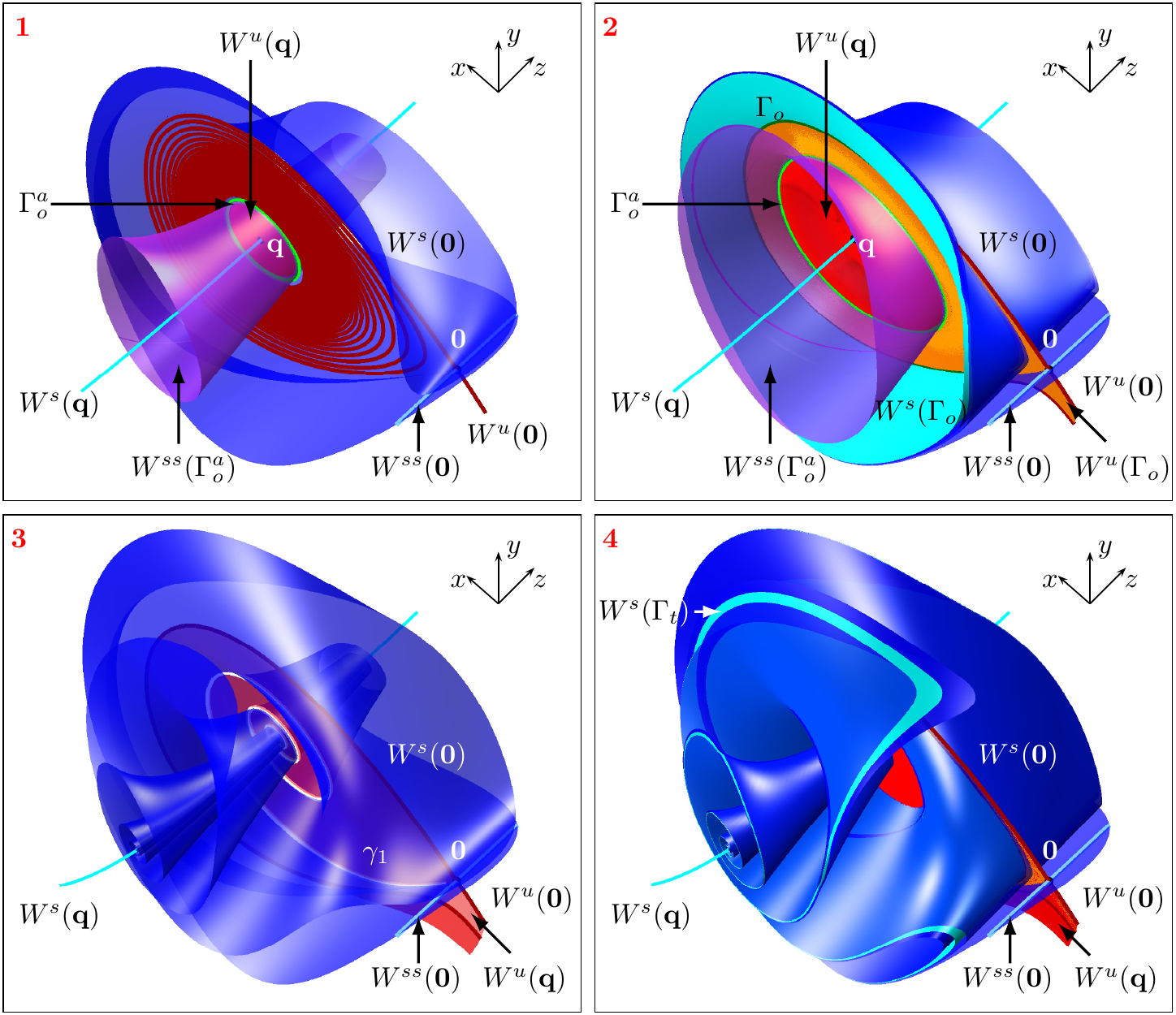}
\caption{Phase portraits of system~\cref{eq:san} in the different
  regions \mRed{1}-\mRed{6} and at the bifurcations $\mathbf{H_o}$ and
  $\mathbf{F}$ of the $(\alpha,\mu)$-plane in \cref{fig:BDInc}. Shown are
  $W^s(\mathbf{0})$ as a dark-blue surface, 
  $W^{ss}(\mathbf{0})$ as a blue curve, $W^{u}(\mathbf{0})$ as a
  pink curve, $W^{u}(\textbf{q})$ as a red surface,
  $W^{s}(\textbf{q})$ as a cyan curve,  $W^s(\Gamma_o)$ and
  $W^s(\Gamma_t)$ as  cyan surfaces, $W^u(\Gamma_o)$ and
  $W^u(\Gamma_t)$  as  orange surfaces, and $W^{ss}(\Gamma^a_t)$ and
  $W^{ss}(\Gamma^a_o)$ as purple surfaces. The $(\alpha,\mu)$-values
  for each panel are  given  in \cref{tab:Inc1} and
  \cref{tab:Inc2}.  See also the accompanying animation ({\color{red} GKO\_Bflip\_animatedFig6-1.gif}).}\label{fig:Sections}
\end{figure}
\setcounter{figure}{5}
\begin{figure}
\centering
\includegraphics{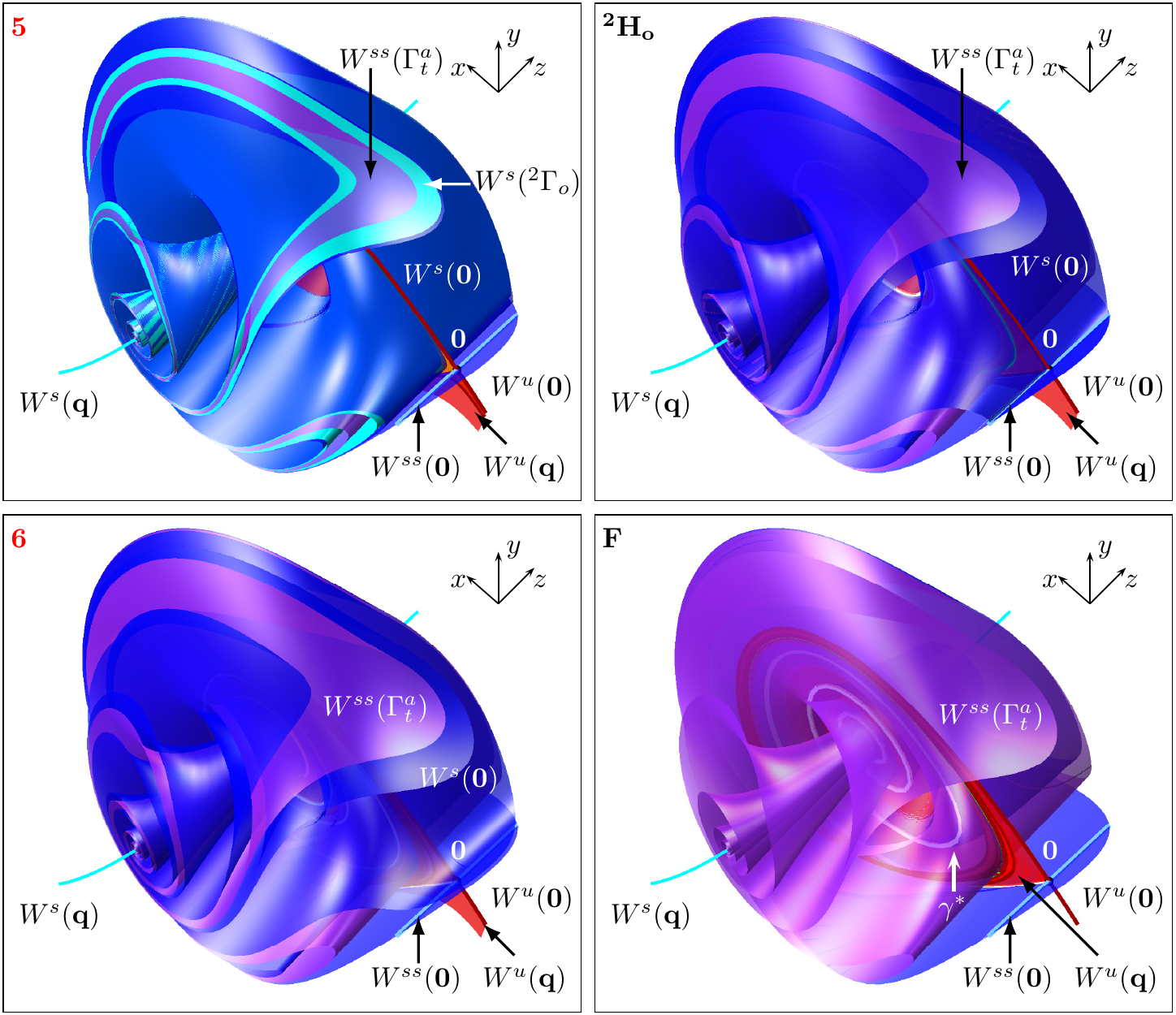}
\caption{Continued. See also the accompanying animation ({\color{red} GKO\_Bflip\_animatedFig6-2.gif}).}
\end{figure}

\subsection{Manifolds in the open regions near $\mathbf{B_I}$}
\cref{fig:Sections} shows phase portraits in each region and
at the bifurcations $\mathbf{^2H_o}$ and \textbf{F}. Specifically,
we show the equilibria $\mathbf{0}$ and $\mathbf{q}$ along with their
stable and unstable manifolds, as well as the periodic orbits and
their manifolds when they exist. In the following, we
cycle around $\mathbf{B_I}$ through the bifurcation diagram 
in \cref{fig:BDInc}, starting from region \mRed{1}, and describe the
transitions on the level of the invariant 
manifolds in phase space. To improve clarity and understanding of
\cref{fig:Sections}, the accompanying animations
({\color{red} GKO\_Bflip\_animatedFig6-1.gif}) and ({\color{red} GKO\_Bflip\_animatedFig6-2.gif}),
show the respective objects rotating clockwise around the $y$-axis.

\subsubsection{Manifolds in region 1}
Region \mRed{1} is characterized by the existence of an orientable
attracting periodic orbit $\Gamma^a_o$.  The corresponding phase
portrait in panel~\mRed{1} of \cref{fig:Sections} shows how one branch
of $W^u(\mathbf{0})$ (red curve) spirals towards $\Gamma^a_o$  (green
curve).  The two-dimensional stable manifold $W^s(\mathbf{0})$ (blue 
surface) folds over $W^u(\mathbf{0})$ and trajectories on
$W^s(\mathbf{0})$ escape towards infinity in backward time.
Furthermore, the two-dimensional unstable manifold 
$W^u(\mathbf{q})$ (red surface) accumulates on $\Gamma^a_o$.  
In fact, $\Gamma^a_o$ is the boundary of $W^u(\mathbf{q})$.
Since the Floquet multipliers of $\Gamma^a_o$ are positive, its
strong stable manifold $W^{ss}(\Gamma^a_o)$
(purple surface) is a topological cylinder.  
We note that the one-dimensional stable manifold $W^s(\mathbf{q})$
(cyan curve) lies in the interior of $W^{ss}(\Gamma^a_o)$.  Therefore,
none of the other stable manifolds outside of  
$W^{ss}(\Gamma^a_o)$ can accumulate onto $W^s(\mathbf{q})$ in backward
time.

\subsubsection{Manifolds in region 2}
The bifurcation curve $\mathbf{H_o}$, between region~\mRed{1} and
\mRed{2}, creates the homoclinic orbit $\mathbf{\Gamma}_{\rm hom}$.  As we transition to region \mRed{2}, the orbit
$\mathbf{\Gamma}_{\rm hom}$ becomes the orientable saddle periodic orbit $\Gamma_o$
(dark-green curve) in panel \mRed{2} of \cref{fig:Sections}. It has
two-dimensional stable and unstable 
manifolds $W^s(\Gamma_o)$  (cyan surface) and $W^u(\Gamma_o)$ (orange
surface), respectively.  Since $\Gamma_o$ is an orientable saddle
periodic orbit, both $W^s(\Gamma_o)$ and $W^u(\Gamma_o)$ are
orientable, but  $W^u(\Gamma_o)$ is bounded by $\Gamma^a_o$ and
$W^u(\mathbf{0})$, while $W^s(\Gamma_o)$ is unbounded.   As shown in
panel~\mRed{2} of \cref{fig:Sections},  the one-dimensional manifold 
$W^u(\mathbf{0})$ no longer accumulates on $\Gamma_o^a$, but one
branch folds over  $W^s(\mathbf{0})$ before both branches move off to
infinity. Furthermore, $W^s(\mathbf{0})$ now accumulates (in backward
time) onto $W^s(\Gamma_o)$.  Note that
$W^s(\mathbf{0})$ intersects $W^u(\Gamma_o)$  transversally; this
implies the existence of a heteroclinic cycle-to-point connecting orbit from
$\Gamma_o$ to $\mathbf{0}$, which exists in the open region \mRed{2} in \cref{fig:BDInc}.

\subsubsection{Manifolds in region 3}\label{sec:Inc1to2}
The transition between regions \mRed{2} and \mRed{3} occurs at the
saddle-node bifurcation of periodic orbits \textbf{SNP}.  At \textbf{SNP}, the periodic orbits $\Gamma_o$, $\Gamma_o^a$ and their
manifolds $W^{s}(\Gamma_o)$ and $W^{ss}(\Gamma^a_o)$ merge, and
disappear as we transition in to region \mRed{3}.  Consequently, $W^{s}(\mathbf{0})$  now spirals towards \textbf{q} and accumulates on
$W^{s}(\mathbf{q})$ in backward time; see panel \mRed{3} of
\cref{fig:Sections}.  The manifold $W^{u}(\mathbf{0})$ is now  the
boundary of $W^{u}(\mathbf{q})$ and the manifolds $W^{s}(\mathbf{0})$ and  
$W^{u}(\mathbf{q})$ intersect transversally  in region \mRed{3}; this
implies the existence of a heteroclinic orbit $\gamma_1$ (white curve)
from $\mathbf{q}$ to $\mathbf{0}$. 

\subsubsection{Manifolds in region 4} \label{sec:Lambda}
Regions \mRed{3} and \mRed{4} are separated by a curve $\mathbf{H_t}$
of codimension-one non-orientable homoclinic bifurcations.  The
homoclinic orbit $\mathbf{\Gamma}_{\rm hom}$ becomes the
non-orientable saddle periodic orbit $\Gamma_t$ in region~\mRed{4}. The stable manifold
$W^{s}(\mathbf{0})$ accumulates onto $W^s(\Gamma_t)$
(cyan surface) in backward time.  In contrast to
region~\mRed{2}, the non-orientable stable manifold $W^s(\Gamma_t)$ is
not a separatrix but spirals towards \textbf{q}
and accumulates in backward time onto $W^s(\mathbf{q})$. Furthermore,
we see intersections between the different manifolds in region~\mRed{4},
although it is hard to appreciate their structure. 

\begin{figure}
\centering
\includegraphics{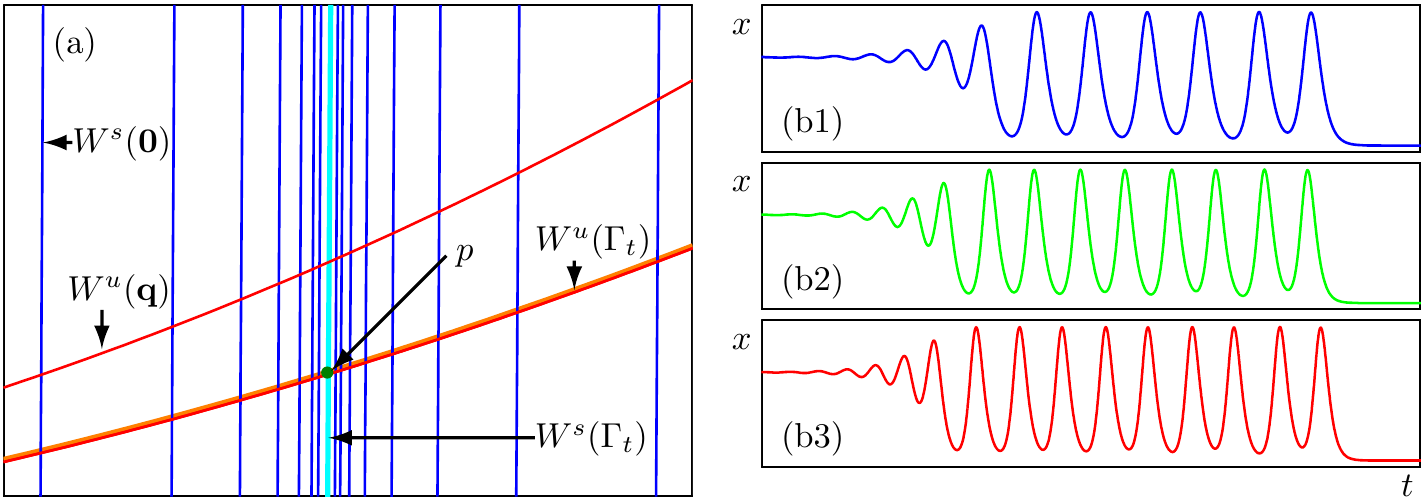}
\caption{Intersection of global manifolds in region \mRed{4}. Panel
  (a) shows the intersections of the manifolds with the plane $\Sigma$; points in
  $\Sigma$ have the same $y$-coordinate as $\mathbf{q}$. Shown are
  $W^s(\mathbf{0})$ as blue curves, $W^s(\Gamma_o)$ as a cyan curve,
  $W^u(\Gamma_t)$ as an orange curve, $W^u(\mathbf{q})$ as red curves and the point $p \in \Gamma_t \cap
  \Sigma$ as a green dot.  Panel (b)
  shows (scaled) time series in $x$ of representative heteroclinic 
  orbits from $\mathbf{q}$ to $\mathbf{0}$, respectively.}\label{fig:infHetero}  
\end{figure}

To illustrate the nature of these manifold interactions we consider their
intersection sets with the plane $\Sigma:= \lp\{ (x,y,z) \in \R^3 :  y =
\mathbf{q}_y \rp\}$, where $\mathbf{q}_y$ is the $y$-component of
$\mathbf{q}$. \cref{fig:infHetero}(a)
shows the intersection sets of $W^s(\mathbf{0})$, $W^u(\mathbf{q})$,
$W^s(\Gamma_t)$ and $W^u(\Gamma_t)$ with $\Sigma$ in a neighborhood of
one of the two points $p \in \Gamma_t \cap \Sigma$. Locally near $p$,
there is a single curve (cyan) representing $W^s(\Gamma_t) \cap
\Sigma$ and another single curve (orange) representing $W^u(\Gamma_t) 
\cap \Sigma$.  Since $W^u(\mathbf{q}) \cap \Sigma$ intersects 
$W^s(\Gamma_t) \cap \Sigma$ and $W^u(\Gamma_t) \cap \Sigma$ intersects
$W^s(\mathbf{0}) \cap \Sigma$, there exist structurally 
stable heteroclinic orbits from  $\mathbf{q}$ to $\Gamma_t$ and
from $\Gamma_t$ to $\mathbf{0}$, respectively. As a consequence of the
$\lambda$-lemma \cite{Taken1,Wigg1}, the intersection sets $W^s(\mathbf{0}) \cap \Sigma$
(blue) and $W^u(\mathbf{q}) \cap \Sigma$ (red) give rise to several
curves in the neighborhood of $p$ that accumulate onto $W^s(\Gamma_t)\cap
\Sigma$ and $W^u(\Gamma_t)\cap \Sigma$, respectively.  Therefore,
there exist transversal intersections between the sets $W^u(\mathbf{q})
\cap \Sigma$ and $W^s(\mathbf{0}) \cap \Sigma$, which imply the
existence of structurally stable heteroclinic orbits from $\mathbf{q}$ to
$\mathbf{0}$. Panels  (b1)--(b3) of \cref{fig:infHetero}
show the evolution of the $x$-variable with respect to time for three heteroclinic
orbits from $\mathbf{q}$ to $\mathbf{0}$ in region \mRed{4}; observe
how they differ in the number of big excursions before converging to 
$\mathbf{0}$, these excursions correspond to intersection points in
$\Sigma$ close to $p$.  The $\lambda$-lemma guarantees the sets
$W^u(\mathbf{q}) \cap \Sigma$ and 
$W^s(\mathbf{0}) \cap \Sigma$ intersect in \emph{an arbitrary small
  neighborhood} of $p$; therefore there exist infinitely many
intersection points. Only a finite number of these intersection points corresponds
to a single heteroclinic orbit from $\mathbf{q}$ to $\mathbf{0}$;
hence, there are indeed infinitely
many heteroclinic orbits from $\mathbf{q}$ to $\mathbf{0}$ in region
\mRed{4}.

We remark that the manifolds $W^u(\mathbf{q})$ and $W^s(\mathbf{0})$ shown in panel \mRed{4}  of
\cref{fig:Sections} are only computed up to the first two of their infinitely many
layers that intersect $\mS^*$; hence, the accumulation of these
manifold with the respective invariant manifolds of $\Gamma_o$ is not visible in panel \mRed{4} of
\cref{fig:Sections}.

\subsubsection{Manifolds in region 5}
When crossing from region \mRed{4} to region \mRed{5} a period-doubling
bifurcation \textbf{PD} occurs. The saddle periodic orbit $\Gamma_t$ becomes
a non-orientable attracting periodic orbit $\Gamma^a_t$, and an orientable saddle periodic
orbit $^2\Gamma_o$  with twice the period of $\Gamma^a_t$  emanates
from the period-doubling bifurcation into region
\mRed{5}.

Panel \mRed{5} in \cref{fig:Sections} shows that $W^s(^2\Gamma_o)$
(cyan) and $W^s(\mathbf{0})$ accumulate onto $\mathbf{q}$ and
$W^s(\mathbf{q})$ in backward time.  The periodic orbit $\Gamma^a_t$
is attracting in region \mRed{5}, but its strong stable manifold
$W^{ss}(\Gamma^a_t)$ (purple) can be viewed as the continuation of
$W^s(\Gamma_t)$.  The portion of $W^{ss}(\Gamma^a_t)$ relative to 
$W^s(^2\Gamma_o)$ suggests that the basin of attraction
$\mathcal{B}(\Gamma^{a}_t)$ of $\Gamma^a_t$ is bounded by 
$W^s(^2\Gamma_o)$.  Indeed, one side of $W^u(^2\Gamma_o)$ accumulates
onto $\Gamma^a_t$, while the other side intersects $W^s(\mathbf{0})$.
Hence, the situation is very similar to that in region \mRed{4}: there
exists  one transversal heteroclinic
orbit from $^2\Gamma_o$ to $\mathbf{0}$, and there exist infinitely many
heteroclinic orbits from $\mathbf{q}$ to $\mathbf{0}$. Furthermore, a two-dimensional submanifold of
$W^u(\mathbf{q})$ lies in the open set $\mathcal{B}(\Gamma^{a}_t)$;
hence this submanifold accumulates on $\Gamma^{a}_t$ and
its boundary corresponds to an intersection of $W^u(\mathbf{q})$
and $W^s(^2\Gamma_o)$, that is, 
there exist transversal heteroclinic orbits from $\mathbf{q}$ to
$^2\Gamma_o$.  Also, as in region \mRed{4}, the
one-dimensional unstable manifold  $W^u(\mathbf{0})$ is contained in
part of the closure of  both $W^u(^2\Gamma_o)$ and $W^u(\mathbf{q})$.

\subsubsection{Manifolds in region 6}
The boundary between region \mRed{5} to \mRed{6} is the curve
$\mathbf{^2H_o}$ of codimension-one orientable homoclinic
bifurcation.  The moment of the homoclinic bifurcation is illustrated
in panel $\mathbf{^2H_o}$  of \cref{fig:Sections}. The limit of the
saddle periodic orbit $^2\Gamma_o$ is the orientable codimension-one homoclinic
orbit $^2\mathbf{\Gamma_{\rm hom}}$ at $\mathbf{^2H_o}$. Note that
$^2\Gamma_o$ and its manifolds have 
disappeared, and so have the heteroclinic orbits connecting
$^2\Gamma_o$ with $\mathbf{0}$ and $\mathbf{q}$. Hence there are no
longer infinitely many codimension-zero
heteroclinic orbits from $\mathbf{q}$ to $\mathbf{0}$.  We find that
$W^s(\mathbf{0})$ interacts non-trivially with $W^u(\mathbf{q})$ in
two transversal  heteroclinic orbits from
$\mathbf{q}$ to $\mathbf{0}$ that persist through the homoclinic
bifurcation $\mathbf{^2H_o}$.  They bound the two-dimensional
submanifold $W^u(\mathbf{q})$ that accumulates on
$\Gamma^{a}_t$. Since these two heteroclinic orbits can be viewed as
the continuation of the two heteroclinic orbits from $\mathbf{q}$ to
$^2\Gamma_o$ in region~$\mRed{5}$. The other
infinitely many heteroclinic orbits from $\mathbf{q}$ to $\mathbf{0}$
all disappear at once in the homoclinic bifurcation $\mathbf{^2H_o}$.

The homoclinic orbit $^2\mathbf{\Gamma_{\rm hom}}$
disappears but the non-orientable attracting periodic orbit $\Gamma^a_t$ and the
two transversal heteroclinic orbits from $\mathbf{q}$ to $\mathbf{0}$
persist in region~\mRed{6}.  In particular,
these two heteroclinic orbits still bound the portion of
$W^u(\mathbf{q})$  that is attracted by $\Gamma^a_t$.  Note that the
branch of $W^u(\mathbf{0})$ that formed $^2\mathbf{\Gamma_{\rm hom}}$
now spirals towards $\Gamma^a_t$. It is
worth noting that the phase space in region~\mRed{6} is
topologically equivalent to that of region~\mRed{3} for case
\textbf{A} in \cite{Agu1}.

\subsubsection{Manifolds in region 1'}
At the curve \textbf{F}, which is the transition from
region~\mRed{6} to region~\mRed{1'}, the two-dimensional manifolds
$W^s(\mathbf{0})$ and $W^u(\mathbf{q})$ lose their  two intersection orbits
in a quadratic tangency; see panel \textbf{F} of
\cref{fig:Sections}. Hence, the two heteroclinic orbits merge to
become the heteroclinic orbit $\gamma^*$, representing the last moment where
$W^u(\mathbf{0})$ is part of the boundary of $W^u(\mathbf{q})$. 

In region $\mRed{1'}$, the manifolds  $W^s(\mathbf{0})$ and
$W^u(\mathbf{q})$ do no longer interact with each other, and
$W^u(\mathbf{q})$ accumulates entirely on $\Gamma^a_t$. The nontrivial Floquet
multipliers of $\Gamma^a_t$ in region~\mRed{1'}  become  equal at the
curve $\mathbf{CC^-}$;  they are then complex conjugates with negative real
part close to $\mathbf{CC^-}$ in region~$\mRed{1^*}$. Hence, there is
not a well-defined strong stable manifold of $\Gamma^a$ in region
$\mRed{1^*}$.  As we approached 
region~\mRed{1} from region $\mRed{1^*}$, the Floquet multipliers of
$\Gamma^a$ cross the imaginary axis and become complex conjugates
with positive real part close to the curve $\mathbf{CC^+}$.   At this
curve, the non-trivial Floquet multipliers are both the same positive real number.
They then become two distinct positive real values  in
region~\mRed{1}, so that $\Gamma^a_o$ has a well-defined strong stable
manifold again. This transition through $\mathbf{CC^-}$ and
$\mathbf{CC^+}$ allows the twisted periodic orbit $\Gamma^a_t$ to
become the orientable $\Gamma^a_o$ that then disappears with $\Gamma_o$ at the 
bifurcation \textbf{SNP} \cite{Hin1}.  Our numerical
computations indicate that the two curves
 $\mathbf{CC^{+}}$ and $\mathbf{CC^-}$ are not tangent to the
homoclinic bifurcation curve at $\mathbf{B_I}$ but approach
this codimension-two point  at a non-zero angle; see \cref{fig:BDInc}.  Since the
manifolds in regions $\mRed{1'}$ and $\mRed{1^*}$  
are qualitatively the same as in region \mRed{1}, except for the
properties of the strong stable manifold $W^{ss}(\Gamma^a)$, we do not
show the respective phase portraits  in \cref{fig:Sections}.

\subsection{Intersections of the invariant manifolds with a sphere} \label{sec:IntSphe}
It is a challenge to extract the precise nature of the phase portraits
in the panels of \cref{fig:Sections} in 
terms of the re-organization of the basins of attracting periodic
orbits.  Therefore, we now study the intersection sets of the respective invariant
manifolds with the sphere $\mS^*$ of radius $R=0.6$ centered at
$c=(c_x,c_y,c_z)=(0.5,0,0)$.  Since $\mS^*$ is a compact set, all
intersection sets of the  manifolds of system~\cref{eq:san} must be
bounded. We consider the intersection sets:
\begin{gather*}
\widehat{W}^s(\textbf{0}) := W^s(\mathbf{0}) \cap \mS^* \text{ , }
\widehat{W}^{ss}(\textbf{0}) := W^{ss}(\mathbf{0}) \cap \mS^* \text{ , }
\widehat{W}^s(\textbf{q}) := W^s(\textbf{q}) \cap \mS^*,\\
\widehat{W}^s(\Gamma_o) := W^s(\Gamma_o) \cap \mS^* \text{ , }
\widehat{W}^{s}(\Gamma_t) := W^{s}(\Gamma_t) \cap \mS^* \text{ and }
\widehat{W}^{ss}(\Gamma^a_{o / t}) := W^{ss}(\Gamma^a_{o / t}) \cap
\mS^*.
\end{gather*}
In particular, the intersection sets of all two-dimensional manifolds
that are transverse to $\mS^*$ are curves, while the
one-dimensional manifolds intersect $\mS^*$ in points. We also determine the
regions on $\mS^*$ that correspond to the intersection sets of the basin of attraction
$\mathcal{B}(\Gamma^{a})$ of $\Gamma^a$;  we denote this set
$\widehat{\mathcal{B}}(\Gamma^{a})$ and color it yellow in the subsequent figures. 

It is convenient to represent these intersection sets in the plane; to
this end, we use stereographic
projection onto the $(x,z)$-plane via the transformation 
\begin{equation}\label{eq:ste}
(x',y',z') \in \mS^* \mapsto \lp( \frac{R(x'-c_x)}{R+(y'-c_y)} ,
\frac{R(z'-c_z)}{R+(y'-c_y)} \rp) \in \R^2.  
\end{equation}
This transformation translates $c$ to $\mathbf{0}$, and then projects
a point on the (translated) sphere $\mS^*$ along the line through
$(0,-R,0)$ to a point on the tangent plane of the sphere at $(0,R,0)$, 
 that is, the  plane parallel to the $(x,z)$-plane through $(0,R,0)$.
\cref{fig:BDI1} shows the intersection sets with $\mS^*$
in each region close to the inclination flip. As in the previous section,
\cref{fig:BDI1} starts with the situation for region \mRed{1} and
cycles through the bifurcation diagram of the inclination flip
bifurcation.  However, now we show also the situation at
region~\mRed{1'}. The left column of \cref{fig:BDI1} shows stereographic projections of
the intersections sets of the manifolds in each region close to the inclination flip 
as computed with \textsc{Auto} \cite{Doe1,Doe2}.   The right column
shows topological sketches of these projections to illustrate and 
accentuate important features. 

\subsubsection{Intersection sets in regions 1 to 3}
\begin{figure}
\centering
\includegraphics{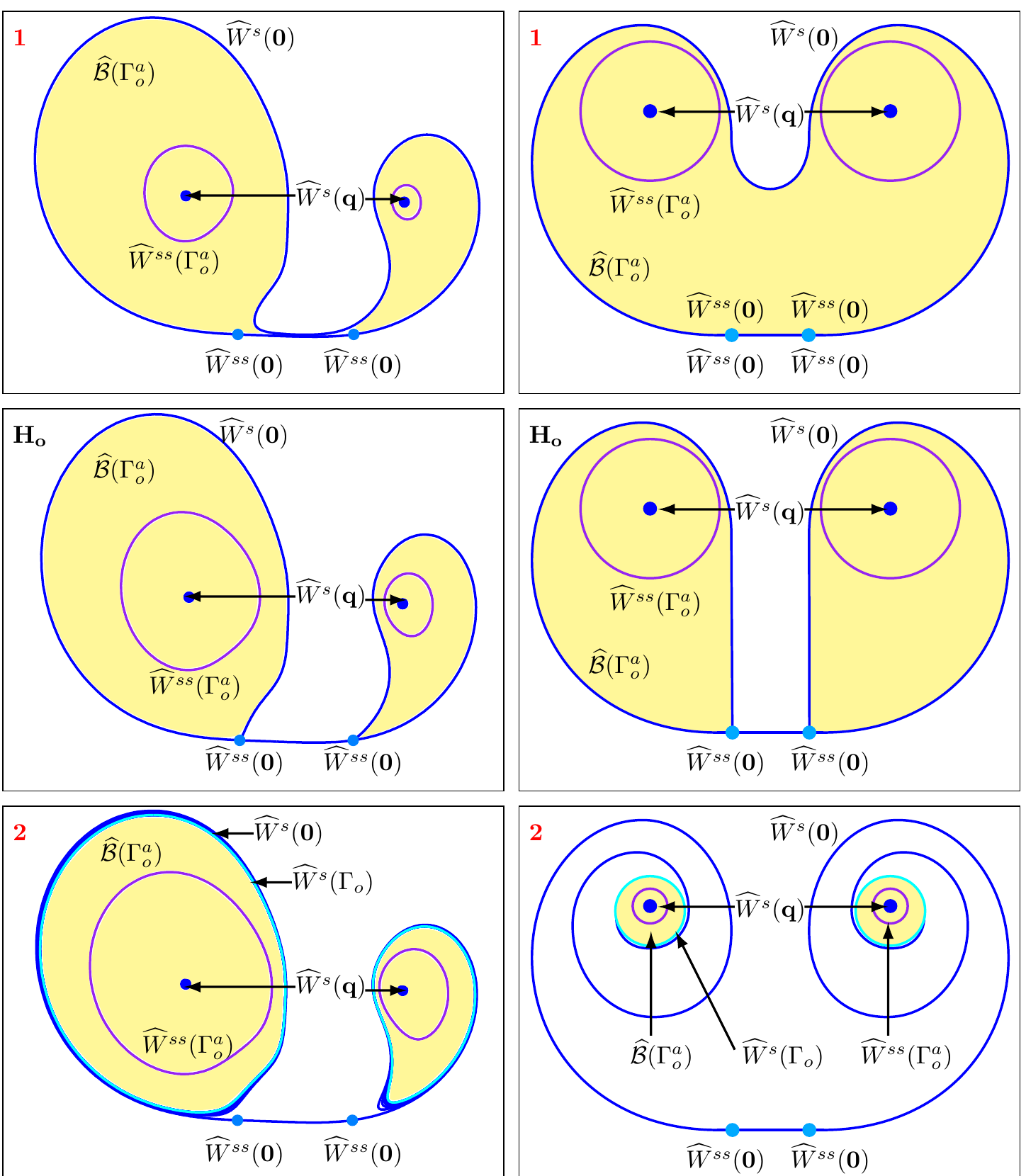}
\caption{Stereographic projections of the intersection  sets of the
  invariant manifolds with $\mS^*$ in the regions and at bifurcations of the
  bifurcation diagram in \cref{fig:BDInc} near the inclination flip $\mathbf{B_I}$; the
  first column shows the computed manifolds of system~\cref{eq:san}
  and the second column are topological sketches.   Shown are
  $\widehat{W}^s(\textbf{0})$ as dark-blue curves,
  $\widehat{W}^{ss}(\textbf{0})$ as light-blue dots and $\widehat{W}^{s}(\textbf{q})$
  dark-blue dots,   $\widehat{W}^s(\Gamma_o)$, $W^{s}(^2\Gamma_o)$ and
  $\widehat{W}^s(\Gamma_t)$ as cyan curves,
  $\widehat{W}^{ss}(\Gamma^a_{o/t})$ as purple curves and
  $\widehat{\mathcal{B}}(\Gamma^{a}_{o/t})$ as a shaded 
  yellow region. For respective
  parameter values see \cref{tab:Inc1} and \cref{tab:Inc2}.}\label{fig:BDI1} 
\end{figure}
\setcounter{figure}{7}

\begin{figure}
\centering
\includegraphics{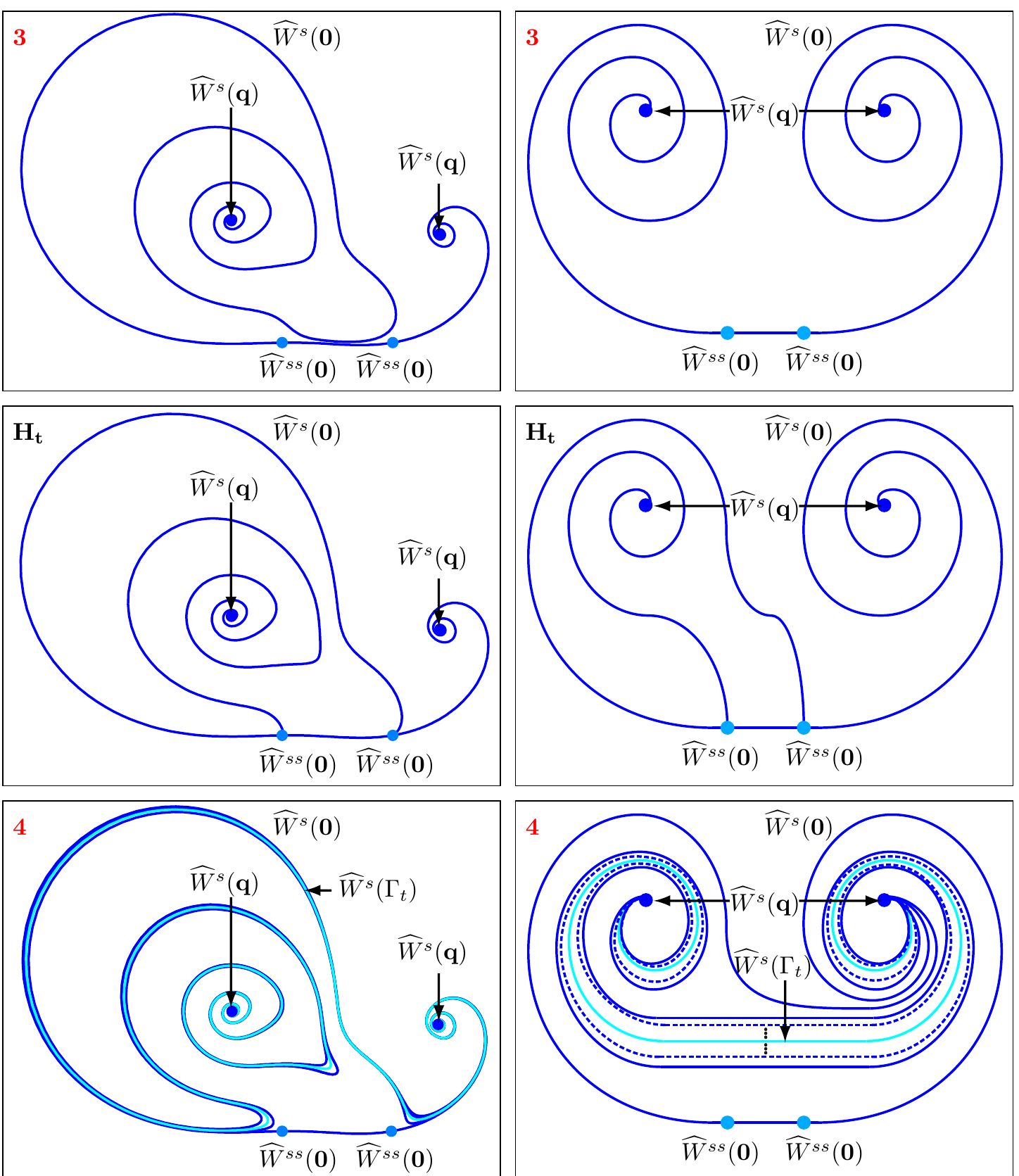}
\caption{Continued.}\label{fig:BDI2}
\end{figure}

\setcounter{figure}{7}
\begin{figure}
\centering
\includegraphics{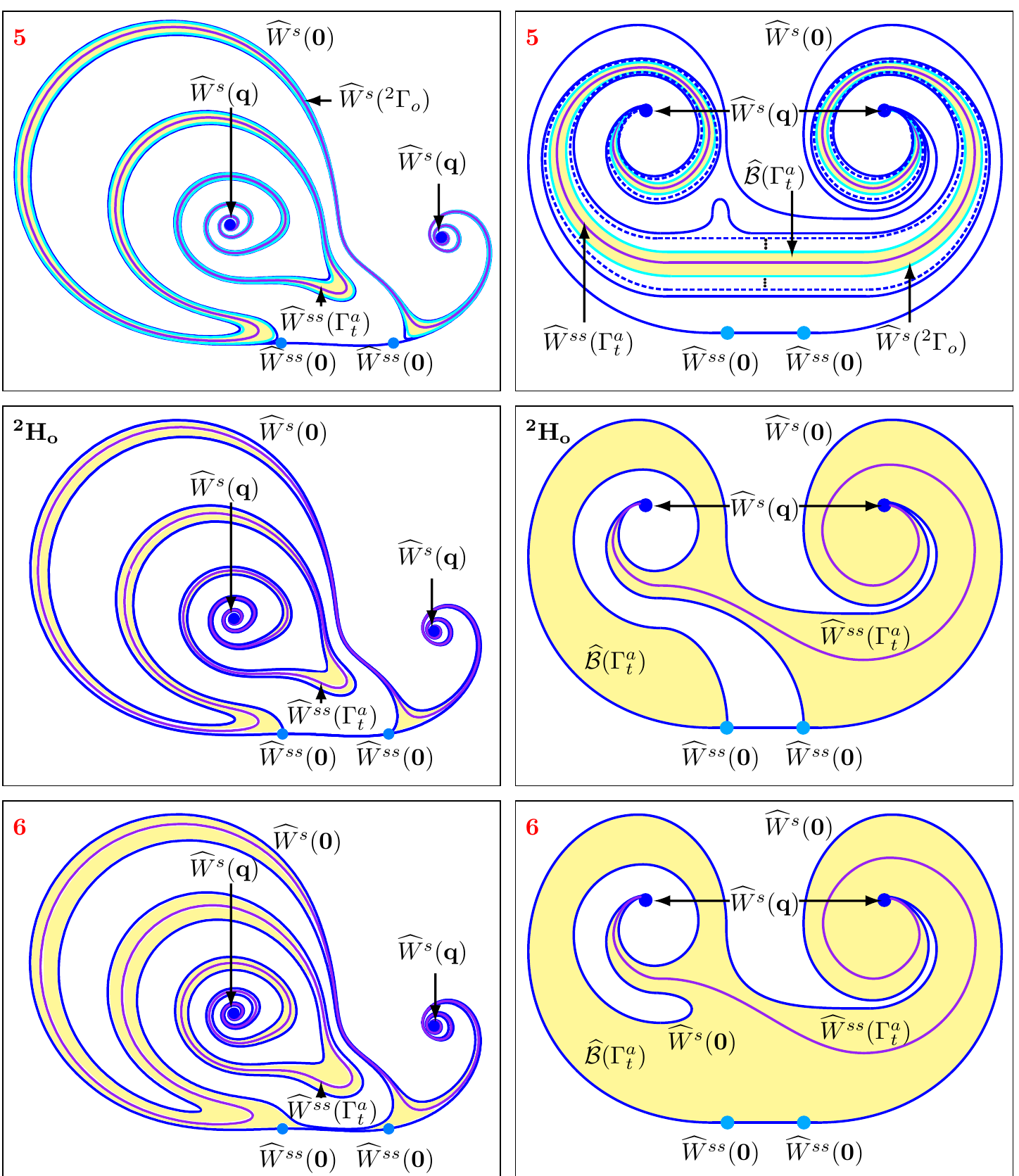}
\caption{Continued.}\label{fig:BDI3}
\end{figure}

\setcounter{figure}{7}
\begin{figure}
\centering
\includegraphics{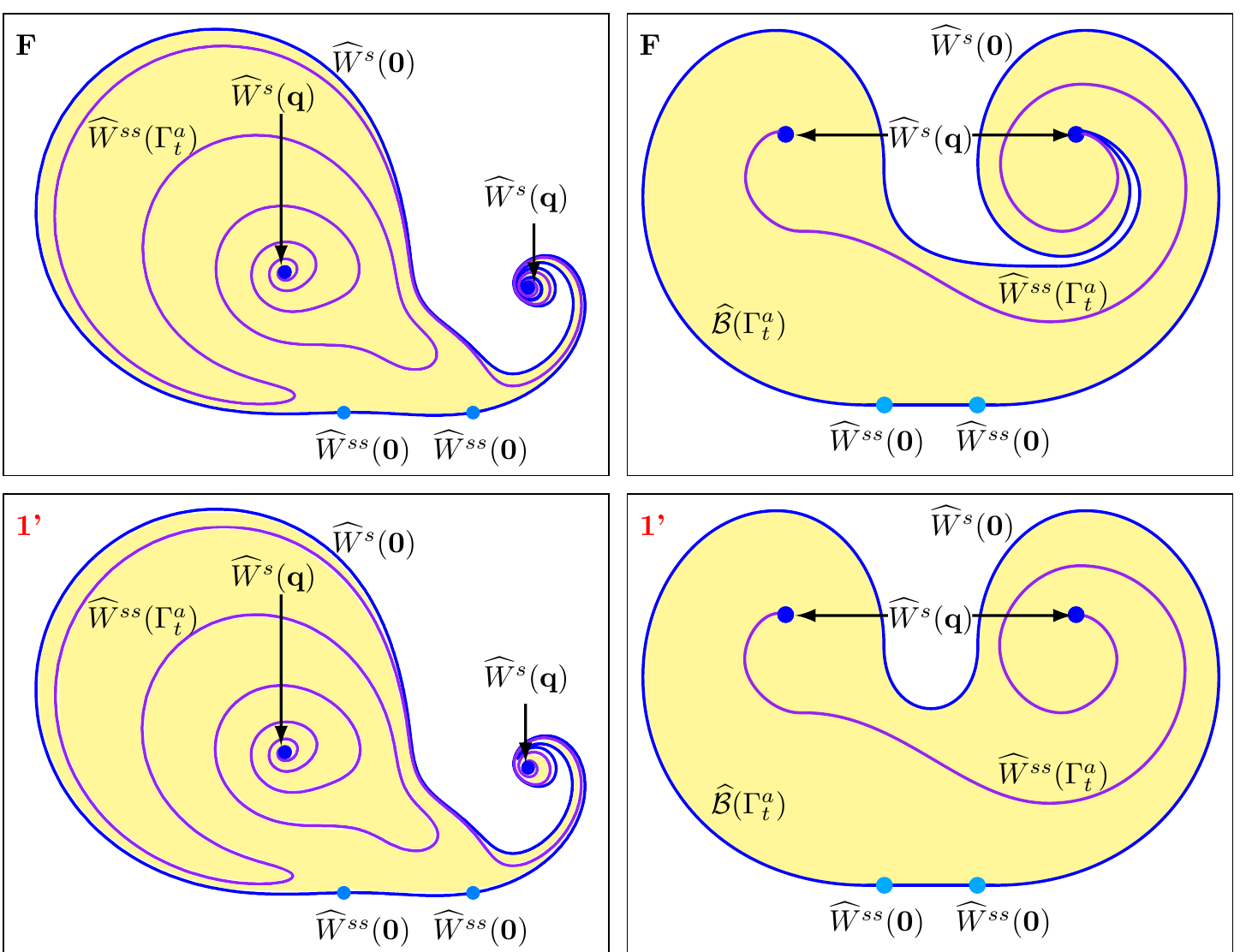}
\caption{Continued.}\label{fig:BDI4}
\end{figure}

In region \mRed{1} the intersection set $\widehat{W}^s(\textbf{0})$ (blue curve) on $\mS^*$ is a single
closed curve; due to the fact that  $W^s(\mathbf{0})$ is a topological
cylinder. The region enclosed by $\widehat{W}^s(\textbf{0})$ 
contains the two points of $\widehat{W}^s(\textbf{q})$ (dark blue).
The orientable attracting periodic orbit $\Gamma^{a}_o$ 
that exists in region \mRed{1} does not 
intersect $\mS^*$. Moreover, the intersection set
$\widehat{W}^{ss}(\Gamma_o^a)$ (purple curve) of its strong stable
manifold intersects $\mS^*$ in two closed
curves. The boundary of the basin 
$\mathcal{B}(\Gamma^a_o)$ is formed by $W^s(\mathbf{0})$ and $W^s(\mathbf{q})$;  hence, $\partial \widehat{\mathcal{B}}(\Gamma^{a}_o) =
\widehat{W}^s(\textbf{0}) \cup \widehat{W}^s(\textbf{q})$.  Note that
its closure $\overline{\widehat{\mathcal{B}}(\Gamma^{a}_o)}$ is topological
single closed disk.

Panel $\mathbf{H_o}$ of \cref{fig:BDI1} shows the homoclinic
bifurcation at the boundary between regions \mRed{1} and \mRed{2},  
where $\widehat{W}^s(\textbf{0})$ closes back on itself along
$\widehat{W}^{ss}(\textbf{0})$.  The basin 
$\widehat{\mathcal{B}}(\Gamma^{a}_o)$ is now disconnected and
$\overline{\widehat{\mathcal{B}}(\Gamma^a_o)}$ is topologically
equivalent to two disjoint disks. Furthermore, not all of
$\widehat{W}^s(\textbf{0})$ is part of $\partial
\widehat{\mathcal{B}}(\Gamma^{a}_o)$ any longer.

In region \mRed{2} the homoclinic orbit $\mathbf{\Gamma_{\rm hom}}$
becomes the orientable saddle periodic orbit $\Gamma_o$. Instead of
$\widehat{W}^s(\textbf{0})$, the intersection set
$\widehat{W}^s(\Gamma_o)$ now forms the outer part 
of the boundary set of  $\widehat{\mathcal{B}}(\Gamma^a_o)$, that is, $\partial
\widehat{\mathcal{B}}(\Gamma^a_o) = 
\widehat{W}^s(\Gamma_o)\cup \widehat{W}^s(\textbf{q})$. Note  that
$\widehat{W}^s(\textbf{0})$ accumulates on 
$\widehat{W}^s(\Gamma_o)$, which consists of two
topological circles, reflecting that $W^s(\Gamma_o)$ is also a
cylinder.  The accumulation of $\widehat{W}^s(\textbf{0})$ on
$\widehat{W}^s(\Gamma_o)$ is a
consequence of the $\lambda$-lemma; the
structurally stable heteroclinic orbit from $\Gamma_o$ to $\mathbf{0}$
forces $\widehat{W}^s(\textbf{0})$ to spiral around
$\widehat{W}^s(\Gamma_o)$. We remark that, as the $\lambda$-lemma  is local in
nature, this accumulation may be lost if a bigger sphere is chosen and $\widehat{W}^s(\Gamma_o)$ becomes tangent to the
sphere. 

The transition from region \mRed{2} to region \mRed{3} is via a
saddle-node bifurcation (\textbf{SNP}) of periodic orbits,  where
$\Gamma_o$ and $\Gamma_a$ merge and disappear. As a consequence,
$\widehat{W}^s(\Gamma_o)$, $\widehat{W}^{ss}(\Gamma^a_o)$ and $\widehat{\mathcal{B}}(\Gamma^a_o)$ are no
longer present in \cref{fig:BDI1} panel 
\mRed{3}. The intersection set $\widehat{W}^s(\textbf{0})$ now
accumulates on $\widehat{W}^s(\textbf{q})$, which reflects the
existence of a structurally stable heteroclinic orbit from \textbf{q} to \textbf{0}.

\subsubsection{Intersection sets in regions 4 and 5}
Panel $\mathbf{H_t}$ of \cref{fig:BDI1} is at the transition between
regions \mRed{3} and \mRed{4}, characterized by a codimension-one
non-orientable homoclinic orbit. As for the orientable homoclinic
orbit, shown in panel $\mathbf{H_o}$,  the
intersection set $\widehat{W}^s(\mathbf{0})$  connects back on itself at
$\widehat{W}^{ss}(\mathbf{0})$, but now 
$\widehat{W}^s(\mathbf{0})$ does not bound two 
open regions. Instead, two segments of $\widehat{W}^s(\mathbf{0})$ accumulate on
the intersection points $\widehat{W}^s(\mathbf{q})$, due to the persistence of the
heteroclinic orbit from $\mathbf{q}$ to $\mathbf{0}$.

In region \mRed{4},  the homoclinic orbit $\mathbf{\Gamma_{\rm hom}}$
becomes in the periodic orbit  $\Gamma_t$. Compare $\widehat{W}^s(\Gamma_o)$ in
panel \mRed{2} with $\widehat{W}^s(\Gamma_t)$ in panel
\mRed{4} of \cref{fig:BDI1}; for the former, $\widehat{W}^s(\Gamma_o)$
is composed of two closed curves, while for the latter, $W^s(\Gamma_t)$
intersects $\mS^*$ in a single curve that accumulates on
$\widehat{W}^s(\mathbf{q})$ as a consequence of the existence of a heteroclinic orbit from
$\mathbf{q}$ to $\Gamma_t$.  The intersection set
$\widehat{W}^s(\mathbf{0})$ consist of many curve segments; there is a
segment that accumulates on a single point  in
$\widehat{W}^s(\mathbf{q})$, while the other curve segments connect
the two intersection points $\widehat{W}^s(\mathbf{q})$. In
\cref{sec:Lambda}, we proved the existence of infinitely many
heteroclinic orbits in region \mRed{4}; as such, there must be
infinitely many curve segments of $\widehat{W}^s(\mathbf{0})$
accumulating  on $\widehat{W}^s(\mathbf{q})$. This is a consequence of
the $\lambda$-lemma when applied to the time-one map of the flow of
system~\cref{eq:san}. Since $\mS^*$ is transverse to $W^s(\mathbf{q})$,
each transverse heteroclinic orbit from   $\mathbf{q}$ to $\mathbf{0}$
creates at least one intersection curve 
$\widehat{W}^s(\mathbf{0})$ whose endpoints are 
$\widehat{W}^s(\mathbf{q})$. Furthermore, this set of curves accumulates
onto $\widehat{W}^s(\Gamma_t)$. In panel \mRed{4} of \cref{fig:BDI1}
we only show three of these infinitely many intersection curves of
$\widehat{W}^s(\mathbf{0})$; the existence of infinitely many curves
is indicated by  dashed blue curves; the three dots illustrates their accumulation on
$\widehat{W}^s(\Gamma_t)$. 

In region \mRed{5}, the period-doubling
bifurcation \textbf{PD} creates  $^2\Gamma_o$ and
$\Gamma^a_t$.  Note that $^2\Gamma_o$ is an orientable
periodic orbit, yet its intersection set
$\widehat{W}^{s}(^2\Gamma_o)$, composed of two open curves, is markedly different
from  $\widehat{W}^{s}(\Gamma_o)$ in region \mRed{2}. This is 
due to the existence of the two heteroclinic orbits from $\mathbf{q}$ to
$^2\Gamma_o$ that force the two curves in
$\widehat{W}^{s}(^2\Gamma_o)$ to accumulate on
$\widehat{W}^s(\mathbf{q})$; see panel \mRed{5} of
\cref{fig:BDI1}. The closure $\overline{\widehat{W}^{s}(^2\Gamma_o)}$
is a topological circle that bounds
$\widehat{\mathcal{B}}(\Gamma^{a}_t)$, namely,
$\partial\widehat{\mathcal{B}}(\Gamma^a_t) = \overline{\widehat{W}^{s}(^2\Gamma_o)}=
\widehat{W}^s(^2\Gamma_o)\cup
\widehat{W}^s(\textbf{q})$. Hence, the manifold $W^s(^2\Gamma_o)$,
together with $W^s(\mathbf{q})$, plays a similar  role as $W^s(\Gamma_o)$ in region \mRed{2}. The
set $\widehat{W}^s(\mathbf{0})$ does not change qualitatively
in the transition from region \mRed{4} to region \mRed{5}, in the
sense that all segments are in one-to-one correspondence with their
counterparts in region \mRed{4}.  The only difference is that
$\widehat{W}^{s}(\mathbf{0})$ now accumulates of
$\widehat{W}^{s}(^2\Gamma_o)$; more precisely, due to the
period-doubled nature of $^2\Gamma_o$, there are two sets of segments
in $\widehat{W}^s(\mathbf{0})$ that accumulate on  different curves of
$\widehat{W}^{s}(^2\Gamma_o)$.

\subsubsection{Intersection sets in regions 6 and 1'}
Panel $\mathbf{^2H_o}$ of \cref{fig:BDI1} 
shows how the intersection set $\widehat{W}^{s}(\mathbf{0})$ meets
itself transversally at $\widehat{W}^{ss}(\mathbf{0})$ 
in this bifurcation.  As $^2\Gamma_o$
becomes $\mathbf{^2\Gamma_{\rm hom}}$, the infinitely many curves of
$\widehat{W}^{s}(\mathbf{0})$ in region \mRed{4} disappear and only two curves
that connect $\widehat{W}^{ss}(\mathbf{0})$ and
$\widehat{W}^{s}(\mathbf{q})$ exist. In addition,  $\partial \widehat{\mathcal{B}}(\Gamma^a_t) \subset
\widehat{W}^s(\mathbf{0})\cup \widehat{W}^s(\textbf{q})$, that is, 
${W}^{s}(\mathbf{0})$ becomes the new separatrix in phase space.   Although  $\mathbf{H_o}$
and $\mathbf{^2H_o}$ are both codimension-one orientable 
homoclinic bifurcations, their intersection sets are not homeomorphic,
as seen in the respective panels of \cref{fig:BDI1}; there exists
a non-trivial intersection between
$W^u(\mathbf{q})$ and $W^s(\mathbf{0})$ at $^2\mathbf{H_o}$.

In region \mRed{6}, the intersection set $\widehat{W}^{s}(\mathbf{0})$ bounds
$\widehat{\mathcal{B}}(\Gamma^{a}_t)$, note that
$\overline{\widehat{\mathcal{B}}(\Gamma^{a}_t)}$ is a topological
annulus. The intersection set $\widehat{W}^{s}(\mathbf{0})$ is
composed of two disjoint 
curves that spiral into the intersection points
$\widehat{W}^{s}(\mathbf{q})$; hence the intersection set
$\overline{\widehat{\mathcal{B}}(\Gamma^{a}_t)}$ is not a simply
connected set, which indicates the persistence of the  
two heteroclinic orbits from $\mathbf{q}$ to $\mathbf{0}$.   

At the fold curve \textbf{F},  the unstable manifold $W^u(\textbf{q})$
intersects $W^s(\mathbf{0})$ tangentially in the heteroclinic orbit
$\gamma^*$; see panel \textbf{F} of \cref{fig:BDInc}.  On the level of
the intersection sets in panel \textbf{F} of \cref{fig:BDI1}, the set
$\widehat{W}^{s}(\mathbf{0})$ is formed by one segment that accumulates
on both sides on a single point in $\widehat{W}^s(\mathbf{q})$.  At this
bifurcation, $W^s(\mathbf{0})$ cannot cross $W^u(\mathbf{q})$, as
they are in tangency, and
$\overline{\widehat{\mathcal{B}}(\Gamma^{a}_t)}$ becomes a simply
connected set in region \mRed{6}. Hence, $\overline{\widehat{\mathcal{B}}(\Gamma^{a}_t)}$ is homeomorphic to a
closed disk, as is the case in region \mRed{1}.  

\begin{figure}
\centering
\includegraphics{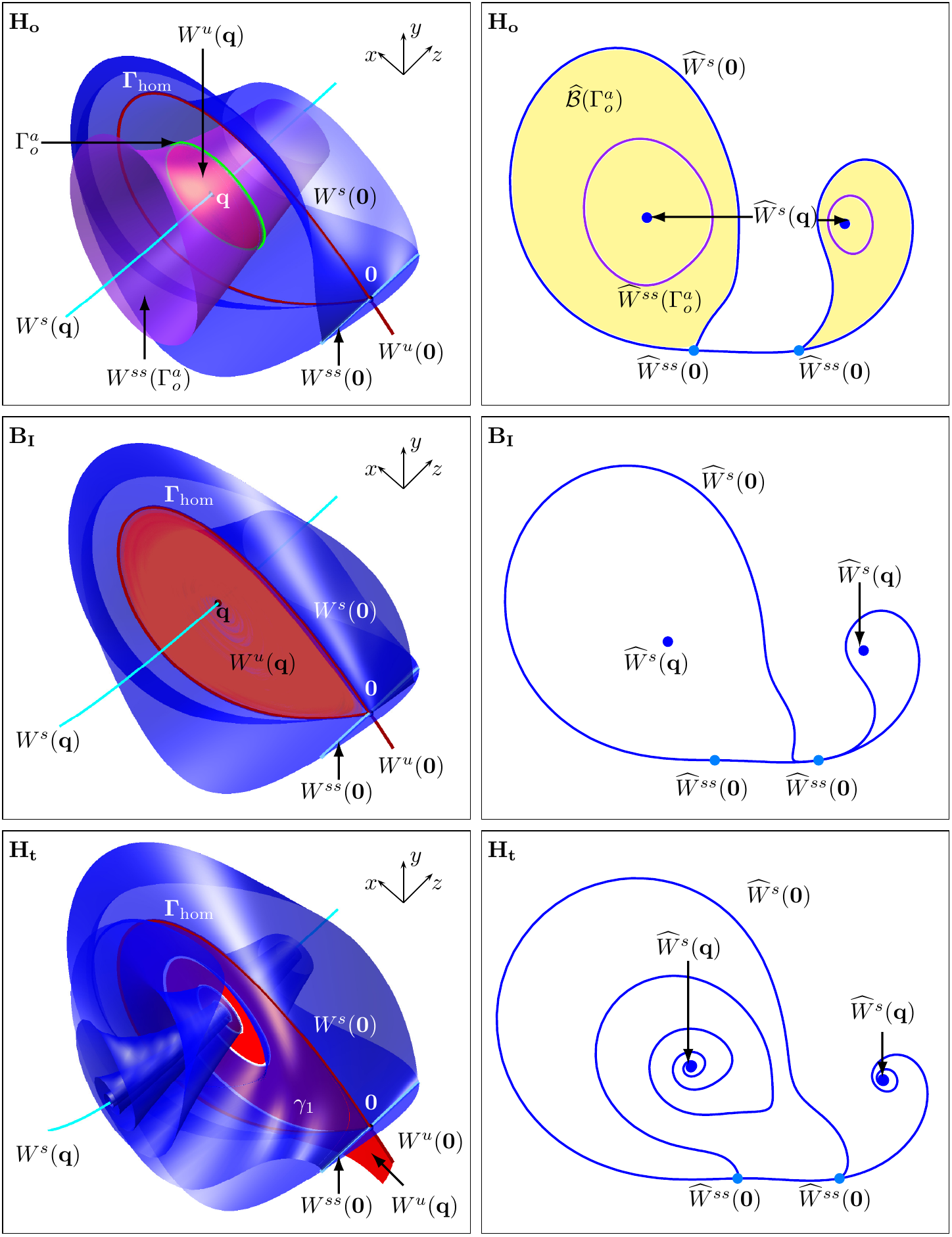}
\caption{Transition along the homoclinic bifurcation through the
  inclination flip bifurcation $\mathbf{B_I}$ of system~\cref{eq:san}. Shown are different manifolds in
  $\R^3$ (left column) and their respective
  stereographic projections (right column).  The color code is in
  \cref{fig:Sections} and the parameter values are given in 
  \cref{tab:Inc2}. See also the accompanying animation ({\color{red}
    GKO\_Bflip\_animatedFig9.gif}).}\label{fig:TranIn}  
\end{figure}

In region \mRed{1'}, $W^s(\mathbf{0})$ no longer intersects 
$W^u(\textbf{q})$ and this is reflected on the level of intersection
sets as a disconnection of $\widehat{W}^{s}(\mathbf{0})$ from both points in
$\widehat{W}^{s}(\mathbf{q})$.  Instead $\widehat{W}^{s}(\mathbf{0})$  encloses $\widehat{W}^{s}(\mathbf{q})$
and bounds $\widehat{\mathcal{B}}(\Gamma^{a}_t)$.  The only difference
with region \mRed{1} is the topological difference between
$\widehat{W}^{ss}(\Gamma^a_t)$ and $\widehat{W}^{ss}(\Gamma^a_o)$,
which are associated with a non-orientable and an orientable attracting
periodic orbit, respectively. As the transition through
$\mathbf{CC^-}$ involves the disappearance of
$\widehat{W}^{ss}(\Gamma^a_t)$, because there is no well-defined strong stable manifold in region~$\mRed{1^*}$,  the two topological
circles of $\widehat{W}^{ss}(\Gamma^a_o)$ appear only after crossing
$\mathbf{CC^+}$ into region \mRed{1}. During these transitions,
the other intersection curves and basin of attraction
do not change qualitatively.  

\subsection{Transition of the homoclinic orbit through the inclination flip}\label{sec:IncFlipTran}
We now focus specifically on the curve of homoclinic bifurcations and
illustrate the transition through the codimension-two 
homoclinic flip bifurcation point $\mathbf{B_I}$.  As illustrated in, e.g., \cite{OldKra1},
the two-dimensional manifold $W^s(\mathbf{0})$ can violate genericity
condition \textbf{(G3)} in two different ways, which depend on
the eigenvalues of the equilibrium; namely, whether
$|\lambda^{ss}|>2|\lambda^s|$ is fulfilled or not. Despite this difference, both
mechanisms unfold in the same way. Condition
$|\lambda^{ss}|<2|\lambda^s|$ was considered in \cite{Agu1} in the transition through the
inclination flip point $\mathbf{A_I}$ of case \textbf{A}. Here, we illustrate the transition
through the inclination flip point $\mathbf{B_I}$ of case \textbf{B}
for the case that $|\lambda^{ss}|>2|\lambda^s|$.

\cref{fig:TranIn} shows the transition through $\mathbf{B_I}$ on the
level of the invariant manifolds in the left
column, and their respective intersection sets with $\mS^*$ in the right
column. We show again the stereographic projections of the intersection
sets at $\mathbf{H_o}$ and $\mathbf{H_t}$ for comparison
purposes. In the accompanying animation
({\color{red} GKO\_Bflip\_animatedFig9.gif}) the phase portraits of \cref{fig:TranIn} are rotated clockwise
around the $y$-axis. At the codimension-one  
orientable homoclinic bifurcation in panel $\mathbf{H_o}$, 
the branch of $W^u(\mathbf{0})$ that spirals towards $\Gamma^a_o$ in 
region \mRed{1} now forms the homoclinic orbit $\mathbf{\Gamma_{\rm
    hom}}$, while the manifold $W^u(\mathbf{q})$ accumulates onto
$\Gamma^a_o$. Note that $\mathbf{\Gamma_{\rm hom}}$ returns to $\mathbf{0}$ along a
direction that is clearly transverse to $W^{ss}(\mathbf{0})$ (light-blue curve) and
$W^s(\mathbf{0})$ closes back on itself along
$W^{ss}(\mathbf{0})$. Furthermore, $W^{s}(\mathbf{0})$ is
topologically  a cylinder; compare 
with \cref{fig:ori}(a1). On the level of intersection sets,
$\widehat{W}^s(\mathbf{0})$ closes on $\widehat{W}^{ss}(\mathbf{0})$,
so that the basin of attraction $\widehat{\mathcal{B}}(\Gamma^a_o)$
is a disconnected set. At the
codimension-two point $\mathbf{B_I}$, the middle of \cref{fig:TranIn},
the surface $W^s(\mathbf{0})$ closes back on itself at 
$W^{ss}(\mathbf{0})$ in such a way that it
makes a quadratic tangency with itself at
$W^{ss}(\mathbf{0})$; the bottom panel of Fig.~2. in \cite{OldKra1} is
  misleading in this respect. Hence, if we follow the tangent plane
of $W^s(\mathbf{0})$ along $\mathbf{\Gamma_{\rm
    hom}}$ as $t \rightarrow -\infty$, it does not contain the strong stable
eigenvector of $\mathbf{0}$; this violates genericity condition
\textbf{(G3)}. As a result, $W^s(\mathbf{0})$ meets and closes along a
single branch of $W^{ss}(\mathbf{0})$. Additionally, the attracting periodic orbit $\Gamma^a_o$
is now the homoclinic orbit $\mathbf{\Gamma_{\rm hom}}$, making it
the boundary of $W^u(\mathbf{q})$ in phase space. On the level of
intersection sets, both parts of $\widehat{W}^s(\mathbf{0})$ have a
tangency with itself at only one of the intersection points of
$\widehat{W}^{ss}(\mathbf{0})$. Finally, at  the codimension-one 
non-orientable homoclinic bifurcation in panel $\mathbf{H_t}$, the stable manifold
$W^s(\mathbf{0})$ makes half a twist before closing along (both branches) of
$W^{ss}(\mathbf{0})$, so that the homoclinic orbit is non-orientable. An
interesting difference between $\mathbf{H_o}$
and $\mathbf{H_t}$ is the existence of the heteroclinic
orbit $\gamma_1$ in panel $\mathbf{H_t}$, caused by the 
transverse intersection of $W^s(\mathbf{0})$ and $W^u(\mathbf{q})$. In fact,
the long excursion of $\gamma_1$ around $\mathbf{q}$ becomes
$\mathbf{\Gamma_{\rm hom}}$ at $\mathbf{B_I}$. On
the level of the intersection sets, $\widehat{W}^s(\mathbf{0})$
consists of two curves that accumulates on
$\widehat{W}^s(\mathbf{q})$, as a consequence of the existence of
$\gamma_1$.

By looking at the stereographic projection in the right column of
\cref{fig:TranIn}, we can see a clearer difference between the  two
conditions. Condition $|\lambda^{ss}|<2|\lambda^s|$, as considered in
\cite{Agu1} for the case \textbf{A},  leads to a limit  at the moment of the inclination flip,
where one end of the intersection set $\widehat{W}^s(\mathbf{0})$
spirals into one of the points in $\widehat{W}^s(\mathbf{q})$; this is 
similar to the right segment in $\mathbf{H_t}$; the other end closes
back on $\widehat{W}^s(\mathbf{0})$, but along the weak direction of $W^s(\mathbf{0})$; see
Fig.~13 of \cite{Agu1}. In contrast, condition
$|\lambda^{ss}|>2|\lambda^s|$ as considered here for case \textbf{B},
leads to a limit at the moment of the 
inclination flip, at which the intersection set $\widehat{W}^s(\mathbf{0})$ closes
tangentially at only one of the intersection points of
$\widehat{W}^{ss}(\mathbf{0})$.

\section{Orbit flip of case \textbf{B}}\label{sec:Or}
A codimension-two orbit flip bifurcation occurs
when condition \textbf{(G2)} is 
violated, that is, the homoclinic orbit $\mathbf{\Gamma_{\rm hom}}$ is
a subset of the strong stable manifold $W^{ss}(\mathbf{0})$. Even
though the mechanism is different from that of the 
inclination flip, the orbit flip also results in a change from an orientable to a 
non-orientable codimension-one homoclinic bifurcation; moreover,
the theoretical unfoldings of both codimension-two points are the same
\cite{san4}.  We now demonstrate that  both bifurcations also have the same
topological organization on the level of the
manifolds involved. Here, we consider case \textbf{B}, meaning that, the
equilibrium $\mathbf{0}$ satisfies the eigenvalue conditions as given
in \cref{sec:hFB}. 

\begin{table}
\begin{center}
\begin{tabular}{|c|r@{.}l|r@{.}l|r@{.}l|r@{.}l|r@{.}l|}
\hline
\textbf{Homoclinic} & \multicolumn{2}{c|}{$\mathbf{H_o}$} &
\multicolumn{2}{c|}{$\mathbf{B_o}$} &\multicolumn{2}{c|}{$\mathbf{H_t}$}  &
\multicolumn{2}{c|}{$\mathbf{^2H_o}$}  & \multicolumn{2}{c|}{$\mathbf{F}$}\\ \hline
$\mu$ & -0 & 150000000 &  0 & 0 & 0 & 150000000 &  0 & 150000000 & 0 & 150000000\\ \hline
$\tilde{\mu}$ & -0 & 062331201 &  0 & 0 & 0 & 062381076& 0 & 069351963 & 0 & 070562587\\
\hline
\end{tabular}
\vspace{2mm}
\caption{Chosen representative parameter values at selected bifurcations in
  \cref{fig:BDO1}.} 
\label{tab:Orb2}
\end{center}
\end{table}

\begin{table}
\begin{center}
\begin{tabular}{|c|r@{.}l|r@{.}l|r@{.}l|r@{.}l|r@{.}l|r@{.}l|}
\hline
\textbf{Region} & \multicolumn{2}{c|}{\textbf{1}} &
\multicolumn{2}{c|}{\textbf{2}} & \multicolumn{2}{c|}{\textbf{3}} &
\multicolumn{2}{c|}{\textbf{4}} & \multicolumn{2}{c|}{\textbf{5}} &
\multicolumn{2}{c|}{\textbf{6}} \\ \hline
$\mu$ & $-0$ & 150 &  $-0$ & 150 & 0 & 150 & 0 & 150 & 0 & 150 & 0 & 150 \\ \hline
$\tilde{\mu}$ & $-0$ & 060 & $-0$ & 065 & 0 & 060 & 0 & 065 & 0 & 069 & 0 & 070 \\ 
\hline
\end{tabular}
\vspace{2mm}
\caption{Chosen representative parameter values for the different open
  regions in \cref{fig:BDO1}.} 
\label{tab:Orb1}
\end{center}
\end{table}

\begin{figure}
\centering
\includegraphics{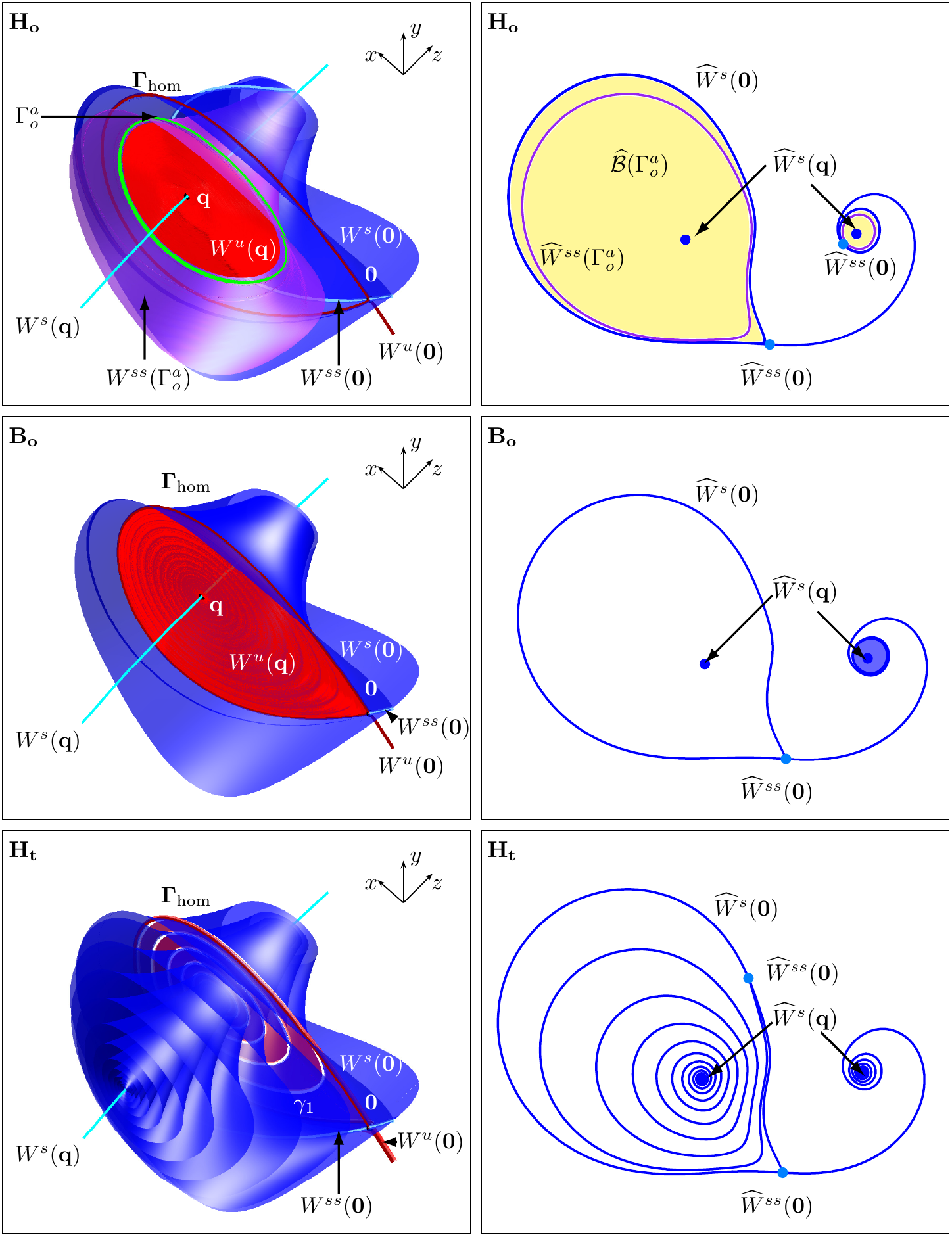}
\caption{ Transition along the homoclinic bifurcation through the orbit flip
bifurcation $\mathbf{B_o}$ of system
  \cref{eq:san}. Shown are different manifolds in
  $\R^3$ (left column) and their respective stereographic projections
  (right column).  The color code is in \cref{fig:Sections} and the
  parameter values are given in  
  \cref{tab:Orb2}.}\label{fig:TranOr} 
\end{figure}

We start by the transition of the homoclinic orbit through the orbit flip bifurcation
$\mathbf{B_o}$. \cref{fig:TranOr} shows the phase portraits and the intersection sets
with $\mS^*$ at the bifurcations $\mathbf{H_o}$, $\mathbf{B_o}$ and
$\mathbf{H_t}$  at the parameter  values given as in \cref{tab:Orb2}.   Note
that the panels $\mathbf{H_o}$ and $\mathbf{H_t}$ are topologically
equivalent to the respective panels in \cref{fig:TranIn} for the
inclination flip, but panel $\mathbf{B_o}$ is different.  At the
moment of the orbit flip, the one-dimensional strong stable manifold
$W^{ss}(\mathbf{0})$ intersects $\mS^*$ in a single point, because the
other branch of $W^{ss}(\mathbf{0})$ is  $\mathbf{\Gamma_{\rm hom}}$.
Hence, only one end of $\widehat{W}^s(\mathbf{0})$ closes back on
itself.  The other end spirals into one point of
$\widehat{W}^s(\mathbf{q})$, but at an algebraic rather than an
exponential rate; we indicate
this accumulation by a light-blue shading. Note that the relative position of the points
in $\widehat{W}^{ss}(\mathbf{0})$ swaps before and after the orbit flip;
see panels $\mathbf{H_o}$ and $\mathbf{H_t}$ in
\cref{fig:TranOr}. Unlike the case $\mathbf{IF}$, the case $\mathbf{OF}$ 
does not have multiple ways of breaking condition \textbf{(G2)} that
depends on additional eigenvalue conditions. Furthermore, the results that we find for the
transition for the orbit flip of type \textbf{B} are topologically equivalent to those
found for case \textbf{A} \cite{Agu1}.  
\begin{figure}
\centering
\includegraphics{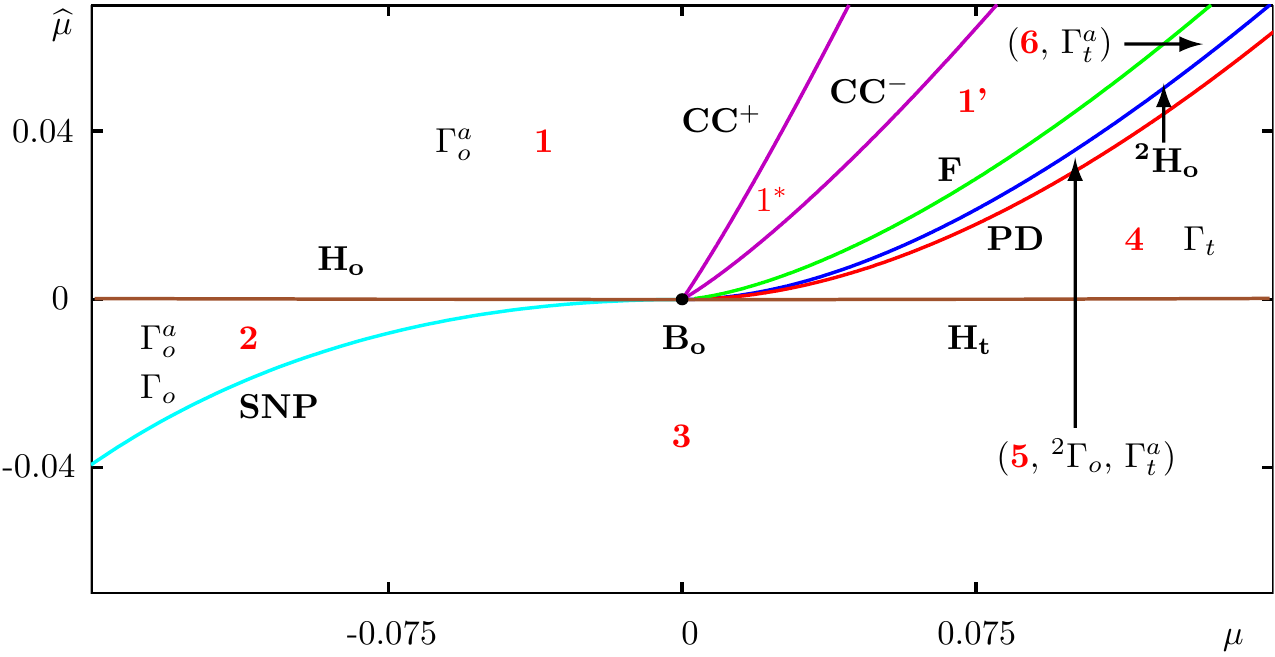}
\caption{Bifurcation diagram in the $(\mu,\widehat{\mu})$-plane, where
  $\widehat{\mu} =  10(\tilde{\mu}- 0.4157\mu)$, near an orbit flip
    bifurcation  $\mathbf{B_o}$ of system~\cref{eq:san} for other
    parameters as given in \cref{sec:san}.  The color  code and nomenclature of the  
  regions is the same as  given in \cref{fig:BDInc}.}\label{fig:BDO1}
\end{figure}
\begin{figure}[h]
\centering
\includegraphics{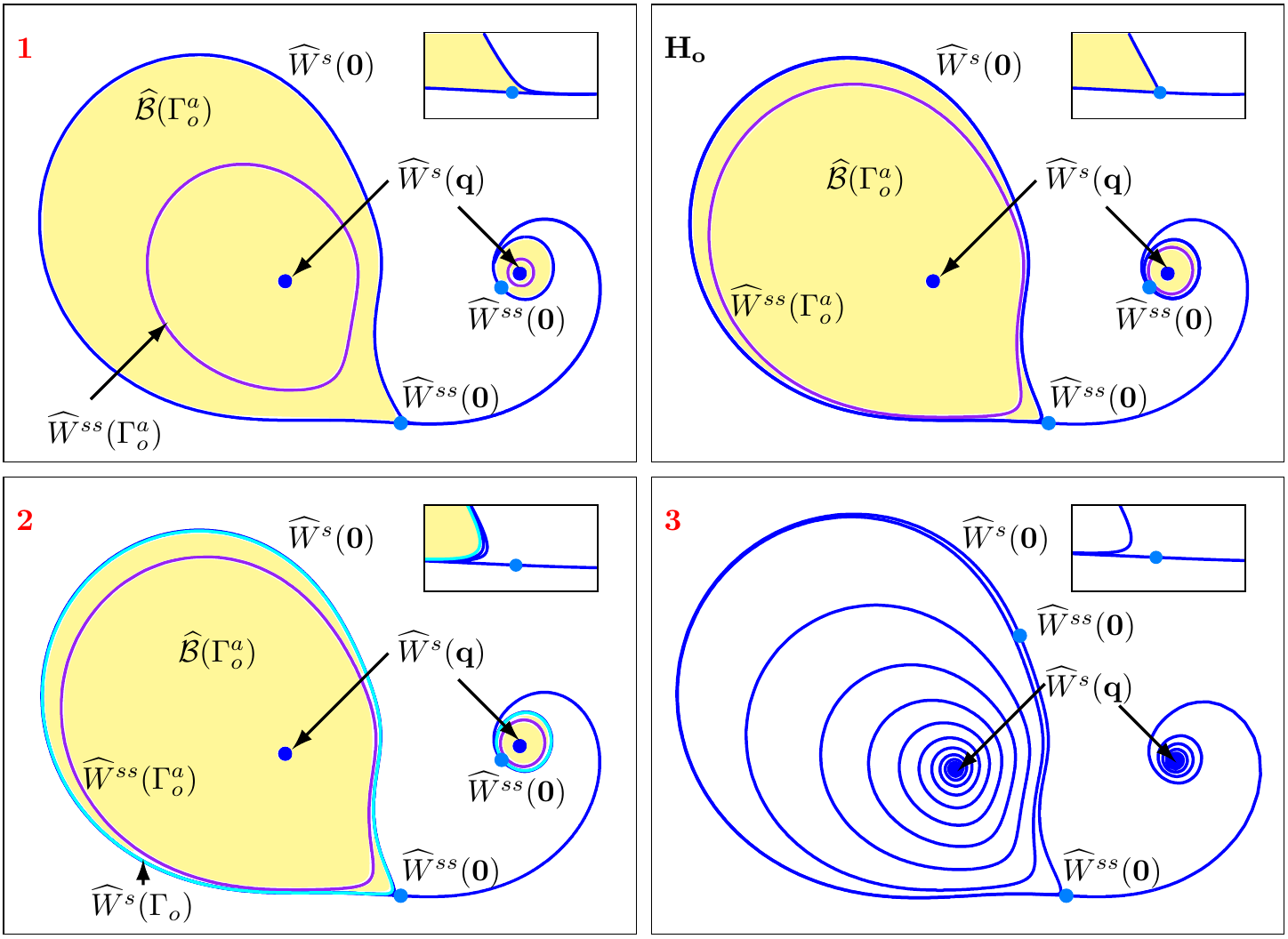}
\caption{Stereographic projections of the  
  intersection sets of the invariant manifolds with $\mS^*$ in the
  regions and at bifurcations of the bifurcation diagram in
  \cref{fig:BDO1} near the orbit flip $\mathbf{B_o}$.  The insets show enlargements
  around one of the points of $\widehat{W}^{ss}(\mathbf{0})$. The
  color  code and nomenclature of the  
  regions is the same as  given in \cref{fig:BDInc} and
  \cref{fig:BDI1}.  For respective parameter values see 
  \cref{tab:Orb1} and \cref{tab:Orb2}.}\label{fig:BDOS}
\end{figure}
\setcounter{figure}{11}
\begin{figure}
\centering
\includegraphics{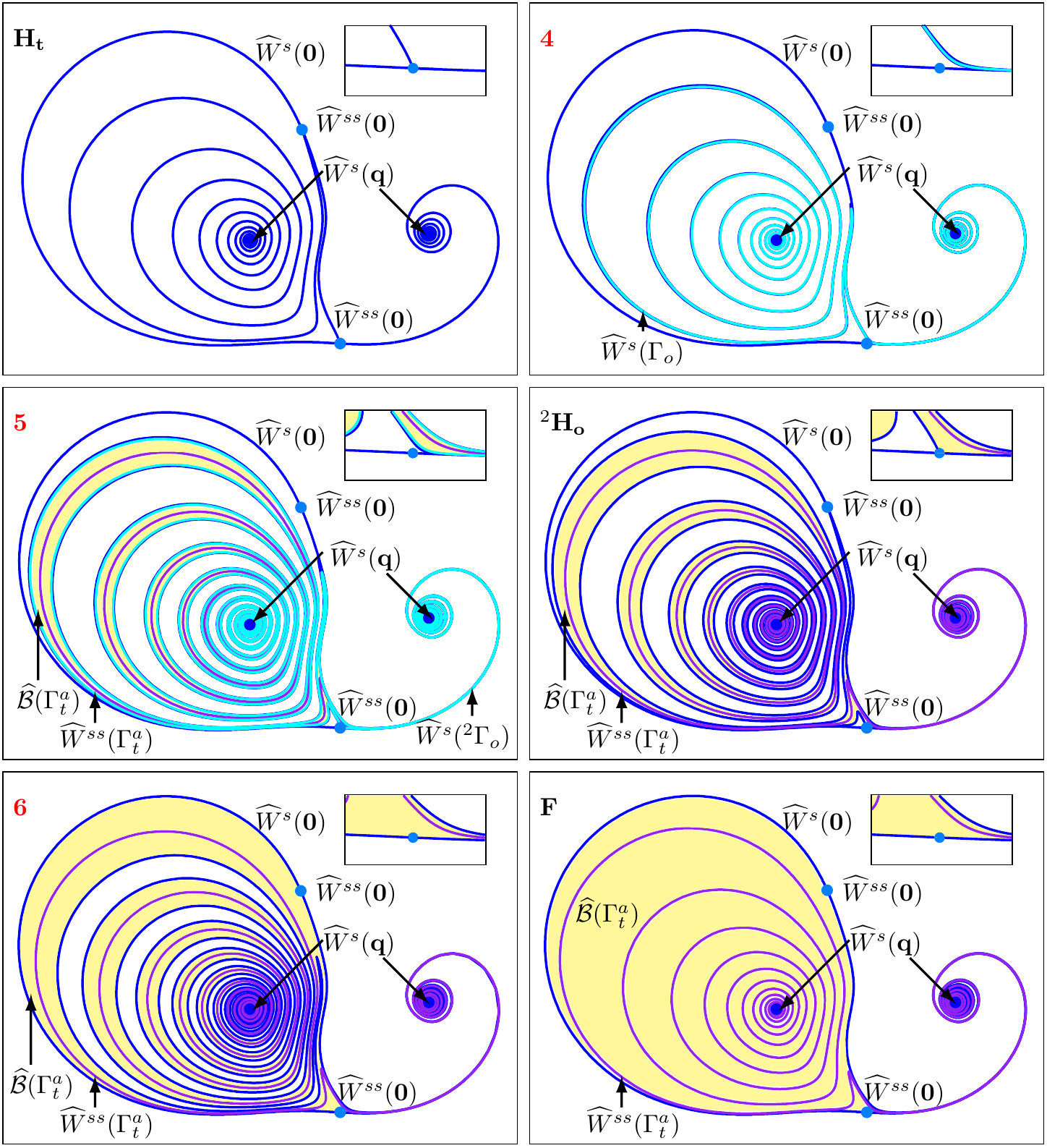}
\caption{Continued.}\label{fig:BDO2}
\end{figure}

We now present the unfolding of the orbit flip bifurcation $\mathbf{B_o}$
with respect to the parameters $\mu$ and $\tilde{\mu}$. \cref{fig:BDO1} shows the
bifurcation diagram locally near the codimension-two orbit flip point $\mathbf{B_o}$ in the
$(\mu,\widehat{\mu})$-plane; here we use the coordinate transformation $\widehat{\mu} := 
10(\tilde{\mu}- 0.4157\mu)$ to improve the
visualization.  Close to $\mathbf{B_o}$, the bifurcation diagram  is
topologically equivalent to the bifurcation diagram of the inclination
flip; see \cref{fig:BDInc}. In particular,  we also find the fold \textbf{F} of
heteroclinic orbits, and the curves $\mathbf{CC^+}$ and $\mathbf{CC^-}$
in the same relative positions with respect to the other
bifurcation curves. As we did for the inclination flip in
\cref{sec:Inc}, we use the bifurcation diagram in 
\cref{fig:BDO1} as a reference to describe the re-organization of the
global manifolds in phase space for system~\cref{eq:san} with
parameters as given in \cref{sec:san}.  Representative points
from each region we selected in the
$(\mu,\tilde{\mu})$-plane are listed in \cref{tab:Orb1}. 

For each point we compute the respective intersection sets with
$\mS^*$, where we also illustrate the manifolds for
points on the codimension-one bifurcation curves 
$\mathbf{^2H_o}$ and \textbf{F}. The parameter values $\mu$ and $\tilde{\mu}$
for at the bifurcation points are given in \cref{tab:Orb2}. \cref{fig:BDOS} shows the
selected stereographic projections of the intersection 
sets of the stable manifolds with $\mS^*$ at all ten representative
points. The insets
are enlargements illustrating the behaviour of
$\widehat{W}^s(\mathbf{0})$ close to one of the points of
$\widehat{W}^{ss}(\mathbf{0})$. As can be checked readily, the stereographic
projections in each panel are  topologically equivalent to the respective topological sketches in 
\cref{fig:BDI1} for the inclination flip. This means that the 
interactions between the manifolds in phase
space for the inclination flip are qualitatively the same as for the
orbit flip case. In particular, we have 
the same conclusions about the number of heteroclinic orbits between equilibria and
for saddle periodic orbits in regions \mRed{4} and \mRed{5}. Moreover,
the role of the separatrix in phase space of the basin of the attracting periodic orbit
switches between the stable manifolds of 
the origin and saddle periodic orbits in precisely the same way.

\section{Discussion} \label{sec:Dis}
We presented a study of invariant manifolds of equilibria
and saddle periodic orbits in the vicinity of a codimension-two homoclinic flip
bifurcation of case \textbf{B}. We characterized the regions
with different equilibria configurations for system~\cref{eq:san} for both inclination and orbit flip
bifurcations of case \textbf{B} by means of Poincar\'e
compactification \cite{Dum1,Vel} and focussed on the parameter region
for which there is only one additional saddle focus
equilibrium \textbf{q}.  We computed the unfoldings of both
inclination and orbit flip bifurcation points in two-parameter
planes and presented representative phase portraits. In this way, we
illustrated the role of the two-dimensional stable 
manifold  $W^s(\mathbf{0})$ of the real saddle equilibrium at the origin
$\mathbf{0}$ and its interaction with other manifolds for the overall
organization of phase space in the vicinity
of the codimension-two bifurcations; in particular, this study included the invariants
manifolds of $\mathbf{q}$, which lies outside the tubular neighborhood
of the homoclinic orbit. Similar to case \textbf{A}, presented in
\cite{Agu1}, we found a  fold \textbf{F}  of heteroclinic orbits from \textbf{q} to \textbf{0} for
case \textbf{B} of the homoclinic blip bifurcation. Furthermore, the presence of saddle periodic
orbits in case \textbf{B} has implications for the interaction of the manifolds of $\mathbf{0}$ and
$\mathbf{q}$: in certain parameter regions, there exist infinitely
many heteroclinic orbits from \textbf{q} to \textbf{0}  close to a
homoclinic flip bifurcation; this phenomenon does not occur 
for case \textbf{A}. Note that 
these heteroclinic orbits are distinguished by their
numbers of large excursions; in particular, large excursions in
periodic orbits can be identified with spiking behaviours as studied
in the Hindmarsh-Rose model that describes the essential
spiking behaviour of a neuron \cite{Lina1}.  

Our approach was to compute $W^s(\mathbf{0})$ as a global
object in phase space to study how it re-arranges itself as the system undergoes 
different bifurcations. Moreover,  we determined the
two-dimensional stable and unstable manifolds of the saddle periodic
orbits that co-exist in certain regions of parameters space
and studied their interaction with $W^s(\mathbf{0})$. We also computed the intersection sets of the
stable manifolds with a suitable sphere $\mS^*$, chosen such that it
contains all compact invariant objects close to $\mathbf{0}$. In particular,
knowledge of the intersection sets on $\mS^*$ 
allowed us to clarify the properties of basins of
attraction. Our numerical results confirm that the local two-parameter
unfoldings of both inclination and orbit flip bifurcations of case
\textbf{B} are the same, even on the level of the interacting global
manifolds including those of $\mathbf{q}$;  the only
difference lies in the phase portraits at the codimension-two points $\mathbf{B_o}$
and $\mathbf{B_I}$.  

Our findings can be summarized as follows:

\textbf{Results} (manifold structure near flip bifurcation of case
\textbf{B}). 
\begin{itshape}
Consider system~\cref{eq:san} near an 
  inclination flip or an orbit flip homoclinic bifurcation of case
  $\mathbf{B}$ at the origin, such that there also exists a nearby
  unstable saddle-focus $\mathbf{q}$. For an inclination and orbit flip the  bifurcation diagram is
  topologically equivalent to the ones shown in \cref{fig:BDInc} and
  \cref{fig:BDO1}; where regions and bifurcations are labelled
  according to \cref{sec:Inc}. The configurations of the manifolds in phase space and their
intersection sets with $\mS^*$ are as follows,
\begin{itemize}
\item[$\mathbf{B_I}$] At the codimension-two inclination flip point
  $\mathbf{B_I}$, the intersection set $\widehat{W}^s(\mathbf{0})$ on $\mS^*$ is a
  closed curve tangent to itself at one of the two points in
  $\widehat{W}^{ss}(\mathbf{0})$.  The stable manifold
  $W^s(\mathbf{0})$ in $\R^3$  closes back on itself along a single branch of
  $W^{ss}(\mathbf{0})$.  More precisely, we find that $W^s(\mathbf{0})$ has a
  quadratic tangency with itself at $W^{ss}(\mathbf{0})$. The unstable manifold $W^u(\mathbf{q})$
  accumulates on $\mathbf{\Gamma_{\rm hom}}$.
\item[$\mathbf{B_o}$] At the codimension-two orbit flip point
  $\mathbf{B_o}$, the intersection set $\widehat{W}^s(\mathbf{0})$
  closes on itself at only one point of
  $\widehat{W}^{ss}(\mathbf{0})$, because the second intersection point
  $\widehat{W}^{ss}(\mathbf{0})$ becomes the homoclinic orbit and does
  not intersect $\mS^*$.  The segment $\widehat{W}^s(\mathbf{0})$ on
  the other side of $\widehat{W}^{ss}(\mathbf{0})$ accumulates on 
  $\widehat{W}^s(\mathbf{q})$.  The homoclinic orbit
  bounds the two-dimensional manifold $W^u(\mathbf{q})$ and part
  of $W^s(\mathbf{0})$ accumulates on $W^s(\mathbf{q})$.
\item[$\mathbf{1}$]  In region \mRed{1}  the intersection set $\widehat{W}^s(\mathbf{0})$ of the
  stable manifold of $\mathbf{0}$  encloses the basin of attraction
  $\widehat{\mathcal{B}}(\Gamma^a_o)$ of the orientable attracting 
  periodic orbit $\Gamma^a_o$; here, $\widehat{\mathcal{B}}(\Gamma^a_o)$ is a
  connected set and its closure  is homeomorphic to a
  disk, and $\widehat{W}^{ss}(\Gamma^a_o)$ is the union of two
  topological circles. The stable manifold $W^s(\mathbf{0})$ is the boundary of
  the  basin of attraction of $\Gamma^a_o$, and
  the unstable manifold $W^u(\mathbf{q})$ is bounded by the attracting
  periodic orbit $\Gamma^a_o$.
\item[$\mathbf{H_o}$] Along the orientable homoclinic curve $\mathbf{H_o}$,
  the intersection set $\widehat{W}^s(\mathbf{0})$ closes on itself at
  $\widehat{W}^{ss}(\mathbf{0})$ and encloses the region
  $\widehat{\mathcal{B}}(\Gamma^a_o)$.  The closure of 
  $\widehat{\mathcal{B}}(\Gamma^a_o)$ is now homeomorphic to two
  disks.  The stable manifold
  ${W}^s(\mathbf{0})$ closes along $W^{ss}(\mathbf{0})$ and
  creates a homoclinic orbit $\mathbf{\Gamma_{\rm hom}}$ that forms
  the boundary of $W^u(\mathbf{q})$.
\item[$\mathbf{2}$] In region \mRed{2} there exists a saddle periodic
  orbit $\Gamma_o$.  The intersection set $\widehat{W}^s(\mathbf{0})$
  spirals towards $\widehat{W}^s(\Gamma_o)$. The 
  closure of $\widehat{\mathcal{B}}(\Gamma^a_o)$ is
  homeomorphic to two disk that are each bounded by a topological
  circle in $\widehat{W}^s(\Gamma_o)$. The stable manifold
  $W^s(\Gamma_o)$ is the boundary of $\mathcal{B}(\Gamma^a_o)$, and
  $W^u(\mathbf{q})$ is contained in $\mathcal{B}(\Gamma^a_o)$ and
  accumulates on $\Gamma^a_o$.  Furthermore, $W^u(\Gamma_o)$ intersects
  $W^s(\mathbf{0})$ in a structurally  stable heteroclinic orbit.
\item[$\mathbf{SNP}$]  At the curve of saddle-node of periodic
  orbit $\mathbf{SNP}$, the periodic orbits $\Gamma^a_o$ and
  $\Gamma_o$ merge into a non-hyperbolic periodic orbit and
  disappear in region~\mRed{3}.
\item[$\mathbf{3}$] In region \mRed{3}, the segment of the
  intersection set $\widehat{W}^s(\mathbf{0})$ 
  spirals towards $\widehat{W}^s(\mathbf{q})$.  There exists
  a structurally stable heteroclinic orbit $\gamma_1$ from $\mathbf{q}$
  to $\mathbf{0}$, and $W^u(\mathbf{q})$ is bounded by the unstable
  manifold $W^u(\mathbf{0})$. 
\item[$\mathbf{H_t}$] Along the non-orientable homoclinic
  curve $\mathbf{H_t}$, the homoclinic orbit
  $\mathbf{\Gamma_{\rm hom}}$  exists, the intersection set
  $\widehat{W}^s(\mathbf{0})$ closes on itself  at
  $\widehat{W}^{ss}(\mathbf{0})$, and segments of it accumulate  on
  $\widehat{W}^s(\mathbf{q})$. The stable manifold ${W}^s(\mathbf{0})$
  closes along $W^{ss}(\mathbf{0})$, while it intersects $W^u(\mathbf{q})$ transversally.
\item[$\mathbf{4}$]  In region \mRed{4} there exists the periodic
  orbit $\Gamma_t$. The intersection set $\widehat{W}^s(\mathbf{0})$
  consists  of infinitely many curves that spiral towards $\widehat{W}^s(\mathbf{q})$
  and  accumulate on $\widehat{W}^s(\Gamma_t)$, which also
  spirals towards $\widehat{W}^s(\mathbf{q})$.  There exists one structurally
  stable heteroclinic orbit from $\Gamma_t$ to $\mathbf{0}$ and one
  from $\mathbf{q}$ to $\Gamma_t$. Furthermore, there are infinitely many structurally stable
  heteroclinic orbits from $\mathbf{q}$ to  $\mathbf{0}$. The unstable
  manifold $W^u(\mathbf{0})$ bounds both $W^u(\Gamma_t)$ and
  $W^u(\mathbf{q})$. 
\item[$\mathbf{PD}$]  Along the period-doubling bifurcation curve
  $\mathbf{PD}$, the periodic orbit $\Gamma_t$ is
  non-hyperbolic. It turns into an attracting periodic orbit
  $\Gamma^a_t$ and creates the period-doubled periodic orbit  $^2\Gamma_o$  in region~\mRed{5}. 
\item[$\mathbf{5}$] In region \mRed{5}, the intersection set $\widehat{W}^s(^2\Gamma_o)$ consists
  of two curves that spiral towards $\widehat{W}^s(\mathbf{q})$ and enclose
   $\widehat{\mathcal{B}}(\Gamma^a_t)$. In a neighborhood of these curves, there are
   infinitely many curves $\widehat{W}^s(\mathbf{0})$ that spiral towards
   $\widehat{W}^s(\mathbf{q})$.  Furthermore,
   $\widehat{W}^{ss}(\Gamma^a_t)$ also spirals towards
   $\widehat{W}^s(\mathbf{q})$. The stable manifold $W^s(^2\Gamma_o)$  is
   the boundary of the  basin of attraction of
   $\Gamma^a_t$. There exist a structurally stable
   heteroclinic orbit from $^2\Gamma_o$ to $\mathbf{0}$ and two from
   $\mathbf{q}$ to $^2\Gamma_o$.  Also, there are
   infinitely many structurally stable heteroclinic orbits from
   $\mathbf{q}$ to $\mathbf{0}$. Moreover, the part of $W^u(\mathbf{q})$
   bounded by the two heteroclinic orbits from $\mathbf{q}$ to
   $^2\Gamma_o$ accumulates on $\Gamma^a_t$, while the other
   part is bounded by $W^u(\mathbf{0})$.
\item[$\mathbf{^2H_o}$] Along the curve $\mathbf{^2H_o}$ the periodic
  orbit  $^2\Gamma_o$ disappears and the homoclinic orbit
  $^2\mathbf{\Gamma_{\rm hom}}$ is created. The
  intersection set $\widehat{W}^s(\mathbf{0})$ consists of curves
  that close along  $\widehat{W}^{ss}(\mathbf{0})$ or spiral
  towards $\widehat{W}^s(\mathbf{q})$. In the process, 
  $\widehat{W}^s(\mathbf{0})$ encloses
  $\widehat{\mathcal{B}}(\Gamma^a_t)$. The stable manifold
  $W^s(\mathbf{0})$ is the boundary of $\mathcal{B}(\Gamma^a_t)$. Furthermore,
  infinitely many heteroclinic orbits from $\mathbf{q}$ to
  $\mathbf{0}$ disappear at once, and only two are preserved.
\item[$\mathbf{6}$] In region \mRed{6} the intersection set $\widehat{W}^s(\mathbf{0})$ 
   forms the boundary of   $\widehat{\mathcal{B}}(\Gamma^a_t)$ and the
   closure of their union is homeomorphic to an annulus.
\item[$\mathbf{F}$] Along the fold $\mathbf{F}$ of heteroclinic orbits, the intersection set
  $\widehat{W}^s(\mathbf{0})$ encloses 
  $\widehat{\mathcal{B}}(\Gamma^a_t)$ but only one curve of
  $\widehat{W}^s(\mathbf{0})$ goes to $\widehat{W}^s(\mathbf{q})$.
  Moreover, the closure of $\widehat{\mathcal{B}}(\Gamma^a_o)$ is
  again a topological disk. The stable manifold $W^s(\mathbf{0})$
  is tangent to $W^u(\mathbf{q})$ at the heteroclinic orbit
  $\gamma^*$; this tangency is quadratic.
\item[$\mathbf{1'}$] In region \mRed{1'}, the intersection set
  $\widehat{W}^s(\mathbf{0})$  no longer accumulates on
  $\widehat{W}^s(\mathbf{q})$; the situation is topologically equivalent to that in region
  \mRed{1} except that $\widehat{W}^{ss}(\Gamma^a_t)$ accumulates on
  $\widehat{W}^s(\mathbf{q})$.  In  phase space, $W^u(\mathbf{q})$
  accumulates on $\Gamma^a_t$, and the heteroclinic orbits 
  between $\mathbf{q}$ and $\mathbf{0}$ have disappeared. 
\item[$\mathbf{CC^-}$] At the curve $\mathbf{CC^-}$, the nontrivial Floquet
  multipliers of $\Gamma^a$ are both the same negative real number,
  meaning that the periodic   orbit $\Gamma^a_t$ becomes $\Gamma^a$. There
  does not exist a well-defined strong 
  stable manifold $W^{ss}(\Gamma^a)$.
\item[$\mathbf{1^*}$] in region~$\mRed{1^*}$ the Floquet multipliers
  of $\Gamma^a$ are complex conjugate and their real part becomes
  positive when approaching the curve $\mathbf{CC^+}$
\item[$\mathbf{CC^+}$] At the curve $\mathbf{CC^+}$, the
  nontrivial Floquet   multipliers of $\Gamma^a$ are both the same
  positive number. In the transition to region~\mRed{1}, the periodic orbit 
  $\Gamma^a$ becomes $\Gamma^a_o$ and there exists a well-defined
  strong stable manifold $W^{ss}(\Gamma^a_o)$.
\end{itemize}
\end{itshape}

As discussed before, the existence of $\mathbf{q}$ induces new
phenomena in the unfolding of an homoclinic flip bifurcation, even
though it does not lie in a tubular neighbourhood of the homoclinic
orbit. It is worth noting that we found parameter regimes of 
system~\cref{eq:san} with none or several additional equilibria. Of
particular interest is the situation where no additional equilibria
exist; since $W^s(\mathbf{q})$ plays an important role in the overall
organization of the two-dimensional global manifolds, we conjecture that
a one-dimensional manifold from infinity then takes on the
role of $W^s(\mathbf{q})$. The compactified version of system~\cref{eq:san} should help
with answering this question. 

In ongoing work we intend to understand the nature
of the global manifolds close to  the most challenging case of a
homoclinic flip bifurcation of case
\textbf{C}.  Its unfolding features  
infinitely many codimension-one homoclinic bifurcations and period-doubling
cascades; this creates horseshoe-regions in the parameter
plane that  are bounded by tangencies of different manifolds. 

\appendix
\section{Poincar\'e Compactification}\label{secApp:Com}
\begin{figure}
\centering
\includegraphics{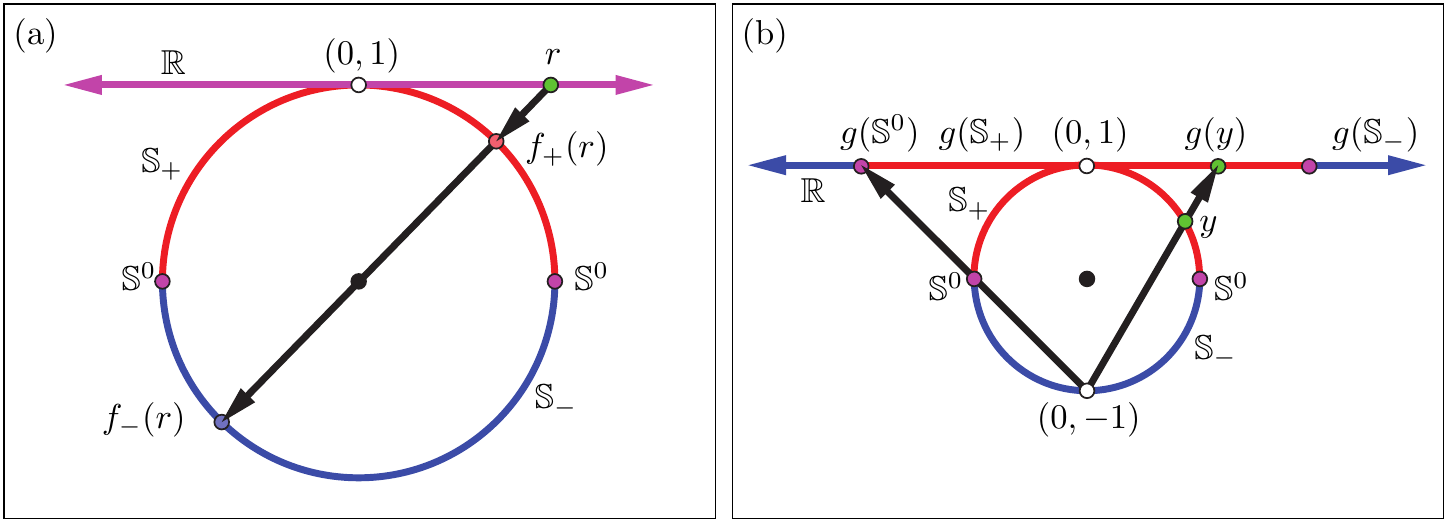}
\caption{Sketches of the transformations for Poincar\'e
  compactification of $\R$.  Panel (a) shows
  how the central projection sends $\R$ (purple curve) to the upper (red)
  and the lower (blue) hemispheres of the circle,  denoted $\mS_+$
  and $\mS_-$, respectively.  In panel (b),  the stereographic projection
  from the south pole is applied to send $\mS_{\pm}$  back to
  the intervals on the real line, as indicated by the corresponding colors.}\label{fig:IdCom} 
\end{figure}

To describe Poincar\'e compactification \cite{Dum1,Vel} of the
three-dimensional vector field \cref{eq:san}.  We first describe 
the Poincar\'e compactification for a one-dimensional system on $\R$. As
is illustrated in \cref{fig:IdCom}(a), the one-dimensional 
phase space (purple curve) is identified with the tangent space of the
one-dimensional sphere, the circle $\mS \subset \R^2$, at its north pole $(0,1)\in\R^2$. Each point $r
\in \R$ (green dot) is related via inverse central projections $f_{\pm}$ to antipodal points
$f_{\pm}(r) \in \mS$, one on the upper half sphere $\mS_+$ (red dot)
and one on the lower half sphere $\mS_-$ (blue dot). In a second step,
shown in \cref{fig:IdCom}(b), The south-pole projection $g$ is used to map the 
northern hemisphere $\mS_+$ to the interval $(-2,2)$ and the equator
$\mS^0$ to its boundary $\{-2,2\}$. Note that $g$ maps $\mS_-
\setminus \{(0,-1)\}$ to the two open intervals  $(-\infty,-2)$ and 
$(2,\infty)$, which constitutes a second transformation of $\R$ that
is not compact. For our purposes, it makes sense to work with $\mS_+$ only.

 To understand how a vector field $X$ defined in $\R$ is transformed
 by the maps $f_{\pm}$ we refer to the following commutative diagram:
\begin{equation*}
\begin{CD}
\R @>f_{\pm}>> \mS_{\pm}\\
@VXVV @VVX_{\mS_{\pm}}V\\
\text{T}\R=\R @>Df_{\pm}>> \text{T}\mS_{\pm}
\end{CD}
\end{equation*}
Here, $\text{T}\R$ and $\text{T}\mS_{\pm}$ are the respective tangent bundles of
$\R$ and $\mS_{\pm}$, the map $Df_{\pm}$ is the Jacobian of
$f_{\pm}$, and the two vector fields $X_{\mS_{+}}$ and $X_{\mS_{-}}$
defined on $\mS_{+}$ and $\mS_{-}$ are conjugate to $X$, respectively. Consequently,
we can write the transformed vector field $X_{\mS_{\pm}}$ as
$X_{\mS_{\pm}}= Df_{\pm} \circ X \circ f^{-1}_{\pm}$. It is possible
to extend the domain of definition of $X_{\mS_{\pm}}$ to the whole of
$\mS$ provided $X$ is polynomial; when doing this, the dynamics at the
equator of $\mS$ are invariant and represent the dynamics at infinity
of $X$ . Let $\rho(X)$ be the extension of $X_{\mS_{\pm}}$ in $\mS$. We have that $g$ generates
the vector field $\overline{X}=Dg \circ \rho(X) \circ g^{-1}$  on $\R$,
where the flow on $(-2,2)$ is conjugate to the flow of $X$ and
$\{-2, 2\}$  represents the dynamics at infinity. This procedure can be generalized to
higher dimensions, the best-known case is the compactification to the
Poincar\'e-disk of polynomial vector fields on $\R^2$
\cite{Dum1,Vel}. We are interested here in $\R^3$, which has been
studied for certain models \cite{Jau1,Jes1, Jes2,Raz1} through the
use of coordinate charts.

\subsection{Compactification of  $\R^3$} \label{sec:compofR3} 
Recall that system~\cref{eq:san} is defined as the polynomial vector
field
\begin{equation*}
X^s(x,y,z):
\begin{cases} 
\dot x = P^1(x,y,z) := ax + by -ax^2+(\tilde{\mu}-\alpha z)x(2-3x), \\
\dot y = P^2(x,y,z) := bx +ay -\frac32 bx^2-\frac32 axy-2y(\tilde{\mu}-\alpha
z), \\
\dot z = P^3(x,y,z) := cz +\mu x +\gamma xz +\alpha \beta  (x^2(1-x)-y^2),
\end{cases}
\end{equation*}
for $(x,y,z) \in \R^3$. We wish to apply a conjugacy  transformation
such that $X^s$ is topologically equivalent to a vector field
$\overline{X}$ when restricted to the open ball $\mB^2(2)\subset\R^3$
of radius $2$.  As the first step, we extend
the system into $\R^4$, that is, we transform system~\cref{eq:san}
such that it is defined on
the unit hypersphere $\mS^3$. Then the two-dimensional sphere 
\begin{equation*}
\mS^2_{\R^4}:=\lp \{ (x_1,x_2,x_3,x_4) \in \mS^3 : x_4=0 \rp \} \subset \R^4
\end{equation*}
is the equator that contains the dynamics at infinity of
system~\cref{eq:san}. Analogous to the one-dimensional case,  we
use the inverse central projections $f_{\pm} : \R^3
\rightarrow \mS_{\pm}^3$ defined by $f_{\pm}(x_1,x_2,x_3)=\pm
(x_1,x_2,x_3,1)/(1+x_1^2+x_2^2+x_3^2)^{1/2}$; note that the radius
$x_1^2+x_2^2+x_3^2$ plays the exact same role as $r^2$ in the
one-dimensional example. We perform a
conjugacy transformation to the vector field $X^s$ on $\R^3$ so that
we obtain the vector fields $X^s_{\mS_{\pm}}=Df_{\pm} \circ X^s \circ
f^{-1}_{\pm}$ defined on the tangent bundle $\text{T}\mS^3_{\pm}$.  It
turns out that both $X^s_{+}$ and $X^s_{-}$ can be expressed as
\begin{equation} \label{eq:comPreSph}
X^s_{\mS_{\pm}} (y)=y_4
\begin{pmatrix}
1-y_1^2 & -y_1y_2 & -y_1y_3 \\[1em]
-y_1y_2  & 1-y_2^2 & -y_2y_3 \\[1em]
-y_1y_3  & -y_2y_3 & 1-y_3^2 \\[1em]
-y_1y_4 & -y_2y_4 & -y_3y_4
\end{pmatrix} \circ
X^s\lp(\dfrac{y_1}{y_4},\dfrac{y_2}{y_4},\dfrac{y_3}{y_4}\rp),
\end{equation}
where $y=(y_1,y_2,y_3,y_4) \in \mS^3_{\pm}$. 

System~\cref{eq:comPreSph} is not well defined on the hyperplane
$y_4=0$. We can salvage this issue via multiplication by a factor 
$y_4^{k-1}$, where $k$ is the maximal degree of the polynomials that
define $X^s$.  Since $k=3$ for system~\cref{eq:san},  we define the
corresponding Poincar\'e compactification on $\mS^3$ as
\begin{equation} \label{eq:Comparse}
\rho(X^s)(y)=y_4^2X^s_{\mS_{\pm}}(y),
\end{equation}  
that is, $\rho(X^s)$ is defined on $\mS^3_{\pm}$ as well as the
equator.  We can think of \cref{eq:Comparse} as a 
vector field on $\R^4$, for which $\mS^3$ is an invariant manifold. Note that,  if $k$ were even,
the dynamics of $\rho(X^s)(y)$ on the  hemisphere $\mS^3_{-}$ would
only be conjugate to $X^s$ by reversing time.  

\subsection{Projection back to $\R^3$} 
As illustrated for the one-dimensional vector field in
\cref{fig:IdCom}(b), we now project $\mS^3 \setminus\{0,0,0,-1\} \subset \R^4$ back to
$\R^3$. We define $g:\mS^3\setminus\{0,0,0,-1\}\rightarrow \R^3$ 
as
\begin{equation*}
g(y_1,y_2,y_3,y_4)=\dfrac{2}{y_4+1}
\lp( y_1, y_2, y_3\rp),
\end{equation*}
which corresponds to the stereographic projection from the south pole
$(0,0,0,-1)$ to the hyperplane tangent to $\mS^3$ at the north pole
$(0,0,0,1)$. The set $g(\mS^3_+)$ is contained in the
three-dimensional sphere $\mS^2(2)\subset\R^3$ with radius two. Its
Jacobian is given by  
\begin{equation*}
Dg(y)=\dfrac{2}{y_4+1}
\begin{pmatrix}
1 & 0 & 0 & \dfrac{-y_1}{y_4+1}\\[1em]
0 & 1& 0 & \dfrac{-y_2}{y_4+1} \\[1em]
0 & 0 & 1 & \dfrac{-y_3}{y_4+1}
\end{pmatrix},
\end{equation*}
and the composition with $Df_{\pm}$ becomes
\begin{equation}\label{eq:Trans2}
Dg \circ Df_{\pm}(y) = 2\dfrac{y_4}{(y_4+1)^2}
\lp ( \begin{array}{ccc}
-y_1^2+y_4+1& -y_1y_2 & -y_1y_3 \\[1em]
-y_1y_2 & -y_2^2+y_4+1 & -y_2y_3 \\[1em]
-y_1y_3 & -y_2y_3 & -y_3^2+y_4+1
\end{array}\rp ).
\end{equation}

Let $\bar{p}=(\bar{x},\bar{y},\bar{z}) \in \R^3$ be a point in the new
compactified phase space.  The inverse of $g$ transforms $\bar{p}$ to
the point
\begin{equation*}
(y_1,y_2,y_3,y_4)=\dfrac{4}{\aNorm{\bar{p}}^2+4} \lp(\bar{x},
\bar{y},\bar{z} ,\dfrac{4-\aNorm{\bar{p}}^2}{4}\rp) \in \mS^3 \setminus\{0,0,0,-1\} ,
\end{equation*}
i.e., $\aNorm{g^{-1}(\bar{p})}=1$. In these coordinates,
\cref{eq:Trans2} becomes
\begin{equation*}
Dg \circ Df (\bar{p}) = \dfrac{4-\aNorm{\bar{p}}^2}{\aNorm{\bar{p}}^2+4}
\lp ( \begin{array}{ccc}
\dfrac{-\bar{x}^2+\bar{y}^2+\bar{z}^2+4}{4} &-\dfrac{\bar{x}\bar{y}}{2} & -\dfrac{\bar{x} \bar{z}}{2} \\[1em]
-\dfrac{\bar{x}\bar{y}}{2} & \dfrac{\bar{x}^2-\bar{y}^2+\bar{z}^2+4}{4} & -\dfrac{\bar{y}\bar{z}}{2} \\[1em]
-\dfrac{\bar{x}\bar{z}}{2} & -\dfrac{\bar{y}\bar{z}}{2} & \dfrac{\bar{x}^2+\bar{y}^2-\bar{z}^2+4}{4} 
\end{array}\rp ).
\end{equation*}
Then the vector field $\overline{X^s}=Dg \circ \rho(X^s) \circ
g^{-1}=y_4^2 \, Dg \circ Df_{\pm} \circ X^s \circ f^{-1}_{\pm}
\circ g^{-1}$, defined on $\R^3$, can be expressed as
\begin{equation} \label{eq:ultVectField}
\overline{X^s}(\bar{p}) = \lp(
\dfrac{4-\aNorm{\bar{p}}^2}{\aNorm{\bar{p}}^2+4}
\rp )^2 Dg \circ Df \circ
X^s\lp( \dfrac{4\bar{x}}{4-\aNorm{\bar{p}}^2},
\dfrac{4\bar{y}}{4-\aNorm{\bar{p}}^2},
\dfrac{4\bar{z}}{4-\aNorm{\bar{p}}^2} \rp).
\end{equation}
Note that the solid sphere of radius two
is invariant under system~\cref{eq:ultVectField} and its interior is
conjugate to the original system~\cref{eq:san} on $\R^3$.

In particular, system~\cref{eq:ultVectField} can be used to track the different
equilibria of system~\cref{eq:san} in \textsc{Auto} as they move
through infinity, to compute the bifurcation diagrams in \cref{fig:BifInf}. 

\subsection{Analytical study of infinity} \label{sec:AnaStInf}
In its general form, system~\cref{eq:ultVectField} is too complex to
study the dynamics at infinity, that is, on the boundary $\mS^2(2)$ of
the compactified phase space.  Instead, we study the dynamics at infinity for the differentiable vector field
$\rho(X^s)$ as given by   \cref{eq:comPreSph}. To this end, we analyze the
dynamics on specific coordinate charts of $\mS^3$ and its equator $\mS^2_{\R^4}$
\cite{Jes1, Jes2}.  We consider three different local charts, namely, $(U_i,
\phi_i)$ for $i=1,2,3$, where $U_i= \{y \in \mS^3 : y_i > 0 \}$ and 
$\phi_i: U_i \rightarrow \R^3$; the transformations $\phi_i$
correspond to the central projections with respect to the
tangent planes at the points $(1,0,0,0)$, $(0,1,0,0)$ and $(0,0,1,0)$, 
respectively, which are similar to the projections $f^{-1}_{\pm}$
used in \cref{sec:compofR3}. The three-dimensional vector fields in these projections
contain subsets of the equator $\mS^2_{\R^4}$ that correspond to invariant
planes. Hence, the problem of studying the dynamics at infinity
 can be simplified to a study of two-dimensional vector fields \cite{Jau1,Jes1, Jes2,Raz1}.  For
 ease of notation, we use the variables $\tilde{x}, \tilde{y}, 
\tilde{z}$ and $\tilde{w}$ interchangeably in the different charts.
Specifically, $\tilde{w}$ represents the proximity to infinity,
that is,  $\tilde{w}=0$ corresponds to the projection of the dynamics infinity in
the corresponding chart.

We show the construction for $U_1$, that is, the half of $\mS^3$ with
$y_1>0$.  Similar to $f_{\pm}$, the inverse central projection with the
hyperplane tangent to $\mS^3$ at $(1,0,0,0)$ adds $1$ as the first
component and normalizes the vector.  Hence, its inverse $\phi_1(y)$
for $y \in U_1$ is defined as $\phi_1(y)=(y_2/y_1,y_3/y_1,y_4/y_1)=:(\tilde{y},\tilde{z},\tilde{w})
\in \R^3$, and the corresponding Jacobian is given by
\begin{equation*}
D\phi_1(y) = \dfrac{1}{y_1} 
\begin{pmatrix}
-\dfrac{y_2}{y_1} & 1 & 0 &0 \\[1em]
-\dfrac{y_3}{y_1} & 0 & 1 &0 \\[1em]
-\dfrac{y_4}{y_1} & 0 & 0 &1
\end{pmatrix}.
\end{equation*}
then the composition with $Df_{\pm}$ becomes
\begin{equation}\label{eq:ProPhi1}
D\phi_1 \circ Df (y) = \dfrac{1}{y_1}
\begin{pmatrix}
-\dfrac{y_2}{y_1} & 1 & 0 \\[1em]
-\dfrac{y_3}{y_1} & 0 & 1 \\[1em]
-\dfrac{y_4}{y_1} & 0  & 0
\end{pmatrix}.
\end{equation}
Rewriting \cref{eq:ProPhi1} with respect 
$(\tilde{y},\tilde{z},\tilde{w})$, we have
\begin{equation}\label{eq:ProPhi2}
D\phi_1 \circ Df (\tilde{y},\tilde{z},\tilde{w}) = \tilde{w}
\begin{pmatrix}
-\tilde{y} & 1 & 0\\
-\tilde{z}  & 0 & 1 \\
-\tilde{w} & 0 & 0
\end{pmatrix}.
\end{equation}

Finally, we use the fact that
$y_4=\tilde{w}/(1+\tilde{y}^2+\tilde{z}^2+\tilde{w}^2)^{1/2}$ and
\cref{eq:ProPhi2} to represent the vector field on $U_1$, that is,
$X^s_{U_1} = D\phi_1 \circ \rho(X^s) \circ \phi_1^{-1}$, as
\begin{equation}\label{eq:Phi1}
X^s_{U_1}(\tilde{y},\tilde{z},\tilde{w})= 
\dfrac{\tilde{w}^2}{1+\tilde{y}^2+\tilde{z}^2+\tilde{w}^2} D\phi_1
\circ Df \circ X^s=\dfrac{\tilde{w}^3}{1+\tilde{y}^2+\tilde{z}^2+\tilde{w}^2}
\begin{pmatrix}
-\tilde{y} P^1+P^2\\
-\tilde{z} P^1+P^3\\
-\tilde{w} P^1
\end{pmatrix},
\end{equation}
where $P^j=P^j(1/\tilde{w},\tilde{y}/\tilde{w},\tilde{z}/\tilde{w})$
for $j=1,2,3$. The dynamics on the chart $U_2$ with
$\phi_2(y)=(y_1/y_2,y_3/y_2,y_4/y_2)=:(\tilde{x},\tilde{z}, \tilde{w}) \in
\R^3$, where $y \in U_2$,  are given by
\begin{equation}\label{eq:Phi2}
X^s_{U_2}(\tilde{x},\tilde{z},\tilde{w})= 
\dfrac{\tilde{w}^2}{1+\tilde{x}^2+\tilde{z}^2+\tilde{w}^2} D\phi_2
\circ Df \circ X^s= \dfrac{\tilde{w}^3}{1+\tilde{x}^2+\tilde{z}^2+\tilde{w}^2}
\begin{pmatrix}
-\tilde{x} P^2+P^1\\
-\tilde{z} P^2+P^3\\
-\tilde{w} P^2
\end{pmatrix},
\end{equation}
where  $P^j=P^j(\tilde{x}/\tilde{w}, 1/\tilde{w},\tilde{z}/\tilde{w})$ for
$j=1,2,3$. Finally, the dynamics on the chart $U_3$ with
$\phi_3(y)=(y_1/y_2,y_3/y_2,y_4/y_2)=:(\tilde{x},\tilde{z},\tilde{w})
\in \R^3$, and $y \in U_3$, are given by
\begin{equation}\label{eq:Phi3}
X^s_{U_3}(\tilde{x},\tilde{y},\tilde{w})= 
\dfrac{\tilde{w}^2}{1+\tilde{x}^2+\tilde{y}^2+\tilde{w}^2} D\phi_3
\circ Df \circ X^s= \dfrac{\tilde{w}^3}{1+\tilde{x}^2+\tilde{y}^2+\tilde{w}^2}
\begin{pmatrix}
-\tilde{x} P^3+P^1\\
-\tilde{y} P^3+P^2\\
-\tilde{w} P^3
\end{pmatrix},
\end{equation}
where  $P^j=P^j(\tilde{x}/\tilde{w}, \tilde{y}/\tilde{w},1/\tilde{w})$ for
$j=1,2,3$.

Note that the denominator term in each of the factors for
\cref{eq:Phi1}, \cref{eq:Phi2} and \cref{eq:Phi3} is strictly
positive.  Hence, this term can be viewed as a time rescaling that does not alter
the dynamics of the vector fields; therefore, it can be omitted. As mentioned
before, $\tilde{w}=0$ is an invariant plane 
for \cref{eq:Phi1}, \cref{eq:Phi2} and \cref{eq:Phi3} that
represents infinity. After substitution of the corresponding polynomials and
simplification of the expressions, we set 
$\tilde{w}=0$ in  \cref{eq:Phi1}, \cref{eq:Phi2} and
\cref{eq:Phi3}, which leads to the following three vector fields that
represent the dynamics of system~\cref{eq:san} at infinity in the
corresponding charts:
\begin{align}
X^s_{U_1^{\infty}}(\tilde{y} ,\tilde{z})&:
\begin{cases} 
\dot{\tilde{y}} = -3 \alpha\tilde{y}\tilde{z}, \\
\dot{\tilde{z}} = - \alpha (3 \tilde{z}^2+ \beta).
\end{cases} \label{eq:UUU1} \\
X^s_{U_2^{\infty}}(\tilde{x} ,\tilde{z})&:
\begin{cases} 
\dot{\tilde{x}} = 3 \alpha \tilde{x}^2\tilde{z}, \\
\dot{\tilde{z}} = - \alpha \beta \tilde{x}^3.
\end{cases}\label{eq:UUU2} \\
X^s_{U_3^{\infty}}(\tilde{x} ,\tilde{y})&:
\begin{cases} 
\dot{\tilde{x}} = \alpha \tilde{x}^2(\beta \tilde{x}^2+3), \\
\dot{\tilde{y}} = \alpha \beta \tilde{x}^3 \tilde{y}.
\end{cases}\label{eq:UUU3}
\end{align}
Systems~\cref{eq:UUU1},  \cref{eq:UUU2} and \cref{eq:UUU3}
highlight that the dynamics at infinity only 
depends on the parameters $\alpha$ and $\beta$, which are
the coefficients of higher powers of the polynomials in system~\cref{eq:san}.
We observe that the three systems each have an integral of motion, namely,
\begin{align}
H_{U_1}(\tilde{y},\tilde{z}) &= \ln \lp(
\frac{\tilde{y}^2}{| 3\tilde{z}^2+\beta |} \rp), \label{eq:H1L}\\
H_{U_2}(\tilde{x},\tilde{z}) &= \beta \tilde{x}^2 + 3
\tilde{z}^2 \qquad \mbox{ and } \label{eq:H2L}\\
H_{U_3}(\tilde{x},\tilde{y}) &= \ln \lp(
\frac{|\beta\tilde{x}^2+3|}{\tilde{y}^2} \rp). \label{eq:H3L}
\end{align}
Given the parameter chosen in \cref{sec:san}, we are
interested in how the dynamics at infinity changes as $\beta$ is varied and
$\alpha > 0$.

Note that, for a complete characterization of $\mS^2_{\R^4}$, one would also
have to study the charts  $(V_i, \sigma_i)$ with $i=1,2,3$, where
$V_i= \{y \in \mS^3 : y_i < 0 \}$ and 
$\sigma_i: V_i \rightarrow \R^3$ are the central projections to the
tangent planes $(-1,0,0,0)$, $(0,-1,0,0)$ and $(0,0,-1,0)$. We do not
study these charts, because the inverse central projections $f_{\pm}$ 
map to antipodal points and the maximum degree of our polynomials is
odd; therefore, the charts $U_i$ and $V_i$ are 
conjugate to each other via the transformation $p \in U_i \rightarrow
-p \in V_i$, for $i=1,2,3$.
\subsubsection{Dynamics at infinity when $\beta < 0$}
\begin{figure}
\centering
\includegraphics{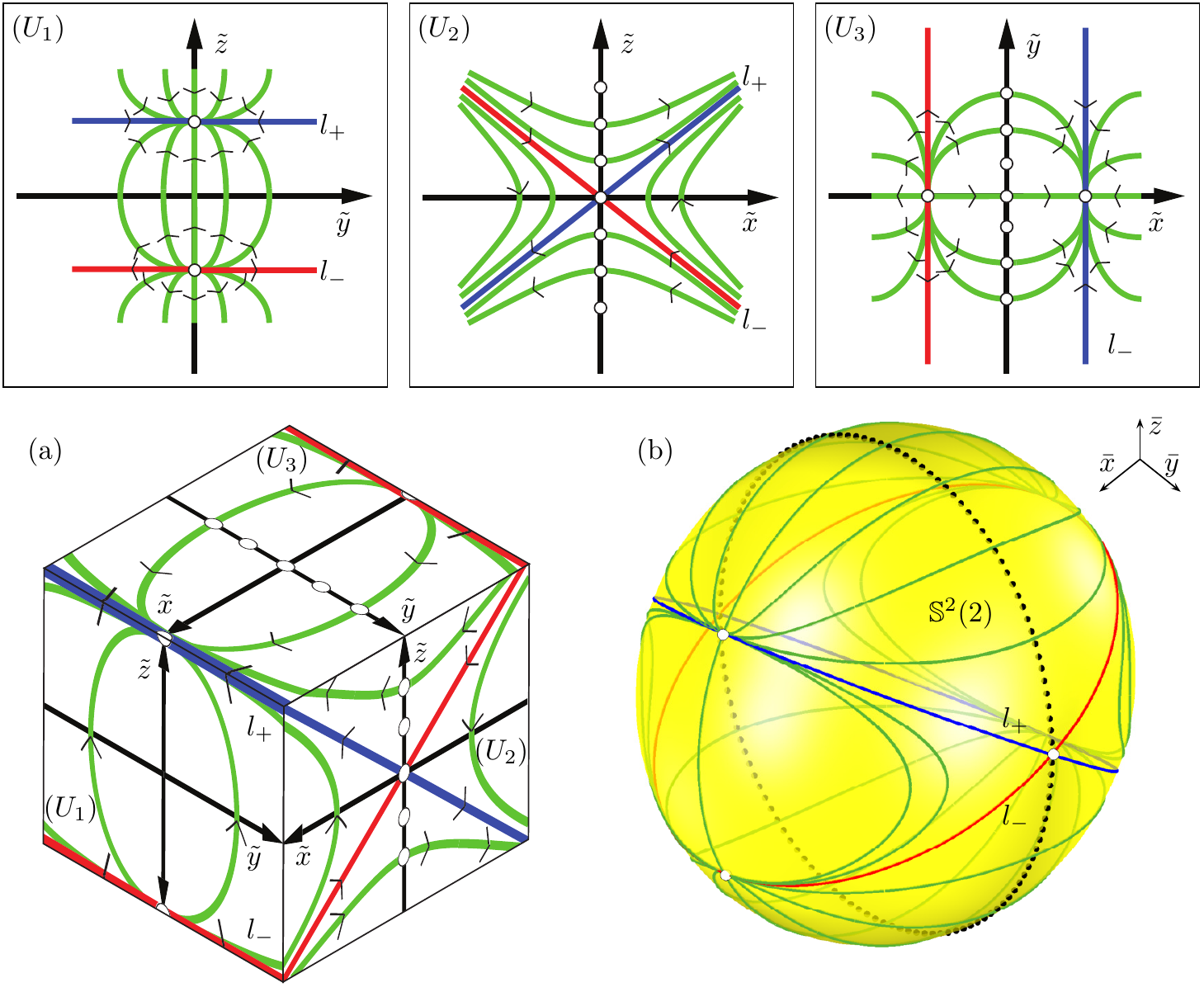}
\caption{Dynamics of system~\cref{eq:san}  
  at infinity when $\beta<0$; the first row shows 
  different coordinate charts that fit
  together as neighboring faces on a cube in panel (a). The
  sphere in panel (b) shows computed trajectories (green) on
  $\mS^2(2)$ for the compactified system~\cref{eq:ultVectField}  for
   $(a,b,c,\alpha,\beta,\gamma,\mu,\tilde{\mu})=(0.22,1,-2,0.65,-1,2,0,0)$.}\label{fig:SteoBL} 
\end{figure}
\cref{fig:SteoBL} illustrates the dynamics of system~\cref{eq:san} at
infinity when $\beta<0$.  The first row shows phase portraits of    
systems~\cref{eq:UUU1}, \cref{eq:UUU2} and \cref{eq:UUU3} on the
charts correspondingly labeled panels~$(U_1)$, $(U_2)$ and $(U_3)$,
respectively. Panel (a)  shows how these 
phase portraits are glued together on a cube and panel(b)
representative computed trajectories of system~\cref{eq:ultVectField}
on $\mS^2(2)$ for
$(a,b,c,\alpha,\beta,\gamma,\mu,\tilde{\mu})=(0.22,1,-2,0.65,-1,2,0,0)$. 

In the local chart $U_1$, there are two equilibria, at
$(\tilde{y},\tilde{z})=(0,\pm\sqrt{-\beta/3})$, and their Jacobian matrix is
diagonal with eigenvalues $\mp 3\alpha \sqrt{-\beta/3})$ and
$\mp 6\alpha \sqrt{-\beta/3}$, respectively; the eigenvectors are $(1,0)$ for the
first and $(0,1)$ for the second eigenvalue, which is the strong
direction. Note that $\tilde{y}=0$ is invariant, so that the strong
(un)stable manifolds are straight lines that coincide. Since the lines $\tilde{z}=\pm
\sqrt{-\beta/3}$ are also invariant, both equilibria have linear weak
(un)stable manifolds as well; we denote these straight lines by
$l_{\pm}$.  The curves $l_{\pm}$ connect the points $(0,2,0)$ and 
$(0,-2,0)$ in the phase space of \cref{eq:ultVectField}.

In the local chart $U_2$, as we see from equation \cref{eq:H2L},
solutions of system~\cref{eq:UUU2} are tangent to the family of 
hyperbolas with asymptotes $\tilde{z}=\pm
\tilde{x} \sqrt{-\beta / 3}$ when $\beta < 0$.  The curves $l_{\pm}$
correspond to these asymptotes in the local chart $U_2$. They are
the only trajectories that converge to the origin, which is the point
$(0,2,0)$ of $\mS^2(2)$.  The $\tilde{z}$-axis is a family of non-hyperbolic
equilibria of system~\cref{eq:UUU2} with one stable direction when
$\tilde{x}<0$ and one unstable direction when $\tilde{x}>0$.   

The phase portrait of system~\cref{eq:UUU3} in the chart $U_3$  with $\beta<0$ is similar to that of
system~\cref{eq:UUU1} after the transformation $(\tilde{y},\tilde{z})
\mapsto (-\tilde{y},\tilde{x})$ and rotation by $-\frac{\pi}{2}$;
panel $(U3)$ of \cref{fig:SteoBL} shows that the $\tilde{x}$-axis is invariant and corresponds to the strong
manifolds of the two equilibria $(\pm\sqrt{-3/\beta},0)$; also, the vertical lines $\tilde{x}=\pm
\sqrt{-3/\beta}$ are invariant and correspond to the
projections of $l_{\pm}$ under $\phi_3$. However, due to the factor
$\tilde{x}^2$ in the equations, the $\tilde{z}$-axis is a set of
non-hyperbolic equilibria just as for system~\cref{eq:UUU2}.

The cube and sphere in the second row of \cref{fig:SteoBL} show
how the curves $l_{\pm}$ connect the two nodes in the chart $U_1$ with the
equilibrium $(0,0)$ in the chart $U_3$. 

\subsubsection{Dynamics at infinity when $\beta = 0$}
\begin{figure}
\centering
\includegraphics{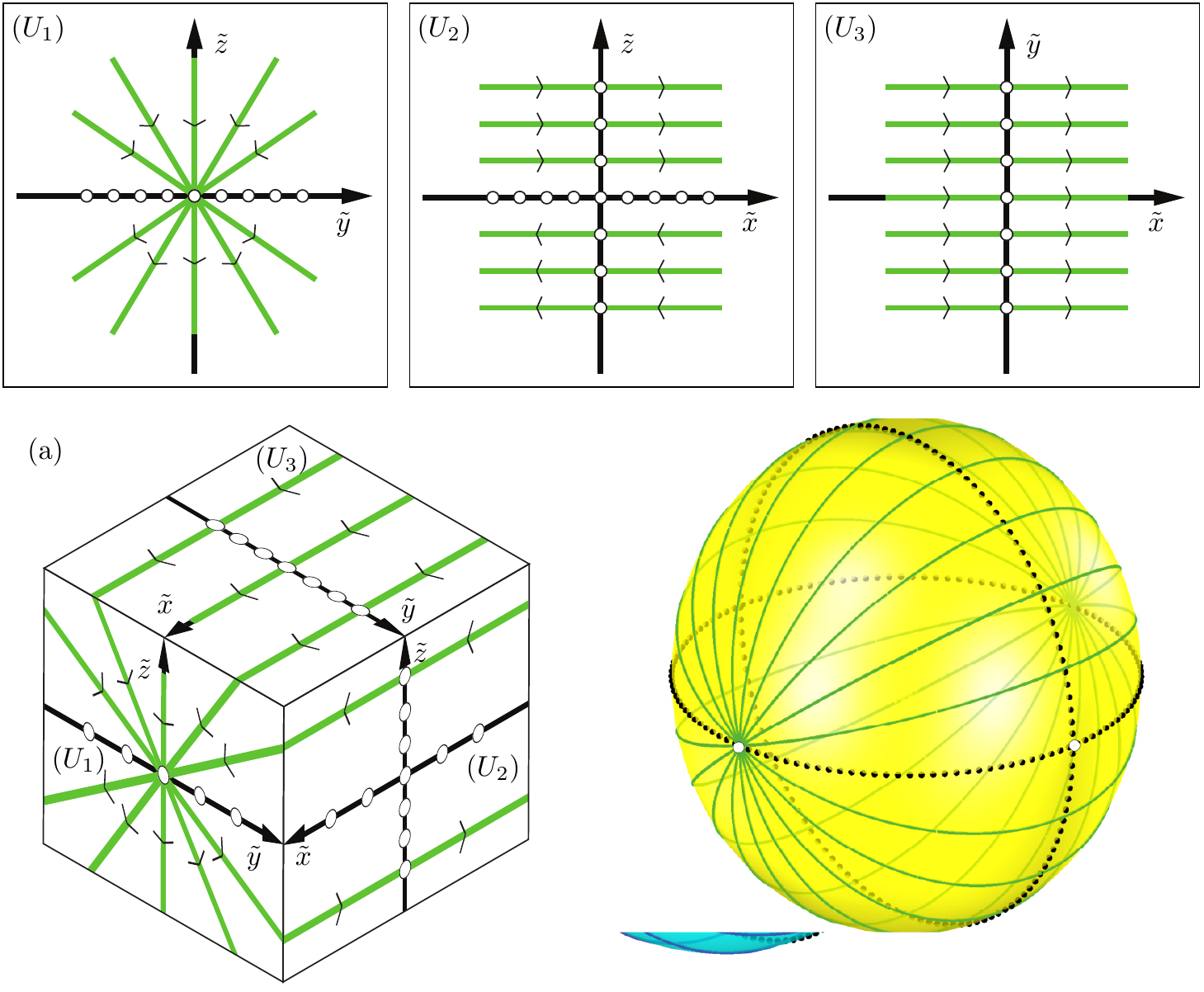}
\caption{Dynamics of system~\cref{eq:san}  
  at infinity when $\beta=0$; the first row shows 
  different coordinate charts that fit
  together as neighboring faces on a cube in panel (a). The
  sphere in panel (b) shows computed trajectories (green) on
  $\mS^2(2)$ for the compactified system~\cref{eq:ultVectField}
   for  $(a,b,c,\alpha,\beta,\gamma,\mu,\tilde{\mu})=(0.22,1,-2,0.65,0,2,0,0)$.} \label{fig:SteoB0} 
\end{figure}
Systems \cref{eq:UUU1}-\cref{eq:UUU3} for $\beta=0$ become
\begin{equation*}
X^s_{U_1^{\infty}}(\tilde{y} ,\tilde{z}):
\begin{cases} 
\dot{\tilde{y}} = -3 \alpha\tilde{y}\tilde{z} , \\
\dot{\tilde{z}} = - 3\alpha \tilde{z}^2,
\end{cases} \quad
X^s_{U_2^{\infty}}(\tilde{x} ,\tilde{z}):
\begin{cases} 
\dot{\tilde{x}} = 3 \alpha \tilde{x}^2\tilde{z}, \\
\dot{\tilde{z}} = 0,
\end{cases} \quad
X^s_{U_3^{\infty}}(\tilde{x} ,\tilde{y}):
\begin{cases} 
\dot{\tilde{x}} = 3 \alpha \tilde{x}^2, \\
\dot{\tilde{y}} = 0.
\end{cases}
\end{equation*}

The corresponding phase portraits are shown in the first row of
\cref{fig:SteoB0}; in the different charts $U_1$, $U_2$ and $U_3$.  Panel $(U1)$ of
\cref{fig:SteoB0} shows a sketch of the phase portrait \cref{eq:UUU1}
for $\beta=0$. Note that the $\tilde{y}$-axis consists of non-hyperbolic
equilibria and the $\tilde{z}$-axis is invariant. In fact,
its first integral of motion \cref{eq:H1L} can be simplified to
\begin{equation*}
H_{U_1}(\tilde{y},\tilde{z}) = \frac{\tilde{y}}{\tilde{z}},
\end{equation*} 
which means that any straight line through the
origin is invariant. We can think of the origin as a saddle-node 
equilibrium; indeed, as $\beta<0$ increases towards  $0$, the two
equilibria from \cref{fig:SteoBL} move closer together and eventually,
at $\beta=0$, merge at the origin. Similarly, the lines $l_{\pm}$ meet
at the $\tilde{y}$-axis for $\beta=0$ and become a family of 
non-hyperbolic equilibria. 

Panel $(U2)$ of \cref{fig:SteoB0} shows a sketch of the phase portrait
\cref{eq:UUU2} for $\beta=0$ and illustrates that the dynamics in the
local chart $U_2$ are reduced to one-dimensional 
dynamics.  The $\tilde{x}$- and $\tilde{z}$-axes are
 families of degenerate equilibria, and the horizontal
lines are invariant.  The non-hyperbolic equilibria on the 
$\tilde{z}$-axis are degenerate saddle-node points.  

The dynamics in system~\cref{eq:UUU3} on $U_3$
are also reduced to one-dimensional
dynamics. The $\tilde{y}$-axis is a set of degenerate saddle-node
equilibria and the horizontal lines are invariant.  In contrast to
the limit argument used for the chart $U_1$,
the two hyperbolic equilibria that exist for $\beta < 0$ do not
disappear in a saddle-node bifurcation at $\beta=0$; instead they
go to infinity in $U_3$.  More precisely, these equilibria disappear at a saddle-node
bifurcation on the local chart $U_1$ and its antipodal chart $V_1$,
which both cannot be seen in the local chart 
$U_3$. 

\cref{fig:SteoB0}(a) shows how the different charts can be glued
together to form a cube. As shown in \cref{fig:SteoB0}(b), we also computed trajectories of system
\cref{eq:ultVectField} on its invariant sphere $\mS^2(2)$;  here,
$(a,b,c,\alpha,\beta,\gamma,\mu,\tilde{\mu})=(0.22,1,-2,0.65,0,2,0,0)$. Note
that the horizontal trajectories in the charts 
$U_2$ and $U_3$ are translated to curves on $\mS^2(2)$ that connect the poles at $(2,0,0)$
and $(-2,0,0)$.  The great circles $\bar{x}=0$ and $\bar{z}=0$ both consist of degenerate
equilibria.
 
\subsubsection{Dynamics at infinity when $\beta > 0$}
\begin{figure}
\centering
\includegraphics{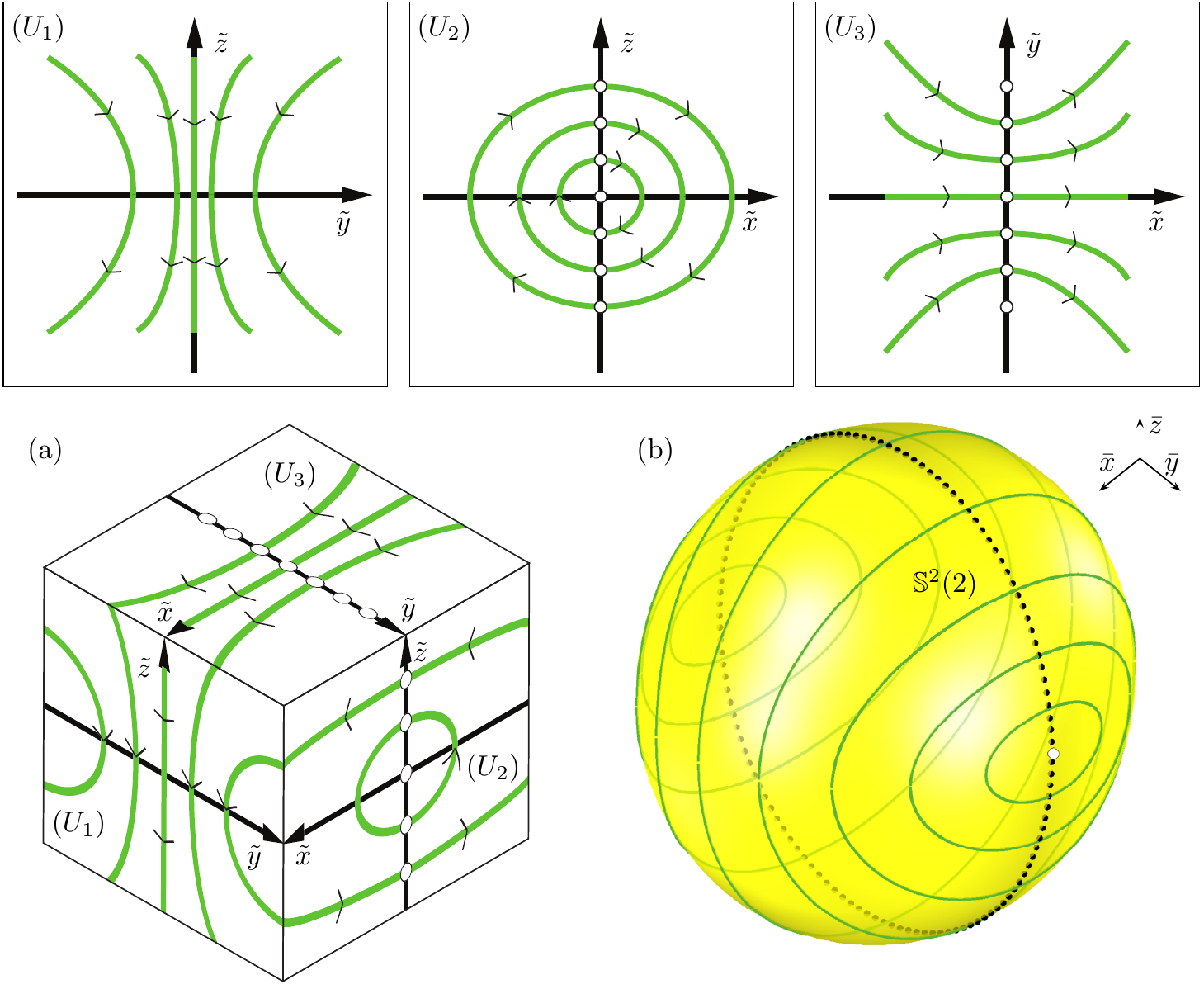}
\caption{ Dynamics of system~\cref{eq:san}  
  at infinity when $\beta>0$; the first row shows 
  different coordinate charts that fit
  together as neighboring faces on a cube in panel (a). The
  sphere in panel (b) shows computed trajectories (green)  on
  $\mS^2(2)$ for the compactified system~\cref{eq:ultVectField}  for
   $(a,b,c,\alpha,\beta,\gamma,\mu,\tilde{\mu})=(0.22,1,-2,0.65,1,2,0,0)$.}\label{fig:SteoBP} 
\end{figure}

\cref{fig:SteoBP} shows the corresponding phase portraits of system
\cref{eq:san} with $\beta > 0$ on the charts at infinity in the
correspondingly labeled panels~$(U_1)$, $(U_2)$ and $(U_3)$,
respectively. On $U_1$  there are no
equilibria when $\beta>0$ and the
$\tilde{z}$-axis is invariant under the flow. The first integral of motion
\cref{eq:H1L} defines the family of radical functions
$\tilde{y} = \pm \sqrt{e^c(3\tilde{z}^2+\beta)}$, for $c \in \R$.
While, on the chart $U_2$, the first integral of motion
\cref{eq:H2L} for system~\cref{eq:UUU2} defines a 
family of ellipses that are oriented clockwise and the $\tilde{z}$-axis
is a family of non-hyperbolic equilibria. Finally, on the chart $U_3$,
system~\cref{eq:UUU3} has a family of non-hyperbolic equilibria on  the
$\tilde{y}$-axis.  As in the local chart $U_1$, its first integral of
motion  \cref{eq:H3L} defines  the family of radical functions 
$\tilde{y} = \pm \sqrt{e^c(\beta\tilde{x}^2+3)}$, for $c \in \R$;
where the $\tilde{x}$-coordinate increases with time. 

Panel (a) of \cref{fig:SteoBP} shows how these projections
fit together on a cube for $\beta>0$. Panel(b) shows computed trajectories of system
\cref{eq:ultVectField} on $\mS^2(2)$.  In particular, note how each
trajectory in $U_1$ is translated to a curve that
connects equilibria  $(0,\bar{y},\bar{z})$ and
$(0,\bar{y},-\bar{z})$. The local charts $U_2$ and $U_3$ imply that
the $\bar{z}$-coordinate of the trajectories with $\bar{x}>0$
decreases with time, while it increases for trajectories with $\bar{x}<0$; this agrees with the
orientation computed for the trajectories on $\mS^2(2)$.

\section{Boundary value problem formulations}\label{secApp:BVP}
More often than not, it is impossible to compute Floquet multipliers of a given saddle periodic
orbit explicitly, let alone approximate the corresponding eigenbundles and 
global stable and unstable manifolds. We employ
continuation of a suitable two-point
boundary value problems (2PBVP) with the software package
\textsc{Auto} \cite{Doe1,Doe2} to solve these problems for a
three-dimensional system of the form \cref{eq:genP}.  The idea behind 2PBVP continuation is to
represent the object of interest as a one-parameter family of
finite-time orbit segments of system~\cref{eq:genP} that satisfy
suitable  boundary conditions; see \cite{Doe4},  for general
background of this  approach.

As discussed in \cite{Agu1,Kra2}
any trajectory of \cref{eq:genP} over the finite-time
interval $[0,T]$ can be represented as an orbit segment $u:[0,1]
\rightarrow \R^n$ over the interval $[0,1]$ that
satisfies equation
\begin{equation} \label{eq:genA}
\dot u = T f(u,\mu),
\end{equation}
which is a time-rescaled version of system~\cref{eq:genP} with
(original) integration time $T>0$.  During the continuation, we impose
additional boundary conditions at $u(0)$ and $u(1)$.

We refer to \cite{Eng1} for an in-depth discussion of the 2PBVP
formulation needed to calculate the Floquet multipliers of a saddle
periodic orbit and their respective tangent bundles.  The 2PBVP formulation and 
computation of two-dimensional stable and unstable manifolds of saddle
periodic orbits (as well as equilibria) and their intersection sets
with a sphere can be found in \cite{Call1,Kra2}. 

In the following sections, we present the 2PBVP formulation for the
following:
\begin{enumerate}
\item The intersection set of the stable manifold of a saddle periodic orbit with a
tubular neighborhood; this is used for the computation of the local
stable and unstable manifolds of the saddle periodic orbits in
\cref{fig:per}.
\item The computation and continuation of a periodic orbit at the
  moment when two Floquet multipliers change from being real to
  complex conjugate;  this allows us to find the  curves $\mathbf{CC^+}$ and $\mathbf{CC^-}$ in
  \cref{fig:BDInc} and \cref{fig:BDO1}.  
\end{enumerate}

\subsection{The intersection set of  the 
  manifold of a saddle periodic orbit with a  tubular section} The
orientation of a two-dimensional stable manifold 
of a saddle periodic orbit can be illustrated by computing a first local portion.  The
approach described here finds this portion as the manifold computed up
to its first intersection with a tubular section of small radius $d
\geq 0$ around the periodic orbit; see rows $1$ and $2$ of \cref{fig:per}.
Our formulation also works particularly well if the Floquet multiplier
associated with the manifold is close to $0$ in magnitude. 

Let $\Gamma$ be a saddle periodic orbit in $\R^3$ and assume that we
wish to compute a first portion of its two-dimensional stable manifold
$W^s(\Gamma)$. We extend the system from three equations to six, so
that we effectively consider two different orbit segments of \cref{eq:genP}:
\begin{equation} \label{eq:BVPstart}
\begin{array}{cc}
\begin{cases} 
\dot{v}_{\Gamma}=T_{\Gamma} f(v_{\Gamma}(t),\mu), \\
\dot{u} =T f(u(t),\mu),
\end{cases}
&
\begin{array}{rcl}
v_{\Gamma}& \in & \Gamma,\\
u& \in & \R^3.
\end{array}
\end{array} 
\end{equation}

The segment $v_{\Gamma}(t) $ is meant to represent $\Gamma$. Hence,
$T_\Gamma$ is the period and we impose the boundary condition: 
\begin{equation}\label{eq:BVPperE}
v_{\Gamma}(1)-v_{\Gamma}(0)=0.
\end{equation}
The idea is that $u(t)$ represents a solution
trajectory with integration time $T$ that is contained in
$W^s(\Gamma)$.  Since $u(t)$ converges to $\Gamma$ as $t$ goes to
infinity, we stipulate that  $u(1)$ lies close to $\Gamma$, in an
approximate one-dimensional fundamental domain 
$\mathcal{F}_{\delta}$ of the linear approximation of $W^s(\Gamma)$.
Every (approximated) trajectory in $W^s(\Gamma)$ 
intersects $\mathcal{F}_{\delta}$  exactly once, and this domain is
parametrized by the variation of $\delta$ in a closed interval; see
\cite{Call1} for details. The parameterized boundary condition:
\begin{equation}\label{eq:BVPfunD}
\begin{array}{cc}
u(1) \in \mathcal{F}_{\delta}
\end{array}
\end{equation}
introduces a free parameter $\delta$ on top of the (free)
parameters $T_\Gamma$ and $T$.

The tubular section with radius $d$ around $\Gamma$ is defined as 
\begin{equation*}
{\rm T}^d_{\Gamma} := \lp\{x \in \R^3 :  \min_{y \in \Gamma} \aNorm{x-y}_{\R^3} =d \rp\},
\end{equation*}
where $\aNorm{\cdot}_{\R^3} $ is the Euclidean norm in $\R^3$.
We are interested in orbit segments $u$ of \cref{eq:BVPstart} that satisfy
\cref{eq:BVPperE} and \cref{eq:BVPfunD}, and also $u(0) \in {\rm T}^d_{\Gamma}$.  The
family of all such orbit segments forms the first portion of $W^s(\Gamma)$ and its end
points $u(\cdot)$ form the intersection set $W^s(\Gamma) \cap {\rm
  T}^d_{\Gamma}$, which is a one-dimensional curve. We parametrize
this set on ${\rm T}^d_{\Gamma}$ via the 
points on $\Gamma$ that achieve the minimum $d$ for the point $u(0)$,
that is, we use the orbit segment $v_\Gamma$ of
system~\cref{eq:BVPstart} to track $W^s(\Gamma)\cap {\rm
  T}^d_{\Gamma}$. We impose the following two boundary conditions:
\begin{align}
\langle f(v_{\Gamma}(0),\mu) \,,\, u(0)- v_{\Gamma}(0) \rangle & = \alpha, \label{eq:BVPtanP}\\
|| v_{\Gamma}(0)- u(0) ||_{\R^3} & = d, \label{eq:BVPend}
\end{align} 
where $\langle\,,\rangle$ is the dot product.  When $\alpha=0$,
condition \cref{eq:BVPtanP} implies that $u(0)$ lies in the plane
normal to $\Gamma$ at $v_{\Gamma}(0)$; this is a necessary condition
for $v_\Gamma(0)$ to achieve the minimal distance of $u(0)$ to $\Gamma$. Condition 
\cref{eq:BVPend} defines the radius of the tubular section, that is, it
ensures that $u(0)$ lies on ${\rm T}^d_{\Gamma}$.  The sequence of
steps to follow in \textsc{Auto} is:
\begin{enumerate}
\item Pre-compute from another run a periodic solution $v_\Gamma(t)$
  with period $T_\Gamma$ and its respective fundamental domain
  $\mathcal{F}_{\delta}$. 
\item Extend the system with the 2PBVP-formulation  \cref{eq:BVPstart}
  to \cref{eq:BVPend}, and define a first solution $u=u(0)=u(1) \in
  \mathcal{F}_{\delta}$ with $T=0$. 
\item Continue in $\alpha$ and let $d$ vary
  until $\alpha=0$.  This step rotates
  $v_\Gamma(0)$ along $\Gamma$ until $\aNorm{u(0)-v_\Gamma(0)}_{\R^3}$
  is minimal. Here, $T=0$ and $d$ are fixed, and $T_\Gamma$ is a
  continuation parameter.
\item Fix $\alpha=0$ and continue in $d$ until a suitable distance is
  reached. Here, both $T$ and $T_\Gamma$ are free, but $T_\Gamma$ will
  remain almost contant and $T$ increases. 
\item Fix $d$ and continue in $\delta$ while $T$ and $T_\Gamma$
  vary. The $\delta$-family of orbit segments computed in this run forms
  $W^s(\Gamma)$ with  the local part of $u(0) \in {\rm T}^d_{\Gamma} $.
\end{enumerate}
If $\Gamma$ is orientable, $W^s(\Gamma) \cap {\rm T}^d_{\Gamma}$
consists of two closed curves; if $\Gamma$ is non-orientable, on the
other hand, $W^s(\Gamma) \cap {\rm T}^d_{\Gamma}$ is a single closed
curve that is found in one continuation run during which
$v_{\Gamma}(0)$ rotates along $\Gamma$ twice.

\subsection{BVP formulation for the computation of the curves $\mathbf{CC^+}$
  and $\mathbf{CC^-}$} 
The strong stable manifold of an attracting periodic orbit $\Gamma^a$
disappears when its Floquet multipliers change from being real to
complex conjugate.  The curve $\mathbf{CC^+}$ and $\mathbf{CC^-}$
represent the moment that two real positive or negative Floquet
multipliers become complex conjugate, respectively.   We use the 2PBVP formulation presented in
\cite{Eng1} to compute the Floquet multipliers and their respective
bundles, and follow these steps in \textsc{Auto}:
\begin{enumerate}
\item Continue the periodic orbit $\Gamma^a$ with one of its Floquet
  multipliers and associated eigenbundle in a system parameter
  $\mu_1$. Here, the period $T_\Gamma^a$ of $\Gamma^a$ varies.
\item The moment when the Floquet multiplier becomes complex is
  detected in  \textsc{Auto} as a fold point. One has to be careful,
  because \textsc{Auto}  also marks an actual saddle-node bifurcation
  of periodic orbits as a fold point, which occurs when the Floquet multiplier is $1$.
\item Compute the locus of the fold point by continuing in $\mu_1$
  and a second system parameter $\mu_2$. Here, the period $T_\Gamma^a$ and
  the value of the Floquet multiplier are free parameters.
\end{enumerate}
The set of points $\mu_1$ and $\mu_2$ computed in step 3 represents
  the curve $\mathbf{CC^{+}}$ or $\mathbf{CC^{-}}$. 

\section*{Acknowledgments}
The authors thank Pablo Aguirre for helpful discussion on
homoclinic flip bifurcations and the computation of manifolds in \textsc{Auto}.
\bibliographystyle{siamplain} 

\end{document}